\documentclass[11pt, twoside, a4paper]{book}

\usepackage[all]{xy}
\usepackage[T1]{fontenc}
\usepackage[english]{babel}
\usepackage[utf8]{inputenc}
\usepackage{setspace}                   
\usepackage{makeidx}                    
\usepackage[nottoc]{tocbibind}          
\usepackage{lmodern}
\usepackage{type1cm}                    
\usepackage{titletoc}
\usepackage[bf,small,compact]{titlesec} 
\usepackage[fixlanguage]{babelbib}
\usepackage[font=small,format=plain,labelfont=bf,up,textfont=it,up]{caption}
\usepackage[usenames,svgnames,dvipsnames]{xcolor}
\usepackage[a4paper,top=2.54cm,bottom=2.5cm,left=2.5cm,right=2.54cm]{geometry} 
\usepackage[pdftex,plainpages=false,pdfpagelabels,pagebackref,colorlinks=true,citecolor=DarkGreen,linkcolor=NavyBlue,urlcolor=DarkRed,filecolor=green,bookmarksopen=true]{hyperref} 
\usepackage[all]{hypcap}                
\usepackage[square,sort,nonamebreak,comma]{natbib}  
\fontsize{60}{62}\usefont{OT1}{cmr}{m}{n}{\selectfont}

\usepackage{nameref}

\usepackage[dvipsnames]{xcolor}
\usepackage{xargs}
\usepackage[colorinlistoftodos,prependcaption,textsize=normalsize]{todonotes}

\usepackage{amsmath}
\usepackage{amssymb}
\usepackage{amsthm}
\usepackage{eucal}
\usepackage{mathrsfs}
\usepackage{skull}
\usepackage{mathtools}
\usepackage{marvosym} 

\usepackage{makeidx}

\usepackage[nottoc]{tocbibind}   
\usepackage[fixlanguage]{babelbib}

\usepackage{graphicx}
\usepackage[svgnames]{xcolor}

\usepackage{setspace}

\renewcommand{\phi}{\varphi}

\newcommand{\img}{\operatorname{Ran}}
\renewcommand{\ker}{\operatorname{Ker}}

\renewcommand{\Re}{\operatorname{Re}}
\renewcommand{\Im}{\operatorname{Im}}

\newcommand{\mc}[1]{\mathcal{#1}}

\newcommand{\bb}[1]{\mathbb{#1}}
\newcommand{\frk}[1]{\mathfrak{#1}}

\newcommand{\Z}{\mathbb{Z}}
\newcommand{\N}{\mathbb{N}}
\newcommand{\C}{\mathbb{C}}
\newcommand{\R}{\mathbb{R}}

\newcommand{\T}{\mathbb{T}}

\newcommand{\Ho}{\mathcal{O}}

\newcommand{\Aut}{\operatorname {Aut} }

\newcommand{\bs}{\backslash}

\newcommand{\ip}[1]{\langle #1 \rangle}
\newcommand{\iip}[1]{( #1 )}

    {\@beginparpenalty\predisplaypenalty
     \@endparpenalty\postdisplaypenalty
     \refstepcounter{equation}%
     \trivlist \item[]\leavevmode
       \hb@xt@\linewidth\bgroup $\m@th
         \displaystyle
         \hskip\mymathindent}%
        {$\hfil 
         \displaywidth\linewidth\hbox{\@eqnnum}%
       \egroup
     \endtrivlist}

\newcommand{\dint}{\ \!\mathrm{d}}
\newcommand{\dext}{\mathrm{d}}

\newcommand{\del}{\partial}
\newcommand{\delbar}{\overline{\partial}}

\newcommand{\ssubset}{\subset \!\!\! \subset}

\newcounter{dummy} \numberwithin{dummy}{section}

\newtheorem{dfn}[dummy]{Definition}
\newtheorem{exa}[dummy]{Example}
\newtheorem{rmk}[dummy]{Remark}
\newtheorem{lem}[dummy]{Lemma}
\newtheorem{prp}[dummy]{Proposition}
\newtheorem{thm}[dummy]{Theorem}
\newtheorem{cor}[dummy]{Corollary}

\newcommand{\SF}[3]{ \mathcal C^\infty_{#3}(#1; \Lambda^{#2}) }

\definecolor{darkred}{RGB}{99,00,00}
\newtheoremstyle{question}{\topsep}{\topsep}{\itshape}{}{\color{darkred}\bfseries}{.}{ }{}
\theoremstyle{question}

\newcommandx{\unsure}[2][1=]{\todo[linecolor=red,backgroundcolor=red!25,bordercolor=red,#1]{#2}}
\newcommandx{\change}[2][1=]{\todo[linecolor=blue,backgroundcolor=blue!25,bordercolor=blue,#1]{#2}}
\newcommandx{\info}[2][1=]{\todo[linecolor=OliveGreen,backgroundcolor=OliveGreen!25,bordercolor=OliveGreen,#1]{#2}}
\newcommandx{\improvement}[2][1=]{\todo[linecolor=Plum,backgroundcolor=Plum!25,bordercolor=Plum,#1]{#2}}
\newcommandx{\thiswillnotshow}[2][1=]{\todo[disable,#1]{#2}}

\newcommand{\LLl}{{\mathrm{L}}}

\newcommand{\Dli}{{\mathcal{D}'}}

\newcommand{\Eli}{{\mathcal{E}}'}

\newcommand{\EE}{\mathrm{E}}

\newcommand{\GG}{{\mathrm{G}}}

\newcommand{\pp}{\mathfrak{p}}
\newcommand{\qq}{\mathfrak{q}}

\newcommand{\supp}{\mathrm{supp}}

\newcommand{\AT}{\!\ ^tA}

\newcommand{\dd}{\operatorname{d}}

\newcommand{\loc}{\small \mathrm{loc}}
\newcommand{\vspp}{\vspace*{5pt}}
\newcommand{\lra}{\longrightarrow}
\newcommand{\sra}{\to}

\newcommand{\VV}{{\mathcal{V}}}

\usepackage{titleref}
\makeatletter
\newcommand*{\currentname}{\TR@currentTitle}
\makeatother

\graphicspath{{./figuras/}}             
\makeindex                              
\fontsize{60}{62}\usefont{OT1}{cmr}{m}{n}{\selectfont}
\normalsize

\begin{document}

\frontmatter
\onehalfspacing

\thispagestyle{empty}
\begin{center}
    \vspace*{2.3cm}
    \textbf{\Large{
        Top-degree solvability for hypocomplex \\
        structures and the cohomology of \\
        left-invariant involutive structures \\
        on compact Lie groups}}

     \vspace*{0.5cm}
    \Large{Max Reinhold Jahnke}

    \vskip 2cm
    \textsc{
    Dissertation \\[-0.25cm]
    presented to the \\[-0.25cm]
    Institute of Mathematics and Statistics\\[-0.25cm]
    of University of São Paulo\\[-0.25cm]
    in partial fulfillment of the requirements\\[-0.25cm]
    for the degree of\\[-0.25cm]
    doctor of philosophy}

    \vskip 1.5cm
    Program: Applied Mathematics\\
    Advisor: Prof. Dr. Paulo Domingos Cordaro\\

   	\vskip 1cm
    \normalsize{During the this work the author was partially supported by CNPq (process 140199/2014-4) and CAPES (PDSE 88881.131905/2016-01).}

    \vskip 0.5cm
    \normalsize{São Paulo, October 2018}
\end{center}

%
%
%
\newpage
\thispagestyle{empty}
    \begin{center}
        \vspace*{2.3 cm}
        \textbf{\Large{Top-degree solvability for hypocomplex structures and the cohomology of left-invariant involutive structures on compact Lie groups}}\\
        \vspace*{2 cm}
    \end{center}

    \vskip 2cm

    \begin{flushright}
	 This is the original version of the doctoral dissertation\\
	 by Max Reinhold Jahnke, as \\
	 submitted to the dissertation committee.
    \end{flushright}

\pagebreak

%
%
\newpage
\thispagestyle{empty}
\vspace*{\fill}
    \begin{center}
        Fire and Ice
        \vskip 0.5cm
	    \textit{Some say the world will end in fire,\\
        Some say in ice.\\
        From what I’ve tasted of desire\\
        I hold with those who favor fire.\\
        But if it had to perish twice,\\
        I think I know enough of hate\\
        To say that for destruction ice\\
        Is also great\\
        And would suffice.\\}
        \vskip 0.5cm
        Robert Frost
    \end{center}
 \vspace*{\fill}

\pagebreak

\chapter*{Acknowledgment}

First, I thank my advisor, Professor Paulo Domingos Cordaro. 
The best teacher I ever had, an excellent advisor and the greatest example of professionalism I know.
In these last years, besides being my advisor, Professor Cordaro became one of my dearest friends. I am very proud of being one of his students.

I also am grateful to Professor Gerardo Mendoza, for supervising me during my visit to Temple University and for providing excellent conditions for me to present my work, make professional contacts and make new friends.
Professor Mendoza pointed out many ways to improve my thesis and helped me to achieve much better results.

I also would like to thank my wonderful wife Adèle Helena Ribeiro for being an excellent partner, for taking care of me, for helping me make plans and keeping me focused on my career all this while she also was pursuing her doctorate degree. I will be always grateful for your dedication and love, Adèle.

A special gratitude goes out to my friend Gabriel C.C.S. de Araújo for all the encouraging, interest and suggestions that helped me during the development and writing of this thesis. Without his help, I would for sure have accomplished less in this work.

A special thanks to my friends Luis Fernando Ragognette, Luca Pallucchini, Farhan Abedin and Professor Gustavo Hoepfner for the company and all the support during my visit to Philadelphia. I am also grateful and to my friends Pedro Henrique Pontes and Sofia Tirado Pontes for all my support during my visit to the New York and New Jersey.

I am also thankful to my dear friends, Nicholas Braun Rodrigues, Antonio Victor, Bruno de Lessa, and Gregorio Chinni, with whom I not only learned a lot of mathematics but had a wonderful time drinking beer, talking, and dealing with the stress and uncertainties of the academic life.

Finally, I am thankful to my father Horst Reinhold Jahnke and my mother Gilda Timóteo Leite. Without the education they gave me and their support during the last years, it would have been very hard to focus on my career. I also thank my brothers Viktor Jahnke and Cristiane Jahnke for, despite being younger than me, always serving as an inspiration and as good examples.

\chapter*{Abstract}

\noindent Jahnke, M. R. \textbf{Top-degree solvability for hypocomplex structures and the cohomology of left-invariant involutive structures on compact Lie groups}. 2018.
Doctoral dissertation - Institute of Mathematics and Statistics, University of São Paulo, São Paulo.\\

    We use the theory of dual of Fréchet-Schwartz ({\bf DFS }) spaces to establish a sufficient condition for top-degree solvability for the differential complex associated to a hypocomplex locally integrable structure. As an application, we show that the top-degree cohomology of left-invariant hypocomplex structures on a compact Lie group can be computed only by using left-invariant forms, thus reducing the computation to a purely algebraic one. In the case of left-invariant elliptic involutive structures on compact Lie groups, under certain reasonable conditions, we prove that the cohomology associated to the involutive structure can be computed only by using left-invariant forms. \\

\noindent \textbf{Keywords:} Involutive Structures, Hypocomplex, Solvability, Lie groups, Cohomology, Left-invariant Cohomology.

\chapter*{Resumo}

\noindent Jahnke, M. R. \textbf{Resolubilidade em grau máximo para estruturas hipocomplexas e a cohomologia de estruturas involutivas invariantes à esquerda em grupos de Lie compactos}. 2018.
Dissertação de doutorado - Instituto de Matemática e Estatística da Universidade de São Paulo, São Paulo.\\

    Usamos a teoria da espaços duais de Fréchet-Schwartz para estabelecer uma condição suficiente para resolubilidade em grau máximo para o complexo associado a estrutuas localmente integráveis hipocomplexas. Como aplicação, provamos que a cohomologia de estruturas hipocomplexas invariantes à esquerda podem ser calculadas usando apenas formas invariantes à esquerda, assim reduzindo o cálculo a um método puramente algébrico. No caso de estruturas invariantes à esquerda, sob certas condições razoáveis, provamos que a cohomologia associada à estrutura pode ser calculada usando apenas formas invariantes à esquerda. \\

\noindent \textbf{Keywords:} Estruturas involutivas, hipocomplexidade, resolubilidade, Grupos de Lie Cohomologia, Cohomologia invariante à esquerda.

\singlespacing

\tableofcontents

\pagenumbering{roman}

\mainmatter

\chapter{Introduction}



This work has two main objectives concerning the cohomology of the differential complex associated to involutive structures.

One objective is to establish a sufficient condition for top-degree solvability for the differential complex associated to a hypocomplex locally integrable structure. Of course, here we are interested in the case in which the underlying manifold has no compact connected component.

The other objective is to extend a result by Chevalley and Eilenberg, published in \cite{chevalley1948cohomology}. That result states that the de Rham cohomology of compact Lie groups and of compact homogeneous manifolds can be computed using only Lie algebras. This led us to question whether similar results could be proved for other involutive structures on compact Lie groups.

Many ideas for finding conditions for top-degree solvability came from the paper by Ramis, Verdier, and Ruget \cite{ramis1971dualite} and the paper by Serre \cite{serre1955theoreme}. These two papers deal with complex structures, which are richier than the context we want to work on: the context of hypocomplex structures. Therefore, we had some limitations to overcome. The main limitation is that we do not always have local solvability in hypocomplex structures. Then we had to assume local solvability in a certain degree of the associated complex. The other limitation is that the differential operator associated to hypocomplex structures do not always have dense range. Thus, we had to find conditions to garantee such topological property.

While pursuing the second objective, we noticed that, to relate certain analytical properties of the involutive structures to the algebraic properties of the underling Lie group, we have to consider involutive structures that are preserved by the action of the left multiplication. These structures are called left-invariant involutive structures.

Left-invariant involutive structures on compact Lie groups as well as structures preserved by the group actions on compact homogeneous manifolds have already been considered by some other authors. For example, Bott proved in \cite{bott1957homogeneous} that, under certain conditions, the Dolbeault cohomology of compact complex manifolds can be computed using only algebraic methods involving the related Lie algebras. In \cite{pittie1988dolbeault}, Pittie not only showed that the result by Bott can be applied to left-invariant Dolbeault cohomologies on compact semisimple Lie groups, but also characterized all left-invariant complex structures on even dimensional compact Lie groups. In \cite{charbonnel2004classification}, Charbonnel and Khalgui provided a classification of all left-invariant CR structures of maximal rank. In this work, we focus on extending their results to left-invariant elliptic involutive structures on compact Lie groups.

To give a more detailed description of the results obtained in this work and to give the reader an idea of the techniques used, we first introduce some notation.

Let $\Omega$ be an orientable manifold and let $\mathcal V \subset \C T \Omega$ be a hypocomplex involutive vector bundle. For convenience of the reader, we reviewed the basic definitions of the theory of involutive structures on Chapter \ref{chp:involutive_structures}. To each involutive structure, there exists an associated differential complex denoted by
\begin{equation}
\label{chp:intro:V_smooth_complex}
(\mathcal C^\infty(\Omega; \Lambda^{p,q}), \dext')
\end{equation}
with cohomology spaces denoted by $H^{p,q}_{\mathcal C^\infty}(\Omega; \mathcal V)$.

This differential operator can obviously be extended to currents and compactly supported currents. In these cases, the respective cohomologies are denoted by $H^{p,q}_{\mathcal D'}(\Omega; \mathcal V)$ and $H^{p,q}_{\mathcal E'}(\Omega; \mathcal V)$.

In terms of cohomology, we want to find conditions under $\mathcal V$ and $\Omega$ so that $H^{p,n}_{\mathcal C^\infty}(\Omega; \mathcal V) = 0$, in which $n$ is the rank of $\mathcal V$ and so is the top-degree of the complex. We require that $\Omega$ has no compact connected components. Also, $\mathcal V$ has to be hypocomplex and the associated operator $\dext'$  has to satisfy some estimates. The details can be found in Theorem \ref{thm:first}.

We proved this result by applying some results of Functional Analysis to the operator $$\dext' : \mathcal E'(\Omega; \Lambda^{p,0}), \dext') \to \mathcal E'(\Omega; \Lambda^{p,1}).$$ The required theorems are stated right at the beginning of Chapter \ref{chp:involutive_structures}. We know two proofs of this theorem. The one presented is the easiest proof we know. However, another proof can be obtained from Theorem \ref{thm:closed_range}.

This theorem gives us conditions so that the range of the operator $$\dext' : \mathcal C^\infty(\Omega; \Lambda^{p,n-1}), \dext') \to \mathcal C^\infty(\Omega; \Lambda^{p,n})$$ is closed. This result was obtained without topological restrictions over $\Omega$, that is, we allow $\Omega$ to have compact components. By dropping this topological assumption, we made this result useful for dealing with compact Lie groups. See Theorem \ref{thm:left_invariance_on_degree} and Remark \ref{rmk:adele_labels_everything}.

The proof of Theorem \ref{thm:closed_range} uses \v Cech cohomology and is an adaptation of ideas from \cite{ramis1971dualite} and \cite{serre1955theoreme}.

Let $G$ be a compact Lie group and let $\frk g$ be the complexification of its Lie algebra. Each complex subalgebra $\frk h \subset \frk g$ defines a left-invariant involutive structure, which we always denote by $\frk h$. In Section \ref{sec:basicdef}, we showed how to construct many examples of Lie subalgebras having properties that are interesting from the point of view of the theory of involutive structures.

Since we are considering involutive structures defined by an algebra $\frk h$, the cohomology of the associated differential operator is going to be denoted by $H^{p,q}(G; \frk h)$. The result by Chevalley and Eilenberg that inspired our work states that, in the de Rham case (when $\frk h = \frk g$), the cohomology spaces $H^{p,q}(G; \frk h)$ can be computed using only the Lie algebra $\frk h$. To show this, they introduced a new complex, today called Chevalley-Eilenberg complex, and proved that the cohomology of this complex is isomorphic to the de Rham cohomology.

It is straightforward to extend the definition of the Chevalley-Eilenberg complex to consider the complex induced by a given subalgebra. The complex induced by the Lie algebra $\frk h$ is denoted by $C^{p,q}_{\frk h}(\frk g)$ with cohomology spaces denoted by $H^{p,q}_{\frk h}(\frk g)$. We have two ways to look into these complexes. In Chapter \ref{chp:involutive_structures_on_compact_lie_groups} we construct a restriction of the usual analytical complexes and their cohomologies, called left-invariant cohomologies, and in Section \ref{sec:cohomology_of_lie_algebras} we define a completely algebraic version of them.

To be explicit, what Chevalley and Eilenberg proved is that, when $\frk h = \frk g$, it is possible to construct an isomorphism
\begin{equation}
\label{chp:intro:the_master_isomorphism}
\phi: H^{p,q}_{\frk h}(\frk g) \to H^{p,q}(G; \frk h).
\end{equation}

Finding conditions for the existence of such isomorphism is exactly the motivation for this part of the work. Specifically, we want to see in what conditions it is possible to construct an isomorphism between the cohomology classes of the complex associated to a left-invariant involutive structures on a compact Lie group and a cohomology of Lie algebras relative to the respective subalgebra.

For a general subalgebra $\frk h \subset \frk g$, it is possible to construct a homomorphism similar to \eqref{chp:intro:the_master_isomorphism} by directly applying the techniques of Chevalley and Eilenberg. We present such techniques in Section \ref{sec:ce_tech}. In addition, we proved in Lemma \ref{inclusion_left_invatiant_cohomology_injective} that the map \eqref{chp:intro:the_master_isomorphism} is always injective.

The result about the injectivity of the map \eqref{chp:intro:the_master_isomorphism}, by itself, is already a useful result. Since the homomorphism is injective, it is clear that the algebraic properties of $\frk h$ and $\frk g$ play a role in the dimension of $H^{p,q}(G; \frk h)$. However, we notice that there are known examples of involutive structures, such as the Example \ref{exa:torus_liouvile}, in which it is impossible for the homomorphism \eqref{chp:intro:the_master_isomorphism} to be surjective.
In the case of Example \ref{exa:torus_liouvile}, since the dimension of the cohomology is infinite, the homomorphism cannot be surjective if the parameter $\mu$ is a rational or a Liouville number. The ratinal case is proved in Lemma \ref{lem:torus_liouville} and the Liouville case is a consequence of Lemma \ref{lem:torus_liouville_2}.

We also notice that, in the general case, finding conditions so that the homomorphism is surjective can be very complicated. Thus, we broke the problem down into simpler problems that can give us insights on how to proceed.

We remark that, in the cases in which the homomorphism is surjective, clearly the cohomology spaces are finite dimensional. Therefore, it makes sense to look for conditions such that the cohomology spaces relative to a subalgebra are zero or at least finite dimensional. To deal with this question, our approach uses Hodge theory. We proved, in Section \ref{sec:hodge_theory}, that the cohomology is finite dimensional given that the structure is subelliptic.

After choosing a Hermitian metric on $G$, we can define the Hodge Laplacian of the differential complex \eqref{chp:intro:V_smooth_complex} (with $\Omega = G$):
$$ \Box^{p,q} : C^\infty(G, \Lambda^{p,q}) \to C^\infty(G, \Lambda^{p,q}).$$
By assuming that the Hodge Laplacian is subelliptic, we can construct an abstract parametrix that allows us to prove that the cohomology of the complex is isomorphic to the kernel of the Hodge Laplacian, that is
\begin{equation}
\label{kohn_laplacian_kernel_iso}
	H^{p,q}(G, \frk h) \cong \ker \Box^{p,q}.
\end{equation}

Since $G$ is a compact manifold and it is supposed that the operator $\Box^{p,q}$ is subelliptic, we have that $\ker \Box^{p,q}$ has finite dimension. This gives us more information on how to find conditions for the homomorphism \eqref{chp:intro:the_master_isomorphism} to be surjective.

The Hodge Laplacian is known to be subelliptic in many cases. It is well known that, in the case of elliptic involutive structures, the associated Hodge Laplacian is elliptic and thus subelliptic. Another important case is related to the Levi form, which, on a Lie group, can be easily explained as follows:

Let $\frk h^\perp_0 = \frk h^\perp \cap \frk g_\R^*$. The \emph{Levi form} at a point $\xi \in \frk h^\perp_0$ is the Hermitian form on $\frk h$, defined by $$\mc L_\xi(L,M) = (-i/2)\xi([L,\overline{M}]), \qquad L,M \in \frk h.$$
We shall say that $\mc L$ is \emph{nondegenerate} if, given any point $\xi \in \frk h^\perp_0$, with $ \xi \neq 0$, and any $L \in \frk h$, with $L \neq 0$, there is some $M \in \frk h$ such that
$$ \mc L_\xi(L,M) \neq 0.$$ If the involutive structure $\frk h$ has nondegenerate Levi form, under some conditions about the eigenvalues of the Levi form, the associated Hodge Laplacian is subelliptic \cite{hanges1995involutive}. 

This approach using Hodge decomposition was our first attempt to find conditions so that the cohomology has finite dimension. In addition, we found an interesting technique by Helgason \cite{helgasondifferential}. This technique is a clever application of Lie derivative to prove Chevalley and Eilenberg's theorem on compact Lie groups. We used the same technique in our case, but noticed that, even for left-invariant elliptic structures, we could not find general conditions so that the homomorphism \eqref{chp:intro:the_master_isomorphism} is surjective. We could obtain some interesting results imposing some restrictions on the degrees of the complex and on the compact Lie group. We describe in more detail the results at the following two paragraphs.

In Proposition \ref{prp:comp_lie:hypo_subelliptic_torus}, we proved that, in the case of left-invariant hypocomplex structures, if we assume the compact Lie group to be a torus, then the homomorphism \eqref{chp:intro:the_master_isomorphism} is, in fact, surjective in the degree in which the Hodge Laplacian is subelliptic.

We also proved that the homomorphism \eqref{chp:intro:the_master_isomorphism} is surjective in top-degree for a general compact Lie group if the involutive structure is hypocomplex. The proof is a combination of Theorem \ref{thm:li:thm:duality_of_left} and Serre duality. This is stated and proved in detail in Theorem \ref{thm:left_invariance_on_degree}.

This naturally leads to the following question: is it possible to find necessary and sufficient conditions such that the involutive structure associated to the Lie algebra $\frk h$ is hypocomplex? An important result by Baouendi, Chang, and Treves \cite{baouendi1983microlocal} says that a sufficient condition for hypocomplexity is that the Levi form associated to the involutive structure has at least one negative eigenvalue and at least one positive eigenvalue at each point.

The algebraic cohomology obviously depends only on the algebraic properties of the Lie algebra and the fixed Lie subalgebra. So, it also makes sense to start looking at the algebraic properties that make some differences on the usual cohomology. For example, how can these algebraic properties impact on the dimension of the cohomology spaces?

Lie theory makes it clear that the topological properties of the group are strongly entangled with the algebraic properties of its Lie algebra. Therefore, it is natural to ask: how are these topological properties related to the dimension of $H^{p,q}(G; \frk h)$?

Since we could obtain some results on the torus, this hinted us that it could be useful to focus on a case that is completely orthogonal to the torus: the case of semisimple Lie groups.

The paper by Pittie shows that the homomorphism \eqref{chp:intro:the_master_isomorphism} is surjective in the Dolbeault's case when the Lie group admits a simply connected and compact covering group, that is, when the original group is semisimple. He promised to prove the general case for every compact Lie group in a future paper, but unfortunatelly this paper has never came out. By studying his paper, we noticed that it is necessary to require that the group is simply-connected for applying Bott's theorem.

We feel obliged to state here that, by using a technique that we did not explore here, Alexandrov and Ivanov proved in \cite[Corollary 5.1]{alexandrov2001vanishing}, that the hypothesis of semisimplicity can be dropped provided that the Lie group admits a bi-invariant metric compatible with the left-invariant complex structure. We believe that it would be a good idea for a further work to study their technique. We think that it can be adapted to general left-invariant elliptic structures, or even left-invariant hypocomplex structures, on compact Lie groups.

Since the de Rham case and the Dolbeault case are extremal cases of elliptic cases, we started looking for techniques to combine both of them to study elliptic structures.

A sufficient condition for using an adapted version of the Chevalley and Eilenberg technique is to assume that a certain subgroup is closed. This subgroup is the Lie group constructed by integrating the real orbits of the left-invariant elliptic structures.

Finding this closed Lie subgroup naturally calls for considering a homogeneous manifold, the quotient of the original Lie group by the closed subgroup. This homogeneous manifold has a natural complex structure inherited from the original left-invariant elliptic structure.

By using Leray-Hirsh Theorem, we are able to combine techniques we learned with Pittie, Bott, Chevalley and Eilenberg to study the original cohomology.

This approach using homogeneous spaces is explored in detail in Chapter \ref{chp:homogeneous}. The theorem that shows that the separate techniques can be combined is Theorem \ref{thm:homogeneous:converges_in_E_2}. Its first consequence is Theorem \ref{thm:omegaIsRieaman}, in which we find conditions so that the Dolbeault cohomolgy is as simple as possible.

We also proved Theorem \ref{thm:homogeneous:homo_on_simply3}, which, in a sense, generalizes both Chevalley-Eilenberg and Pittie's Theorems. In this case, we do not have dimensional restrictions on the left-invariant elliptic structure, but we have to assume that the compact Lie group is semisimple. The problem in full generality is still open.

\chapter{Involutive structures}
\label{chp:involutive_structures}

In this chapter, we introduce the basic definitions regarding involutive structures along with some tools to study them. The main concepts introduced here are the concept of involutive structure, the associated differential operators, and their cohomology spaces, which are the main objects of study in this work.

\section{Preliminaries on functional analysis}

As pointed out before, we apply a suitable characterization of the homomorphisms with dense image defined between Fr\'echet-Schwartz spaces.
The advantage of this class relies on the fact that the strong dual of a Fr\'echet-Schwartz space is an inductive limit of a sequence of Banach spaces
with compact inclusions. The proof of the homomorphim theorem makes use of its well known version valid for the class of the Fr\'echet-Montel spaces, whose proof can be found in  \cite[Section 33]{kothe19679topological}:

\begin{thm}\label{thm:FM}
 Let $E$, $F$ be Fr\'echet-Montel spaces and let
$A:E\sra F$ be a continuous linear map with $A(E)$ dense in $F$. The following properties are equivalent:
\begin{enumerate}
\item $A$ is a homomorphism;
\item $A(E)=F$;
\item $\AT$ is a strong homomorphism;
\item $\forall {\cal{B}} \subset F'$,  $\AT({\cal{B}})\subset E'$  strongly bounded $\Rightarrow$ ${\cal{B}}$ strongly bounded;
\end{enumerate}\end{thm}

A particular class of Fr\'echet-Montel spaces is the class of the so-called Fr\'echet-Schwartz ({\bf FS}) spaces (cf. \cite{komatsu1967projective}).
When $E$ is a ({\bf FS}) space its strong dual is an inductive limit of a sequence of Banach spaces
$E'=\mbox{ind}\,( E'_p)$,
where, for each $p\in\N$, $E'_p$ is a subspace of $E'_{p+1}$, the injection map being compact.

\begin{thm}\label{thm:FS}
Let $E$, $F$ be Fr\'echet-Schwartz spaces and let
$A:E\sra F$ be a continuous linear map with dense image. Write the corresponding strong duals as inductive limits of Banach spaces
$E'=\mbox{ind}\,( E'_p)$, $F'=\mbox{ind}\,( F'_p)$,
where, for each $p$, the inclusions $E'_p\hookrightarrow E'_{p+1}$, $F'_p\hookrightarrow F'_{p+1}$ are compact. The following properties
are equivalent:
\begin{enumerate}
\item $A(E)=F$;
\item For each $p\in\N$, there is $q\in\N$ such that, if ${\cal{B}} \subset F'$ is such that $\AT({\cal{B}})$ is contained and bounded in $E'_p$,
then ${\cal{B}}$ is contained and bounded in $F'_q$;
\item For each $p\in\N$, there is $q\in\N$ such that, if $u\in F'$ is such that $\AT(u)\in E'_p$, then $u\in F'_q$.
\end{enumerate}
\end{thm}

\begin{proof} Since every strongly bounded subset of $E'$ is contained and
bounded in one of the Banach spaces $E'_p$ it is clear that (2) implies (4) of
Theorem \ref{thm:FM} and hence (1) holds. Conversely, if (1) holds then, by Theorem
\ref{thm:FM}, $\AT$ is a strong homomorphism. By \cite[Theorem 7']{komatsu1967projective},
$\AT(F')=\mbox{ind}\,(E'_p\cap \AT(F'))$, that is, $\AT(F')$ can also be described as the inductive limit of the sequence of Banach spaces
$G_p\doteq E'_p\cap \AT(F')$ with the corresponding compact inclusions $G_p\stackrel{\iota_p}{\hookrightarrow} G_{p+1}$. For each $p\in\N$, we set
$\lambda_p=(\AT)^{-1}|_{G_p}:G_p\sra F'$. Since $\lambda_p= \lambda_{p+1}\circ \iota_p$, it follows that each $\lambda_p$ is a compact linear map. Consequently,
if $U_p$ denotes the unit ball in $G_p$, then $\lambda_p(U_p)$ is strongly bounded in $F'$ and hence it is contained in $F'_q$ for some $q\in\N$
(and ``a fortiori'' $\lambda_p:G_p\sra F'_q$ is a continuous linear map between Banach spaces). We have
proved the following assertion: for each $p\in \N$, there is $q\in\N$ such that $(\AT)^{-1}|_{G_p} \in L(G_p, F'_q)$, from which (2) follows immediately.

Next, it is clear (2) implies (3) after taking ${\cal{B}}=\{u\}$. We are then left with the proof that (3) implies (2).
For this, it suffices so show that
the following strengthen form of (3) is satisfied (in the argument, we shall denote by $|\cdot|_p$ the norm in $E'_p$ and by $\|\cdot\|_p$ the norm in $F'_p$):
\begin{itemize}
 \item[($\bullet$)] {\sl Given $p\in\N$,  there are $q\in \N$ and $C>0$ such that,
if $u\in F'$ and $\AT u\in E'_p$, then $u\in F'_q$ and $\|u\|_q \leq C |\AT u|_p$}.
\end{itemize}

Indeed, for each $p\in \N$, take $q\in\N$ as in (3).
The space ${\mathfrak{E}}$ of all $u\in F'_{q+1}$ such that $\AT u\in E_p$
is a Banach space when we consider the norm $u \mapsto  \|u\|_{q+1} + |\AT u|_p$.
Property (3) gives ${\mathfrak{E}}\subset F'_q$ and the closed graph theorem implies the existence of a constant $\mu>0$ such that
$$   \|u\|_q \leq \mu \left( |\AT u|_p + \|u\|_{q+1}\right), \qquad u\in {\mathfrak{E}}. \leqno{(\sharp)} $$

If ($\bullet$) were not true, it would exist a sequence $\{u_j\}\subset F'_q$ such that $\{\AT u_j\}\subset E'_p$ and
$\|u_j\|_{q} > j|\AT u_j|_p$.
If we define $v_j=u_j/\|u_j\|_{q}$, then $\{v_j\}\subset E'_q$, $\|v_j\|_{q}=1$ for every $j$ and $\AT v_j\sra 0$ in $F'_p$. Since the inclusion
$F'_q\subset F'_{q+1}$ is compact, some subsequence $\{v_{j_k}\}$ of $\{v_j\}$ must converge in $F'_{q+1}$. Denoting by
$v\in F'_{q+1}$
the limit of this subsequence, we have $\AT v=0$, and hence $v=0$ since $A(E)$ is dense in $F$. On the other hand ($\sharp$) implies
$1 \leq \mu \left\{ |\AT v_{j_k}|_{p} + \|v_{j_k}\|_{q+1} \right \}$, $k=1,2,\ldots,$
and consequently $\|v\|_{q+1}\geq 1/\mu$, which is a contradiction. \end{proof}
\vspp

\section{The associated differential complex}

Let $\Omega$ be a smooth manifold and denote by $N$ its dimension. We always assume that $N \geq 2$. An involutive structure on $\Omega$ is a smooth subbundle $\mathcal V$ of the complexified tangent bundle $\C T \Omega$ of $\Omega$ that satisfies the involutive, or Frobenius, condition $[\mathcal V, \mathcal V] \subset \mathcal V$. We always denote the rank of $\mathcal V$ by $n$. We denote by $\Lambda^k$ the bundle $\Lambda^k \C T^* \Omega$ and by $\mathcal C^\infty(\Omega; \Lambda^k)$ the space of sections of $\Lambda^k$ with smooth coefficients. 


There is a natural differential complex associated to each involutive structure $\mathcal{V}$. We briefly recall the construction of such complexes. For each $x \in \Omega$, we define $$T'_x = \{ u \in \C T^*_x \Omega: u(X) = 0 ,~ \forall X \in \mathcal V_x \}$$ and we have that
$ T' = \bigcup_{x \in \Omega} T'_x $
is a vector subbundle of $\C T^*\Omega$. For each $p, q \in \Z_+$, we denote by $T^{p,q}_x$ the subspace of $\Lambda^{p+q} \C T_x^*\Omega$ consisting of linear combinations of exterior products $u_1 \wedge \ldots \wedge u_{p+q}$ with $u_j \in \C T_x^*\Omega$ for $j = 1, \ldots, p+q$ and at least $p$ of these factors belonging to $T'_x$.
Notice that $T^{p+1,q-1}_x \subset T^{p,q}_x$ and so we can define $\Lambda^{p,q}_x \doteq  T^{p,q}_x / T^{p+1,q-1}_x$ and consequently $\Lambda^{p,q} = \bigcup_{x \in \Omega} \Lambda^{p,q}_x $ is a smooth vector bundle over $\Omega$. Since $\mathcal{V}$ is involutive, we can easily verify that the exterior derivative take smooth sections of $T^{p,q}$ into smooth sections of $T^{p,q+1}$ and hence there exist a unique operator $\dext'_{p,q}$ such that the diagram
$$
\xymatrix{
 T^{p,q} \ar[r]^{\dext} \ar[d]_{\pi_{p,q}}& T^{p,q+1} \ar[d]^{\pi_{p,q+1}} \\
\Lambda^{p,q} \ar[r]_{\dext'_{p,q}} & \Lambda^{p,q+1}
}
$$
is commutative, with $\pi_{p,q} : T^{p,q} \to \Lambda^{p,q}$ being the quotient map.

When there is no risk of confusion, we simplify the notation and write simply $\dext '$
for the operator $\dext'_{p,q}$. We refer as $\dext'_{\mathcal V}$ when it is necessary to emphasize the associated involutive structure.

We denote by $\mathcal C^\infty(\Omega; \Lambda^{p,q})$ 
the space of sections of $\Lambda^{p,q}$ with smooth coefficients. 
Notice that the operator $\dext'$ maps smooth sections of $\Lambda^{p,q}$ into smooth sections of $\Lambda^{p,q+1}$ and also that $\dext' \circ \dext' = 0$. Therefore, for each $p \geq 0$, the operator $\dext'$ defines a complex of $\C$-linear mappings
\begin{equation}
\label{h_complex}
 \mathcal C^\infty(\Omega; \Lambda^{p,0}) \xrightarrow{\dext'} \mathcal C^\infty(\Omega; \Lambda^{p,1})  \xrightarrow{\dext'} ... \xrightarrow{\dext'} \mathcal C^\infty(\Omega; \Lambda^{p,q}) \xrightarrow{\dext'} \mathcal C^\infty(\Omega; \Lambda^{p,q+1}) \xrightarrow{\dext'} \ldots
\end{equation}


Of course we obtain another differential complex after replacing $\mathcal C^\infty$ by $\mathcal C^\infty_c$, where $c$ stands for compact supports.

For each $p \geq 0$, we denote the set of the $(p,q)$-cocycles elements by
$$ Z^{p,q} (\Omega; \mathcal V) = \ker \left( \dext' : \mathcal C^\infty(\Omega; \Lambda^{p,q}) \to \mathcal C^\infty(\Omega; \Lambda^{p,q+1}) \right),$$
the set of $(p,q)$-coboundaries by
 $$ B^{p,q} (\Omega; \mathcal V) = \img \left( \dext' : \mathcal C^\infty(\Omega; \Lambda^{p,q-1}) \to \mathcal C^\infty(\Omega; \Lambda^{p,q}) \right),$$
and the $(p,q)$-cohomology classes by
 $$  H^{p,q} (\Omega; \mathcal V) = \frac{Z^{p,q} (\Omega; \mathcal V)}{B^{p,q} (\Omega; \mathcal V)}$$
 (we set $B^{p,0}(\Omega)=\{0\}$).

Given any open set $V \subset \Omega$ we can restrict the involutive structure $\mathcal V$ to $V$, which we still denote by $\mathcal V$, and construct a complex $\left( \mathcal C^\infty(V; \Lambda^{p,q}), \dext' \right)$. And, by restriction, we have a homomorphism of differential complexes
$$ \left( \mathcal C^\infty(\Omega; \Lambda^{p,q}), \dext' \right) \to \left( \mathcal C^\infty(V; \Lambda^{p,q}), \dext' \right).$$ This then  allows us to define, for a given $x \in V$,  $$ \mathcal H^{p,q} (\{x\}; \mathcal V) = \varinjlim_{V \ni x} H^{p,q}(V; \mathcal V) $$ the inductive limit of $H^{p,q}(V; \mathcal V)$ taken when  $V$ varies in the set of all open neighborhoods of $x$.

\begin{exa}
By taking $\mathcal V = \C T \Omega$, we have the simplest involutive structure. In this case, the operator $\dext'$ is just the usual exterior derivative. The space $H^p_{\mathcal C^\infty} (\Omega; \mathcal V)$ is just the de Rham cohomology space.
\end{exa}

\begin{exa}
Suppose that $J : T \Omega \to T \Omega$ is a vector bundle isomorphism satisfying $J \circ J = -\operatorname{Id}$ and
\begin{equation}
\label{chp:inv:j_complex}[JX, JY] = [X,Y] + J[JX, Y] + J[X,JY].
\end{equation} We extend $J$ complex-linearly to the complexified tangent bundle $\C T \Omega$ and we take $\mathcal V = \ker(J - i{\operatorname{Id}})$. The condition \eqref{chp:inv:j_complex} implies that  $\mathcal V$ is involutive. This structure is called a \emph{complex structure} over $\Omega$ and, in this case, the operator $\dext'$ is the Doubeault operator $\delbar$.
\end{exa}
It is also important to consider the preceeding complexes with more general coefficients. In order to do so, we assume from now on that
$$   \mbox{$\Omega$ is orientable.} $$
Under this hypothesis, we can also construct a complex similar to \eqref{h_complex}
where now the coefficients are distributions over $V$. In this case, we denote the cohomology space by
$   H^{p,q}_{\mathcal D'} (V; \mathcal V).$ In general, if $F(V)$  is a subspace
of $\Dli(V)$, we denote by $F(\Omega; \Lambda^{p,q})$ 
the space of sections of $\Lambda^{p,q}$ with coefficients in $F$ and we denote the cohomology of $\dext'$ with coefficients in $F(V)$ by $ H^{p,q}_{F} (V; \mathcal V).$
Furthermore, by \cite[Propositions VIII.1.2 and VIII.1.3]{treves1992hypo}, there is a natural bracket
which turn the spaces
$$ \mathcal C_c^\infty(V;\Lambda^{p,q})\mbox{ and } \Dli(V;\Lambda^{m-p,n-q})
\,\, (\mbox{resp. }   \mathcal C^\infty(V;\Lambda^{p,q})\mbox{ and }
\mathcal{E}'(V;\Lambda^{m-p,n-q})) $$
into the dual of one another and in such a way that the transpose of $\dd'$ is also $\dd'$ (up to a sign).

We can endow each space $\mathcal C^\infty(V;\Lambda^{p,q})$ with a locally convex structure of an ({\bf{FS}}) space. Its dual $\mathcal{E}'(V;\Lambda^{m-p,n-q})$ is then a ({\bf{DFS}}) space and a sequence of definition for its topology can be taken by the  sequence
$$     \GG_j(V;\Lambda^{m-p,n-q}) := \{u\in H^{-j}_{\tiny{\mbox{loc}}}(V,\Lambda^{m-p,n-q}): \,\mbox{supp}\, u \subset K_j \}, $$
where $\{K_j\}$ is an exhaustion of $V$ by compact sets (with the inclusions $$\GG_j(V;\Lambda^{m-p,m-q}) \hookrightarrow \GG_{j+1}(V;\Lambda^{m-p,m-q})$$ being compact by the Relich lemma).

We conclude this section by stating a direct consequence of Theorem \ref{thm:FS}.

\begin{thm}\label{thm:main}
Let $\Omega$ be an orientable differentiable manifold endowed with a formally integrable structure. Assume that
\begin{equation}\label{density}
\dd' : \mathcal{E}'(\Omega;\Lambda^{m-p,0}) \lra \mathcal{E}'(\Omega;\Lambda^{m-p,1})
\end{equation}
is injective. Then the map
\begin{equation}\label{1}
 \dd':  \mathcal C^\infty(\Omega;\Lambda^{p,n-1}) \lra \mathcal C^\infty(\Omega;
\Lambda^{p,n}) \end{equation}
is surjective (which is the same as saying that $H^{p,n}(\Omega; \VV)=0$) if and only if the following property holds:
\begin{itemize}
\item[$(\dagger)$] Given $K\subset \Omega$ compact and $s\in \R$, there are
 $K_1\subset \Omega$  compact and  $t\in \R$ such that the following holds for an arbitrary $u\in \mathcal{E}'(\Omega;\Lambda^{m-p,0})$:
$$ \mathrm{supp}\,\dd'u \subset K,\,\, \dd' u \in H^s_{\tiny{\mathrm{loc}}}(\Omega,\Lambda^{m-p,1}) \Longrightarrow  \mathrm{supp}\, u\subset K_1,\,\, u\in H^t_{\tiny{\mathrm{loc}}}(\Omega,\Lambda^{m-p,1}). $$
\end{itemize}
\end{thm}

Indeed, Property \eqref{density} is equivalent to the fact that \eqref{1} has dense image and hence Theorem \ref{thm:FS} applies immediately.

Under an additional property, we can strength the result stated in Theorem \ref{thm:main}. We introduce the following definition:

\begin{dfn}\label{dfn:contun} We shall say that $\VV$ satisfies property $(\star)$ if given $V\subset \Omega$ open and $u\in \Dli(V;\Lambda^{p,0})$ satisfying $\dd' u=0$ in a nonempty  connected open subset
$\omega$ of $V$ then $u$ vanishes identically in the component of $V$ that contains $\omega$.
\end{dfn}

\begin{prp}\label{prp:imp} Let $\Omega$ be an orientable differentiable manifold endowed with a formally integrable structure. Suppose that $\Omega$ has no compact connected component and that $\VV$ satisfies $(\star)$. Then
\begin{enumerate}
 \item  The kernel of the map $\dd':\mathcal E'(\Omega;\Lambda^{m-p,0})\sra \mathcal E'(\Omega;\Lambda^{m-p,1})$ is trivial;
\item If $K\subset \Omega$ is compact and if $u\in \Eli(\Omega;\Lambda^{m-p,0})$ is such that $\mbox{supp}\,\dd 'u\subset K$ then $\mbox{supp}\, u\subset \tilde K$,
where $\tilde K$ is the compact set obtained as the union of $K$ with the relatively compact components of $\Omega\setminus K$.
\end{enumerate}
\end{prp}

The proof is immediate.
\vspp

From Theorem \ref{thm:main} and Proposition \ref{prp:imp} we obtain:

\begin{cor}\label{cor:Main} Let $\Omega$ be an orientable differentiable manifold endowed with a formally integrable structure. Suppose that $\Omega$ has no compact connected component and that $\VV$ satisfies $(\star)$.
If for a given $K\subset \Omega$ compact
and for a given $s\in \R$ there is $t\in \R$ such that
$$ u\in \mathcal E'(\Omega,\Lambda^{m-p,0}),~ \supp\, \dext' u\subset K,~ \dext' u\in H^s_{\tiny\mathrm{loc}}(\Omega;\Lambda^{m-p,1}) \Longrightarrow u\in
H^t_{\tiny\mathrm{loc}}(\Omega;\Lambda^{m-p,0}) $$
then $H^{p,n}(\Omega;\VV)=0$.
\end{cor}


\section{Hypocomplex structures}

From now on, we assume the stronger property that the structure $\mathcal V$ is {\sl locally integrable}. This means that, in a neighborhood of an arbitrary point $\pp\in \Omega$, there are defined $m=N-n$ smooth functions whose differentials span
$T'$ at each point in a neighborhood of $\pp$. Notice that each of these functions is annihilated by the operator  $\dd'$, that is, they are solutions for $\VV$.
In this section, we are going to introduce the concept of hypocomplex structures, a class of structures with some
good analytical properties, similar to the class of holomorphic functions. Such analytical properties are going to
be used to find conditions so that the operator $\dext'$ in maximum degree has closed range.

Let $\pp \in \Omega$. We denote by $\mathcal S(\pp)$ the ring of germs of solutions of $\mathcal V$ at $\pp$, that is,
$$ \mathcal S(\pp) = \{ f \in \mathcal C^\infty(\pp) : \dext' f = 0 \},$$
and we denote by $\mathcal O^m_0$ the ring of germs holomorphic functions defined in a neighborhood of $0$ in $\C^m$.

Now, consider $Z_1, \ldots, Z_m$ solutions for $\VV$, defined in an open neighborhood of $\pp$ with linearly independent differentials, and denote by $Z$ the map $Z = (Z_1, \ldots, Z_m)$ defined in a neighborhood of $\pp$. Of course, we assume that $Z(\pp)=0$.

We define a ring homomorphism
$$ \lambda : h\in \mathcal O^m_0 \mapsto  h \circ Z \in  \mathcal S(\pp). $$
It is not difficult to prove that $\lambda$ is injective. We shall say that $\mathcal V$ is \emph{hypocomplex} at $\pp$ if $\lambda$ is surjective, that is, if for every $f \in \mathcal S(\pp)$ there exist a holomorphic function $h \in \mathcal O^m_0$ such that $f = h \circ Z$. 

The next result follows directly from the definitions:

\begin{prp}\label{prp:hypoxstar}
If $\VV$ is hypocomplex at each point, then $\VV$ satisfies Property $(\star)$.
 \end{prp}

A particular but important class of hypocomplex structures is that of the elliptic structures. Recall that the {\sl characteristic set of} $\VV$ is the set $C(\VV) = \VV^\perp \cap T^*\Omega$ and
that the structure is {\sl elliptic} if $C(\VV)=0$. Every elliptic structure, being hypocomplex ([T4, Proposition III.5.1]),
satisfies Property ($\star$) and, moreover,
given $\pp\in\Omega$ and $s\in\R$, the following is true:
if $u\in{\cal{D}}'(\Omega,\Lambda^{m-p,0})$ and $\dd' u\in H^s_{\loc}(\{\pp\};\Lambda^{m-p,1})$
then $u\in H^{s+1}_{\loc}(\{\pp\},\Lambda^{m-p,0})$. From Corollary \ref{cor:Main} follows:

\begin{thm}
Suppose that $\Omega$ has no compact connected component and that $\VV$ is elliptic.
Then $H^{p,n}(\Omega;\VV)=0$.
\end{thm}

Another important class of hypocomplex structures is the one described in the following. Assume that $\Omega$ is endowed with a locally integrable structure $\VV$. Let
 $T^0=  T' \cap \C T^* \Omega$ denote the characteristic set of $\VV$. Given $(\pp,\xi)\in T^0$, $\xi\neq 0$, the \emph{Levi form} at  $(\pp,\xi) \in T^0$ is the Hermitian form on $\mathcal V_x$ defined by $$\mc L_{(\pp,\xi)}(v,w) = (1/2i)\xi([L,\overline{M}]_\pp), \quad v,w \in \mathcal V_\pp,$$
 in which $L$ and $M$ are any smooth local sections of $\mathcal V$ in a neighborhood of $\pp$ such that $L_\pp = v$ and $M_\pp = w$. We shall say that $\mc L$ is \emph{nondegenerate} if given any point $(\pp,\xi) \in T^0$, with $ \xi \neq 0$, the hermitian form $\mc L_{(\pp,\xi)}$ is nondegenerate.

The following result is due to Baouendi-Chang-Treves \cite{baouendi1983microlocal}:

\begin{thm}\label{thm:BCT}
Let $\Omega$ be endowed with a locally integrable structure for which the Levi form at
each $(\pp,\xi)\in T^0$, $\xi\neq 0$, has one positive and one negative eigenvalue. Then
$\VV$ is hypocomplex.
\end{thm}

Assume now that the locally integrable structure $\VV$ has rank $N-1$, that is, $T'$ is a complex line bundle. For such, structures hypocomplexity is completely characterized \cite{treves1992hypo}:

\begin{thm}\label{thm:open}
If the locally integrable structure $\VV$ has rank $N-1$, then $\VV$ is hypocomplex if and only if given any $\pp\in \Omega$ and given any local solution $Z$ for $\VV$ near
$\pp$ with $\dd'Z_\pp\neq 0$ then $Z$ is open at $\pp$.
\end{thm}


\section {A remark on local solvability}

Let $\Omega$ be endowed with a locally integrable structure $\VV$. We shall now make use of the local representation for $\VV$ as described in \cite{treves1992hypo}, \cite{berhanu2008introduction}. Any point $\pp\in\Omega$ is the center of a coordinate system $(x,t) = (x_1, \ldots, x_m, t_1, \ldots, t_n)$ defined on an open neighborhood of $p$ denoted by $U$ such that $(x(\pp),t(\pp)) = (0,0) \in \R^m \times \R^n$.
On $U$ it is defined a smooth function $\phi(x,t) = (\phi_1(x,t), \ldots,$ $\phi_m(x,t))$ satisfying $\phi(0,0) = 0$ and $\phi_x(0,0) = 0$ such that the differentials of the functions
\begin{equation}
\label{chp:inv:first_integrals}
Z_k(x,t) = x_k + i\phi_k(x,t), \quad k = 1, \ldots, m
\end{equation}
 are linearly independent and span $T'$ over $U$. If we define the vector fields $$M_k = \sum_{j=1}^m \mu_{kj}(x,t) \frac{ \del }{ \del x_j}, \quad k = 1, \ldots, m,$$
characterized by the rule $M_kZ_j = \delta_{k,j}$ for $k,j = 1, \ldots, m$, then the complex vector fields
$$ L_j = \frac{ \del }{ \del t_j} - i\sum_{k=1}^m \frac{ \del \phi_k }{ \del t_j }(x,t)M_k, \quad j = 1, \ldots, n, $$
span $\mathcal V$ over $U$. The vector fields $L_j$ and $M_k$ commute pairwise.

The differentials $\dext t_k$, for $k = 1, \ldots, n$, span a bundle over $U$, which is supplementary to $T'|_U$. Therefore we can adjoin them with $\dext Z_k$, for $k = 1, \ldots, m$, to get a smooth basis for $\C T^* \Omega$. So, if $f$ is a smooth section of $\Lambda^{p,q}$ over $U$, we have the following unique representation for $f$:
$$ f = \sum_{|J| = p} \sum_{|K| = q} f_{J,K}(x,t) \dext Z_J \wedge \dext t_K.$$

The standard representation of $\dext'f$ over $U$ is then given by
$$ \dext' f = \sum_{|J| = p} \sum_{|K| = q} \sum_{l=1}^n \left( L_l  f_{J,K} \right) (x,t) \dext t_l \wedge \dext Z_J \wedge \dext t_K. $$

The goal of this section is to prove the following result:

\begin{prp}
\label{prp:equiv} Let $\Omega$ be endowed with a locally integrable structure $\mathcal V$ and let $\pp \in \Omega$. The following properties are equivalent, for a fixed $q\in\{1,\ldots,n\}$:
\begin{enumerate}
\item \label{hypo:prp:equiv1} $ \mathcal H^{0,q}_{\mathcal C^\infty} (\{\pp\}, \mathcal V) = 0$;
\item \label{hypo:prp:equiv2} $ \mathcal H^{0,q}_{\mathcal D'} (\{\pp\}, \mathcal V)= 0$;
\item $ \mathcal H^{p,q}_{\mathcal C^\infty} (\{\pp\}, \mathcal V) = 0$ for all $p = 0, 1, 2, \ldots, m$;
\item $ \mathcal H^{p,q}_{\mathcal D'} (\{\pp\}, \mathcal V)= 0$ for all $p = 0, 1, 2, \ldots,m$.
\end{enumerate}
\end{prp}

\begin{proof} From the discussion that precedes, it is clear that it suffices to prove the equivalence between (1) and (2).
We first prove that \eqref{hypo:prp:equiv1} implies \eqref{hypo:prp:equiv2} by using in part an argument in \cite{andreotti1981absence}. We keep the notation just established
and apply the well known Grothendieck argument as in \cite[Theorem VII.6.1]{berhanu2008introduction}:
given $U_0\subset U$ an open neighborhood of the origin, there corresponds another open neighborhood $V_0\ssubset U_0$ such that, if we set
$$    \EE := \{(f,u)\in Z^{0,q}(U_0)\times C^\infty(V_0,\Lambda^{0,q-1}) :
\dd' u =f|_{V_0} \},$$
then the projection map $\pi:\EE \sra Z^{0,q}(U_0)$ (a continuous linear map between Fréchet spaces) is surjective and hence open. If we follow \cite{andreotti1981absence}, we can ensure:

\begin{itemize}
\item Given $K'\subset V_0$ compact and $\nu'$ a positive integer, there are $K\subset
U_0$ compact, $\nu$ a positive integer and $C>0$ such that the following is true: if $f\in Z^{0,q}(U_0)$, there is $u\in C^\infty(V_0,\Lambda^{0,q-1})$ such that
$$     \|u\|_{\nu',K'}\leq C\|f\|_{\nu,K}. $$
\end{itemize}
We obtain
$$ \|f\|_{\nu,K} \doteq \sup_{|\alpha| \leq \nu} \sup_K |D^\alpha f|. $$

In order to conclude (2), we make use of \cite[Theorem II.5.2]{treves1992hypo}. It suffices to establish the following solvability property:

\begin{itemize}
\item Given an open neighborhood $U_1$ of the origin, there is $\mu$ such that, if $f\in \mc C^\mu(U_1;\Lambda^{0,q})$ satisfies $\dd'f=0$, we can find $u\in \Dli(\{0\};\Lambda^{0,q-1})$ solving $\dd'u =f$.
\end{itemize}

We first apply the approximate Poincaré lemma \cite[Section II.6]{treves1992hypo}: there is
an open neighborhood $U_0\ssubset U_1$ of the origin such that the conclusion of \cite[Theorem II.6.1]{treves1992hypo} holds. Now we refer to the discussion presented before: taking a sufficiently small closed ball
$K_0$ centered at the origin, there are $\mu$ and $C>0$ such that, for every
$g\in Z^{0,q}(\overline{U_0};\Lambda^{0,q})$, there is
$v\in C^\infty(K_0,\Lambda^{0,q-1})$ such that $\dd' v = g$ in $K_0$
and
$$       \|v\|_{K_0,0}\leq C \|g\|_{\overline{U_0},\mu}. $$
This defines our sought $\mu$. Given $f\in C^\mu(U_1,\Lambda^{0,q})$, we can find $f_j\in Z^{0,q}(\overline{U_0})$ such that
$f_j\sra f$ in $\mc C^\mu (\overline{U_0};\Lambda^{0,q})$. Let $u_j \in C^\infty(K_0,\Lambda^{0,q-1})$ such that $\dd' u_j = f_j$ in $K_0$
and
$$       \|u_j\|_{K_0,0}\leq C \|f_j\|_{\overline{U_0},\mu}. $$
This shows that $u_j$ is bounded in $L^\infty(K_0;\Lambda^{0,q-1})$ and hence some subsequence of $u_j$ converges weakly to some $u$ defined in the interior of $K_0$, which a fortiori satisfies $\dd'u=f$.

For the proof that (2) implies (1), we again start by applying the Grothendieck argument, this time in the following way. Let $U_0$ be a neighborhood of the origin, let
$\omega_j\ssubset U_0$ be a fundamental system of neighborhoods of the origin
and define
$$    \EE_j \doteq \{(f,u)\times Z^{0,q}(U_0)\times H^{-j}_{\mbox{\small loc}}(\omega_j):
 \dd 'u=f|_{\omega_j} \}. $$
By the standard argument, we can find $V_0\subset U_0$ an open neighborhood of the origin and $\sigma\in\R$ such that the following is true:
\begin{itemize}
\item Given $f\in Z^{0,q}(U_0)$, there is $u\in H^{-\sigma}_{\mbox{\small loc}}(V_0,\Lambda^{0,q-1})$ such that $\dd'u=f|_{V_0}$.
\end{itemize}

The proof completes  after applying the argument presented in the proof of  \cite[Theorem VIII.9.1]{treves1992hypo}
\end{proof}

\begin{rmk}\label{rmk:1}
Inspection of the proof of Proposition \ref{prp:equiv} shows that Properties (1) to (4)
in Proposition \ref{prp:equiv} are still equivalent to the following one:
\begin{itemize}
\item[5.] Given $s\in\R$, there is $t\in\R$ such that, if $f\in H^{s}_{\mbox{\small loc}}(\{\pp\},\Lambda^{p,q})$ satisfies $\dd'f=0$, then there is $u\in H^{t}_{\mbox{\small loc}}(\{\pp\},\Lambda^{p,q-1})$ such that $\dd' u=f$, $p=0,1,\ldots,m$.
\end{itemize}
\end{rmk}

Next, we state the main theorem of this chapter:

\begin{thm}\label{thm:first}  Let $\Omega$ be orientable with no compact connected component. Assume also that $\Omega$ is endowed with a hypocomplex structure $\mathcal V$. If the equivalent conditions in Proposition \ref{prp:equiv} hold for $q=1$, then
$$   H^{p,n}_{\mathcal C^\infty}(\Omega) = 0,\qquad p = 0, 1, 2, \ldots.$$
\end{thm}

\begin{proof} We apply Corollary \ref{cor:Main}. Let $K\subset \Omega$ be compact, $s\in\R$ and assume that $u \in \mc E '(\Omega;\Lambda^{m-p,0})$
satisfies $\supp~ u\subset K$, $\dd'u \in H^{s}_{\mbox{\small loc}}(\Omega;\Lambda^{m-p,1})$. By Remark \ref{rmk:1}, for each $\pp\in \supp~ u$ we can find $t\in \R$ and $v_\pp\in
H^{t}_{\mbox{\small loc}}(\{\pp\},\Lambda^{p,0})$ with $\dd'v_\pp= \dd'u$ near $\pp$.
By compactness, we can select $t$ independent of the choice of $\pp$. Furthermore, since $\VV$ is hypocomplex, it follows that $u-v_\pp$ is smooth near $\pp$ and hence
$u\in  H^{t}_{\mbox{\small loc}}(\{\pp\},\Lambda^{p,0})$. Since $\pp$ is arbitrary, Corollary \ref{cor:Main} gives the result. \end{proof}

\section{Examples}

We now apply Corollary \ref{cor:Main} and Theorem \ref{thm:first} to derive some examples:



\begin{exa}\label{exa:2}
Let ${\cal{M}}$ be an oriented, real-analytic manifold of dimension $n\geq 1$ and set
$\Omega\doteq \R\times{\cal{M}}$. Given a real-analytic function $\Phi:{\cal{M}}\sra \R$, we shall consider the locally
integrable structure $\VV$ on $\Omega$ whose orthogonal $T'$ is spanned by the differential of the function $Z:\Omega\sra \C$,
$Z(x,\qq)=x+i\Phi(\qq)$, $x\in\R$, $q\in {\cal{M}}$. Notice that, if $(t_1,\ldots,t_n)$ are local coordinates on an open subset
$W$ of ${\cal{M}}$, then $\VV$ is spanned, over $\R\times W$, by the vector fields
$$ \LLl_j = \frac{\partial}{\partial t_j}-i\frac{\partial \Phi}{\partial t_j}(t)\frac{\partial}{\partial x}, \quad j=1,\ldots, n.$$
We show that, under these conditions, if $\Phi$ is an open map then $H^{q,n}_{\mc C^\infty}(\Omega;\VV)=0$.
\end{exa}

\begin{proof}
 According to Theorem \ref{thm:open}, the structure $\VV$ is hypocomplex. Moreover, from
[JT, Theorem 1.5], the following property holds: given $(x_0,\qq_0)\in\Omega$,
if $u\in{\cal{D}}'(\Omega;\Lambda^{q,0})$ and $\dd' u\in H^s_{\loc}(\{(x_0,\qq_0\};\Lambda^{q,1})$,
then $u\in H^{s-n/2}_{\loc}(\{x_0,\qq_0\};\Lambda^{p,0})$. The conclusion follows again from Corollary \ref{cor:Main}.
\end{proof}

\begin{exa} Let $\Omega$ be orientable with no compact connected component. Suppose that $m = 1$ and suppose that $\mathcal V$ satisfies:
\begin{enumerate}
    \item $\VV$ is hypocomplex (cf. Theorem \ref{thm:open});
    \item $\VV$ satisfies Condition $(\star)_0$ (as discussed in \cite{cordaro2001local}) at every point.
\end{enumerate}
Then $H^{p,n}(\Omega; \VV) = 0$.
\end{exa}

\begin{exa}
\label{exa:2'}
Let $\Omega$ be an analytic orientable manifold with no compact connected component. Suppose that $\VV$ is analytic and suppose that, for $(\pp, \xi) \in T^0$ with $\xi \neq 0$, we have that $\mc L_{\pp, \xi}$ is non-degenerate and has signature different from $n$ and $n-2$. Then $H^{p,n}(\Omega; \VV) = 0$.
\end{exa}

The following example is a particular case of Example \ref{exa:2'}.

\begin{exa}\label{exa:3}
Let $\mc Q : \C^k \times \C^k \to \C^d$ be a quadratic form and consider the (non compact, connected) quadric submanifold $\Omega = \{(x + iy, w) \in \C^d \times \C^k : y = \mc Q(w, \overline{w})\}$. We endow $\Omega$ with the involutive structure $\mc V$ generated by the vector fields
$$ L_j =  \frac{\del}{\del w_j} + 2i \sum_{l=1}^d \frac{ \del \mc Q_l}{ \del w_j} \frac{\del}{\del z}, \qquad 1 \leq j \leq k.$$ This involutive structure is a CR structure (we refer to \cite{boggess1991cr} to a detailed description of this class of CR manifolds and the computation of its Levi form). Since the characteristic set of $\VV$ is a real line bundle, it suffices to compute it at $(\pp,\xi)$, with $\xi \in T^0_0$.
Up to a sign, this Levi form equals $\mc Q$. If $\mc Q$ is non-degenerate and has at least one positive and one negative eigenvalue, then $\mc V$ is hypocomplex. Furthermore,
 if its signature is different from $n-2$, then the equivalent conditions in Proposition \ref{prp:equiv} hold when $q=1$. Hence, we can infer that $H^{p,n}(\Omega;\VV)=0$ when the signature of $\mc Q$ is different from $n$ and $n-2$. In Chapter 3, we shall see that this example corresponds to a left-invariant hypocomplex CR structure defined in a suitable Lie group.
 \end{exa}

\section{Sheaf cohomology}

In this section, we are going to be assuming that $\mathcal V$ is hypocomplex, but we make no restrictions to the manifold, that is, $\Omega$ can have compact connected components. Let $\mathcal S^p$ be the sheaf $(p,0)$-forms which, are solutions of $\mathcal V$ over $\Omega$, and denote by $\Phi$ a paracompactifying family of supports in $X$. We denote by $H^q_\Phi(\Omega; \mathcal S^p)$ the \emph{cohomology group over $\Omega$ with coefficients in $\mathcal S^p$ and support in $\Phi$}. We are interested in the family of supports consisting of all compacts subsets of $\Omega$ which is going to be denoted by $c$. 

Following \cite{gunning1965analytic}, we recall some notions related to \v Cech cohomology of a covering with values in a sheaf which allows us to use functional analysis to study the operator $\dext'$ and the space $H^1_c(\Omega; \mathcal S^p)$. 

This sheaf cohomology can be related to the cohomology introduced in the beginning of this chapter  via a spectral sequence (cf. \cite{godement1958topologie} Théorème 5.2.4). A particular case of this can be done in locally trivial cohomologies, which is the case of elliptic structures, for example. In fact, if the involutive structure $\mathcal V$ is elliptic, by the Poincaré Lemma for elliptic structures (cf.  \cite{berhanu2008introduction} Theorem VIII.3.1), we know that $\dext'$ is locally exact, so we have a fine resolution of $\mathcal S$. Therefore, we have, for $q > 0$,

$$ H^q(\Omega; \mathcal S^p) = \frac { \ker \left[ \dext' : \mathcal C^\infty(\Omega; \Lambda^{p,q}) \to \mathcal C^\infty(\Omega; \Lambda^{p,q+1}) \right] }{ \img \left[ \dext' : \mathcal C^\infty(\Omega; \Lambda^{p,q-1}) \to \mathcal C^\infty(\Omega; \Lambda^{p,q}) \right] } $$
and
$$ H^0(\Omega; \mathcal S^p) =  \ker \left[ \dext' :\mathcal  C^\infty(\Omega; \Lambda^{p,0}) \to \mathcal  C^\infty(\Omega; \Lambda^{p,1}) \right]. $$

That is, $H^q(\Omega; \mathcal S^p)$ is isomorphic to $H^{p,q}_{\mathcal C^\infty}(\Omega; \mathcal V)$. We note that, in this context, we also have that $\dext' : \mathcal E' (\Omega; \Lambda^{p,q}) \to \mathcal E' (\Omega; \Lambda^{p,q+1})$ has closed range. We are going to show that, under certain conditions, we have similar results in the context of hypocomplex structures. In order to do this, we need some technical results.

We start with some remarks. If $Y \subset \Omega$ is open, the space $S^q(Y)$ of sections over $Y$ is a closed subspace of $\mc C^q(Y)$ and thus an ({\bf{FS}}) space. If $Y'\subset Y$ is also open and with compact closed contained in $Y$, then the restriction map $S^q(Y) \to S^q(Y')$ is compact: this is a consequence of \cite[Corollary III.5.5]{treves1992hypo}. If now $K \ssubset \Omega$ is compact, then $S^q(K)$ denotes the inductive limit of the spaces $S^q(Y)$ when $Y$ runs over the set of all open neighborhoods of $K$. This is a ({\bf{DFS}}) space. If $K'$ is a compact subset of the interior of $K$, then the natural map $S^q(K) \to S^q(K')$ is also compact.

Let $\frk K$ be a cover of $\Omega$ satisfying the following properties:

\begin{itemize}
\item $\frk K$ is countable, with elements indexed as $\frk K = \{K_j\}_{j \in \N}$;
\item $\frk K$ is locally finite;
\item Each $K \in \frk K$ is a compact set of $\Omega$.
\end{itemize}

For $q \in \Z_+$, let $\sigma = (K_{j_0}, \ldots, K_{j_q})$ be a $q$-simplex of $\frk K$, that is, an ordered collection of $q+1$ sets of $\frk K$ such that $|\sigma| \doteq K_{j_0} \cap \ldots \cap K_{j_q} \neq \emptyset.$ Sometimes we denote the simplex $\sigma$ just by the indexes, that is $\sigma = (j_0, \ldots, j_q)$.

We denote by $N^q(\frk K)$ the set of all $q$-simplexes and by $N(\frk K) = \bigcup_{q \in \Z_+} N^q(\frk K)$ the nerv of $\frk K$. 
We recall that for each $q \in \Z_+$, a $q$-cochain of $N(\frk K)$ with values in $\mathcal S^p$ is a function $f$ such that, for each $q$-simplex $\sigma$, there is an associated local section $f_\sigma \in \mathcal S^p(|\sigma|)$. We denote by $C^q(\frk K, \mathcal S^p)$ the set of all $q$-cochains with values in $\mathcal S^p$.

The coboundary operators $\delta^q : C^q(\frk K, \mathcal S^p) \to C^{q+1}(\frk K, \mathcal S^p)$, $q \in \Z_+$, are defined by
$$ (\delta^q f)_\sigma = \sum_{j=0}^{q+1}(-1)^j \rho^{|\sigma_j|}_{|\sigma|} f_{\sigma_j},$$
in which $\sigma_j \doteq (U_0, \ldots, \hat {U_j}, \ldots, U_{q+1} )$ and the mapping $\rho^{|\sigma_j|}_{|\sigma|}$ associates each $p$-form in $|\sigma_j|$ to its restriction to $|\sigma|$. By definition, we have that $\delta^{q+1} \circ \delta^{q} = 0$ for all $q \in \Z_+.$

We define the cohomology of the complex $\delta^q : C^q(\frk K, \mathcal S^p) \to C^{q+1}(\frk K, \mathcal S^p)$, $q \in \Z_+$, for $q \geq 1$, as

\begin{equation}
\label{eq:chech_coho}
\check H^q(\frk K, \mathcal S^p) = \frac{ \ker \left[ \delta^q : C^q(\frk K, \mathcal S^p) \to C^{q+1}(\frk K, \mathcal S^p)\right] }{ \img \left[ \delta^{q-1} : C^{q-1}(\frk K, \mathcal S^p) \to C^{q}(\frk K, \mathcal S^p)\right] },
\end{equation} and, for $q = 0$, as
\begin{equation}
\label{eq:chech_coho_}
\check H^0(\frk K, \mathcal S^p) = \ker \left[ \delta^0 : C^0(\frk K, \mathcal S^p) \to C^{1}(\frk K, \mathcal S^p)\right] = \mathcal S^p(\Omega).
\end{equation} We defined in \eqref{eq:chech_coho} and \eqref{eq:chech_coho_} the \emph{\v Cech cohomology of $\frk K$ with values in $\mathcal S^p$}.

We denote by $\mathcal C^q_c(\frk K, \mathcal S^p)$ the set of all \emph{finite} $q$-cochains with values in $\mathcal S^p$, that is, for a given $f \in \mathcal C^q_c(\frk K, \mathcal S^p)$, there is only a finite number of $q$-simplexes of $\frk K$ such that $f_\sigma$ is non-zero. We clearly have $\delta^q( C^q_c(\frk K, \mathcal S^p) ) \subset C^{q+1}_c(\frk K, \mathcal S^p)$ for all $q \in \Z_+$, and so we have a complex whose cohomology is denoted by $\check H^q_c(\frk K, \mathcal S^p)$ for all $q$. And now we have what we call \emph{\v Cech cohomology of $\frk K$ with values in $\mathcal S^p$ with compact support}.

We note that $C^q_c(\frk K; \mathcal S^p)$ is a closed subspace of a direct sum of ({\bf DFS}) spaces and so it is also a ({\bf DFS}) space. Moreover, since $\delta^q$ is a linear continuous map between ({\bf DFS}) spaces, each $\ker \delta^q$ is also a ({\bf DFS}) space.

Let $A \subset \N$ be a finite subset. We denote by $C^q_A(\frk K; \mathcal S^p)$ the set of all cochains $s \in C^q_c(\frk K; \mathcal S^p)$ such that $s_{j_0 j_1 \ldots j_q} = 0$ if $\{j_0, j_1, \ldots, j_q\} \cap A = \emptyset.$

Since $\delta^q( C^q_A(\frk K; \mathcal S^p)) \subset C^{q+1}_A(\frk K; \mathcal S^p)$ we have a well defined \emph{complex of cochains with supports contained in $A$}. The corresponding cohomology spaces are denoted by $H^q_A(\frk K; \mathcal S^p)$ with $q = 0, 1, \ldots$. We have natural homomorphisms
$$ H^q_A(\frk K; \mathcal S^p) \to H^q_B(\frk K; \mathcal S^p) $$ if $A \subset B$ and clearly $$ H^q_A(\frk K; \mathcal S^p) \to H^q_c(\frk K; \mathcal S^p).$$

\begin{prp}
\label{hypo:aprime} Suppose that $\Omega$ is endowed with a hypocomplex involutive structure $\mathcal V$. For every finite subset $A$ of $\N$, there is another finite subset $A'$ of $\N$ with $A \subset A'$ such that the kernel of the homomorphism $H^1_A(\frk K; \mathcal S^p) \to H^1_c(\frk K; \mathcal S^p)$ is contained in the kernel of the homomorphism $H^1_A(\frk K; \mathcal S^p) \to H^1_{A'}(\frk K; \mathcal S^p)$.
\end{prp}

\begin{proof}
Let $L \subset \Omega$ be a compact such that $Y \doteq \Omega \bs L$ is connected and $\bigcup_{j \in A} K_j \subset L$. We define
$ A' = \{ j \in \N : K_j \cap L \neq \emptyset \}.$
Since $\frk K$ is locally finite, we have that $A'$ is finite. We also have, by the choice of $L$, that $A \subset A'$.

Let $s \in C^1_A(\frk K; \mathcal S^p)$ be such that $\delta^1 u = s$ for some $u \in C^0_c(\frk K; \mathcal S^p)$.
We write $u = \{u_j\}$ with $u_j \in \mathcal S^p(K_j)$.
We have $u_j - u_k = s_{jk}$ as elements in $\mathcal S(K_j \cap K_k)$ if $K_j \cap K_k \neq \emptyset.$ Let $p \in Y$ and assume that $p \in K_j \cap K_k$.
Notice that $\{j,k\} \cap A = \emptyset$  implies $s_{jk} = 0$ and consequently the germs of $u_j$ and $u_k$ coincide at $p$. Since $\VV$ is hypocomplex, we have that there is $f \in \mathcal S^p(Y)$ such that the germs of $f$ and $u_j$ coincide at $p$ if $p \in K_j$. Since all but a finite number of the $u_j$'s are zero, it follows that $f$ vanishes in an open neighborhood of $\del \Omega$ and thus in $Y$.
In particular, $u_j = 0$ if $K_j \subset Y$, that is, $u_j =0$ if $j \notin A'$. 
\end{proof}

\begin{prp}
\label{hypo:rho1_injective} Let $\frk L = \{L_j\}$ be a countable and locally finite refinement of $\frk K$ with each $L_j$ compact. Then, the refinement homomorphism
$$ \rho^q : C^q_c(\frk K; \mathcal S^p) \to C^q_c(\frk L; \mathcal S^p), \qquad q = 0, 1, 2 \ldots $$
induces homomorphisms
$$ \rho^q_* : H^q_c(\frk K; \mathcal S^p) \to H^q_c(\frk L; \mathcal S^p), \qquad q = 0, 1, 2 \ldots $$
and $\rho^1_*$ is injective.
\end{prp}


\begin{proof}
Let $ \gamma : \N \to \N$ be such that $L_j \subset K_{\gamma(j)}$ for every $j$. Let $ s = \{s_{jk}\} \in C^1_c(\frk K; \mathcal S^p)$ satisfying $\delta^1 s = 0$ and suppose that $\rho^1 s = \delta^0 v$ for some $v = \{v_j\} \in C^0_c(\frk L; \mathcal S^p).$ Then
$ v_k - v_j = s_{\gamma(j) \gamma(k)}$
as elements of $ \mathcal S^p(L_j \cap L_k) $ and, since $s$ is a cocycle, we have that $s_{kr} - s_{jr} + s_{jk} = 0$ as elements of $\mathcal S^p(K_j \cap K_k \cap K_r)$ and consequently
$$ v_k + s_{\gamma(k)r} = v_j + s_{\gamma(j)r}$$
as elements of $\mathcal S^p(L_j \cap L_k \cap K_r)$. Notice that the left hand side does not depend on $j$ and the right hand side does not depend on $k$. By varying $j$ and $k$, we see that there is $u_r \in \mathcal S^p(K_r)$ such that $u_r = v_j + s_{\gamma(j)r}$ as elements of $\mathcal S^p(L_j \cap K_r)$ whenever $K_j \cap K_r \neq \emptyset.$

Let $I,J \subset \N$ be such that $v_j = 0$ unless $j \in I$ and $s_{jk} = 0$ unless $j$ and $k$ belong to $J$. The set $I' \doteq J \cup \left\{ r \in \N : \left( \bigcup_{j \in I} L_j \right) \cap K_r \neq \emptyset\right\}$
is finite and $u_r = 0$ if $r \notin I'$. Indeed, if $r \notin I'$, then $r \notin J$ and, if $j$ is such that $L_j \cap K_r \neq \emptyset$, then $j \notin I$. These two facts imply $u_r = 0$ in $K_r \cap L_j $. Thus $u \doteq \{u_j\}$ belongs to $C^0_c(\frk K; \mathcal S^p)$. Finally, as elements of $\mathcal S^p(L_j \cap K_r \cap K_t)$, we have
$$ u_r - u_t = (v_j + s_{\gamma(j)r}) - (v_j + s_{\gamma(j)t}) = s_{\gamma(j)r} - s_{\gamma(j)t} = - s_{tr}.$$ And the proof is completed. \end{proof}

From now, on we assume that $\frk K$ and $\frk L$ satisfy the following properties:

\begin{equation}
\label{hypo:p1} L_j \subset \operatorname{int}(K_j) \text{ for all } j \in \N.
\end{equation}
\begin{equation}
\label{hypo:p2}
\text{\parbox{.85\textwidth}{ If $(j_0, j_1)$ is a 1-simplex for $\frk K$, then it is also a 1-simplex for $\frk L$. And if $(j_0, j_1, j_2)$ is a 2-simplex for $\frk K$, then it is also 2-simplex for $\frk L$.}}
\end{equation}

For this, it suffices to take $\frk L$ to be a locally finite countable covering of $\Omega$ such that each $K_j$ is compact and $\frk K = \{K_j\}$ is the covering defined by
$$ K_j = \{ p \in \Omega : \operatorname{dist}(p, L_j) \leq \epsilon_j\}.$$

It is clear that $\frk K$ also is locally finite and countable and that each $K_j$ is compact if $\epsilon_j > 0$ is chosen suitable small.



Notice that Property \eqref{hypo:p2} guarantees that the maps $\rho^p$ are injective for $p = 0, 1, 2$. Moreover, if $A \subset \N$ is finite, then $A$ induces homomorphisms
$$ \rho^p_A : C^p_A(\frk K; \mathcal S^p) \to C^p_A(\frk L; \mathcal S^p)$$
which are now, thanks to Property \eqref{hypo:p1}, compact operators.

\begin{thm}
\label{hypo:thm:dfstop} Assume that $\Omega$ is endowed with a hypocomplex structure. The map $$ \delta^0 : C^0_c(\frk K; \mathcal S^p) \to C^1_c(\frk K; \mathcal S^p) $$ has closed range. 
\end{thm}

\begin{proof} By Theorem 6' of \cite{komatsu1967projective}, it suffices to show that $\img(\delta^0) \cap C^1_A(\frk K; \mathcal S^p)$ is closed in $C^1_A(\frk K; \mathcal S^p)$ for all finite subsets $A \subset \N.$

Let $A \subset \N$ be finite and let $A'$ be as in Proposition \ref{hypo:aprime}. We define the following ({\bf DFS}) spaces:
$$ V \doteq \{ (u,s) \in C^0_{A'}(\frk K; \mathcal S^p) \times C^1_{A}(\frk K; \mathcal S^p): \delta^0 u = s\},$$
$$ W \doteq \{ (u,s) \in C^0_{A'}(\frk L; \mathcal S^p) \times C^1_{A}(\frk K; \mathcal S^p): \delta^0 u = \rho^1s\},$$
and let $\sigma : V \to W$ be the linear map given by $ \sigma(u, s) = (\rho^0_{A'}u, s).$

We need the following lemma:

\begin{lem}  Assume that $\Omega$ is endowed with a hypocomplex structure. The map $\sigma$ is a bijection.
\end{lem}
\begin{proof}
Since $\delta^0$ is injective, we have that $\sigma$ is injective. In fact, if $\sigma(u,s) = (0,0)$, then $s = 0$ and, since $(u,s) \in V$ we have $\delta^0u = 0$.

Now we are going to prove that $\sigma$ is surjective. Let $(v,s) \in W$. We have that $\delta^1(\rho^1 s) = \delta^1(\delta^0 v) = 0$ and, since $\rho^2(\delta^1 s) = \delta^1(\rho^1 s) = 0 $ and $\rho^2$ is injective, it follows that $\delta^1 s = 0$. Hence, the class of $s$ defines an element in $H^1_A(\frk K; \mathcal S^p)$ which belongs to the kernel of $\rho^1_*$.
By Proposition \ref{hypo:rho1_injective}, the class of $s$ vanishes in $H^1_c(\frk K; \mathcal S^p)$ and so Proposition \ref{hypo:aprime} implies that the class of $s$ vanishes in $H^1_{A'}(\frk K; \mathcal S^p)$. In other words, there is $u \in C^1_{A'}(\frk K; \mathcal S^p)$ such that $(u,s) \in V$.
To finish the proof, we need to show that $\sigma(u,s) = (v,s)$. This follows from the facts that $\delta^0(v - \rho^0_{A'}u) = 0$ and that $\delta^0$ is injective.
\end{proof}
We denote by $E$ the ({\bf DFS}) space $C^0_{A'}(\frk L; \mathcal S^p) \times C^1_A(\frk K; \mathcal S^p)$; we remark that $W$ is a closed subspace of $E$ and thus also a ({\bf DFS}) space. Let $\imath$ be the inclusion map $W \subset E$. It follows that $\imath \circ \sigma : V \to E$ is a topological embedding.

We introduce the following continuous linear maps
 $$ \tau_1(u,s) = (\rho^0_{A'}, 0), \qquad \tau_2(u,s) = (0,s).$$

Since $\tau_1$ is compact and $\tau_2 = \imath \circ \sigma - \tau_1$, it follows from Theorem 1 of \cite{schwartz1956} that $\tau_2$ has closed range in $E$. This is the same as saying that the map $p_2 : V \to C^1_A(\frk K; \mathcal S^p)$, given by the projection on the second factor, has closed range, but according to Proposition \ref{hypo:aprime} we have
$$ \img{p_2} = \img(\delta^2) \cap C^1_A(\frk K; \mathcal S^p),$$
which completes the proof.
\end{proof}

\begin{thm}
\label{thm:closed_range} Let $\Omega$ be orientable and endowed with a hypocomplex structure $\mathcal V$. Assume that the equivalent conditions in Proposition \ref{prp:equiv} hold for every point of $\Omega$. Then 
\begin{equation}
\label{hypo:eq:dprimecpctsupp} \dext' : \mathcal C^\infty (\Omega; \Lambda^{p,n-1}) \to \mathcal C^\infty (\Omega; \Lambda^{p,n})
\end{equation}
has closed range.
\end{thm}

We emphasize that we are not precluding that $\Omega$ has compact connected components.

\begin{proof}
According to \cite[Page 18]{kothe19679topological}, it suffices to show that $\dext' : \mathcal E' (\Omega; \Lambda^{m-p,0}) \to \mathcal E' (\Omega; \Lambda^{m-p,1})$ has strongly sequentially closed range. The argument follows by the conjunction of Remark \ref{rmk:1}, the proof of Theorem \ref{thm:first}, and Theorem \ref{hypo:thm:dfstop}.

Let $\{u_j\} \subset \mathcal E'(\omega; \Lambda^{m-p,0})$ be such that $\dext' u_j \to f$ in $E'(\omega; \Lambda^{m-p,1})$. Hence, there are a compact set $K \subset \Omega$ and a real number $s \in \R$ such that $\dext' u_j \in H^s(K, \Lambda^{m-p,1})$, for all $j$.

By the hypothesis and by Remark \ref{rmk:1}, there exists a covering $\frk U = \{K_j\}$ of $\Omega$ by compact sets such that $\{K_1, \ldots, K_\nu\}$ is a covering of $K$ and there exists an $s \in \R$ such that the following is true:

For every $l=1, \ldots, \nu$, there are $U_l \supset K_l$ open and $v_{l,j} \in H^t_{\tiny{loc}}(U_l; \Lambda^{m-p,0})$ such that $\dext' v_{l,j} = \dext' u_j$ in $U_l$ and $\{v_{l,j}\}_{j=1}^\infty$ is a Cauchy sequence if $H^t_{\tiny{loc}}(U_l; \Lambda^{m-p,0})$ (this last property follows from the proof of Proposition \ref{prp:equiv}).
We define $v_l = \lim_j v_{l,j} \in H^t_{\tiny{loc}}(U_l; \Lambda^{m-p,0}).$ For $l > \nu$ we have $\dext' u_j = 0$ for all $j$ and so we take $v_{l,j} = 0$.

Now, suppose that $K_{l_1} \cap K_{l_2} \neq \emptyset$. Then $\dext' (v_{l_1,j} - v_{l_2,j}) = 0$ in a neighborhood of $K_{l_1} \cap K_{l_2}$ and so $V_j = \{v_{l_1,j} - v_{l_2,j}\}$ defines, for each $j$, an element in $C^1(\frk U; \mathcal S^{m-p})$.

Moreover, a standard argument of comparison of topologies, taking into account the property of hypocomplexity, shows that $V_j \to V$ in $C^1(\frk U; \mathcal S^{m-p})$. We have that $V = \{v_{l_1} - v_{l_2}\}_{l_1, l_2}$. On the other hand, $W_j = \{v_{lj} - u_j\}_l \in C^0(\frk U; \mathcal S^{m-p})$ and $\delta_0 W_j = V_j \to V$.

By Theorem \ref{hypo:thm:dfstop}, there is $W \in C^0(\frk U; \mathcal S^{m-p})$ such that $\delta_0 W = V$. Now we write $W = \{w_l\} \in C^0(\frk U; \mathcal S^{m-p})$, and then we obtained
$w_{l_1} - w_{l_2} = v_{l_1} - v_{l_2}$ near $K_{l_1} \cap K_{l_2}$. Now we define $u$ by setting $u = v_l - w_l$ near $K_l$ and we have $d'u = d'v_l = f$ near $K_l$. By construction, $u$ has compact support and the proof is completed.
\end{proof}

\section{Serre duality}

We recall some basic facts about Fréchet spaces. We remember that, if $E$ and $F$ are Fréchet spaces and $u : E \to F$ is a linear continuous map, we say that $u$ is a homomorphism if $u(E)$ is closed. By the Open Mapping Theorem, we have that $u$ is a homomorphism if and only if $u': E / \ker u \to F$ is an isomorphism.

The following result is important for us:

\begin{prp}
\label{prp:finite_dim_and_closedness} Let $E$ and $F$ be Fréchet spaces and $u : E \to F$ be a linear continuous map. If the codimension of $u(E)$ is finite, then $u$ is a homomorphism.
\end{prp}

\begin{prp} Let $E$, $F$ and $G$ be Fréchet spaces. Consider the mappings
$$ E \xrightarrow{u} F \xrightarrow{v} G .$$ We assume that $u$ and $v$ are homomorphism and that $u(E) \subset \ker v$. Then $ \ker v / u(E) $ is a Fréchet space and we have a canonical isomorphism
$$ \left( \frac{\ker v}{u(E)} \right)^* \cong \frac{ \ker u^t}{v^t(G^*)} .$$
\end{prp}

Combining the results above with the discussion that precedes Theorem 2.2.3, we obtain the following version of Serre's duality theorem:

\begin{thm}
\label{chp:inv:thm:serre_duality}
Suppose that $\Omega$ is orientable and that the following linear maps are homomorphisms
$$ \mathcal C^\infty(\Omega; \Lambda^{p,q-1}) \xrightarrow{\dext'} \mathcal C^\infty(\Omega; \Lambda^{p,q}) \xrightarrow{\dext'} \mathcal C^\infty(\Omega; \Lambda^{p,q+1}).$$
Then $ H^{p,q}_{\mathcal C^\infty}(\Omega;\mathcal V)$ is a Fréchet space and there is a canonical isomorphism
$$  H^{p,q}_{\mathcal C^\infty}(\Omega; \mathcal V)^* \cong  H^{m-p,n-q}_{\mathcal E'}(\Omega; \mathcal V).$$
\end{thm}


\section{Hodge decomposition}
\label{sec:hodge_theory}

In this section, we are going to explain how to use the Hodge decomposition for the operator $\dext'$. With this decomposition, we are going to find some simple conditions so that the cohomology is finite dimensional on a specific degree. Usually, Hodge decomposition is introduced in the context of elliptic differential operators between Hermitian vector bundles. Here we made a small improvement: we assumed that the operator is subelliptic. See Definition \ref{dfn:subellipticity_and_stuff}.

Suppose that $\Omega$ is a compact manifold with a positive $\mathcal C^\infty$-density $\mu$ and let $E_j$, for $j = 1, 2, 3$ be three Hermitian vector bundles over $\Omega$. We denote by $\ip{,}_j$ the hermitian inner product on $E_j$.
Let $$A_j : \mathcal C^\infty(\Omega; E_j) \to \mathcal C^\infty(\Omega; E_{j+1}), \qquad j = 1, 2$$ be partial differential operators of order $m$ satisfying $A_2 \circ A_1 = 0$. That is, we have
$$ \mathcal C^\infty(\Omega; E_1) \xrightarrow{A_1}{} \mathcal C^\infty(\Omega; E_{2})  \xrightarrow{A_2}{}  \mathcal C^\infty(\Omega; E_3).$$

We endow $\mathcal C^\infty(\Omega; E_j)$ with the inner product
$$\iip{u,v}_j \doteq \int_\Omega \ip{u(x),v(x)}_j \dint \mu(x),$$
for $u, v \in \mathcal C^\infty(\Omega; E_j)$, and we define the operators $A_j^* : \mathcal C^\infty(\Omega; E_{j+1}) \to \mathcal C^\infty(\Omega; E_{j})$ by
$$ \iip{A_j u, v }_{j+1} = \iip{u, A^*_j v}_{j+1},$$
for $u \in \mathcal C^\infty(\Omega; E_{j})$ and $v \in \mathcal C^\infty(\Omega; E_{j+1}).$

Now we have everything we need to introduce the Hodge operator:

\begin{dfn}
For $u \in \mathcal C^\infty(\Omega; E_2)$, we define $\Box u \doteq (A_1 A_1^* + A_2^* A_2) u$ and, in the sense of distributions, we define $$D(\Box) = \{u \in L^2(\Omega; E_2): \Box u \in L^2(\Omega; E_2) \}.$$ The operator $$\Box : D(\Box) \subset L^2(\Omega; E_2) \to L^2(\Omega; E_2)$$ is densely defined and is called \emph{Hodge operator}, \emph{Hodge Laplacian} or just \emph{Box operator}.
\end{dfn}


The operator $\Box$ has closed graph. In fact, let $u_j, v_j \in L^2(\Omega; E_2)$ be a sequence satisfying $\Box u_j = v_j$ and suppose that $u_j \to u \in L^2(\Omega; E_2)$ and $v_j \to v \in L^2(\Omega; E_2)$.
We want to show that $u \in D(\Box)$ and $\Box u = v$. Since $u_j \to u$ and $v_j \to v$ in $L^2(\Omega; E_2)$, we have this convergence in $\mc D'(\Omega; E_2)$ and, since $\Box$ is a differential operator, it is continuous in $\mc D'(\Omega; E_2)$. Thus, we can apply the Closed Graph Theorem and conclude that $\Box u = v$ in the sense of distributions. By hypothesis, we have $u,v \in L^2(\Omega; E_2)$ and so $u \in D(\Box)$.

We also remark that
$$
\begin{aligned}
\iip{\Box u, u}
	& = \iip{A_1 A_1^*u, u} + \iip{A_2^* A_2u, u} \\
	& = \iip{A_1^*u, A_1^*u} + \iip{A_2 u, A_2 u} \\
	& = \|A_1^*u\|_1^2 + \|A_2 u\|_3^2
\end{aligned}
$$
and so we have the following result:

\begin{prp}
\label{prp:kernel_Box_decomposition} The kernel of the operator $\Box$ is
$$ \ker \Box = \{u \in L^2(\Omega; E_2): A_1^*u = 0 \text{ and } A_2 u = 0\}. $$
\end{prp}

In the following definition, we introduce two concepts regarding the regularity of differential operators. Under certain conditions, these concepts provide sufficient conditions so that the kernel of the differential operator has finite dimension.

\begin{dfn}
\label{dfn:subellipticity_and_stuff} Let $E$ be a vector bundle over $\Omega$ and consider a linear partial differential operator $P : \mc D' (\Omega; E) \to \mc D'(\Omega; E)$. For an $\epsilon > 0$, the operator $P$ is said to be \emph{$\epsilon$-subelliptic} if $Pu \in H^s(\Omega;E)$ implies that $u \in H^{s+\epsilon}(\Omega; E)$ for all $u \in H^s(\Omega; E)$.
The operator $P$ is said to be \emph{globally hypoelliptic} if, given $u \in \mc D'(\Omega; E)$ such that $Pu \in \mathcal C^\infty(\Omega;F)$, it holds that $u \in \mathcal C^\infty(\Omega;E)$.
Notice that subellipticity implies global hypoellipticity.
\end{dfn}

On the above definition, it is not necessary to assume that the manifold $\Omega$ is compact.

The following proposition is a well known result about differential operators on compact manifolds.

\begin{prp} If $P : \mc D' (\Omega; E) \to \mc D'(\Omega; F)$ is globally hypoelliptic, then $$\dim \ker P < \infty.$$
\end{prp}

\begin{proof}
Since $P$ is globally hypoelliptic, we have that $\ker P \subset \mathcal C^\infty(\Omega; E)$. Consider $\ker P$ with the topology induced by $L^2(\Omega; E)$. Notice that this induced topology agrees with the topology induced by $H^s(\Omega; E)$ for any $s \in \R$. Take $B = \{ u \in \ker P: \|u\|_1 \leq 1 \}$ and note that $B$ is closed and bounded. By the Rellich lemma, we know that the inclusion $H^1(\Omega; E) \to H^0(\Omega; E) = L^2(\Omega; E)$ is compact and so $B$ is compact.
\end{proof}

In this work, we are interested in the case in which $P$ is the Hodge operator. Therefore, we need to establish some conditions so that the Hodge operator is globally hypoelliptic. But first we need to finish our decomposition. We need the next lemma which is a basic result on Hilbert spaces.

%

This next proposition shows another important result regarding subellipticity of the Hodge operator. By itself, this proposition is interesting, but it is also necessary for the decomposition.

\begin{prp} If the operator $\Box$ is $\epsilon$-subelliptic, then there exists $C > 0$  such that
\begin{equation}
    \|u\|_0 \leq C \|\Box u\|_0, \qquad u \in (\ker \Box)^\perp \cap D(\Box).
\end{equation}
\end{prp}

\begin{proof} 
Defining the norm $$ u \mapsto \left\{ \|u\|_0^2 + \|\Box u\|_0^2 \right\}^{1/2} $$ on $$ D(\Box) = \{ u \in L^2(\Box; E_2) : \Box u \in L^2(\Box; E_2) \}$$ gives a Hilbert space and subellipticity gives $D(\Box) \subset H^\epsilon(\Omega; E_2)$. By the Closed Graph Theorem, there exists a constant $C > 0$ such that
$$ \|u\|_\epsilon \leq C (\|\Box u\|_0 + \|u\|_0),$$
for all $u \in D(\Box)$. Hence, if $u \in (\ker \Box)^\perp \cap D(\Box)$, we have $$\|u\|_\epsilon \leq C \|\Box u\|_0.$$
Indeed, suppose that this last inequality is not valid. Then, for every $j$, there exist a $u_j$ such that $\|u\|_\epsilon \geq j\|\Box u _j\|_0$. Notice that, if we take $v_j = u_j/\|u\|_\epsilon\|$, we have $\|v_j\|_\epsilon = 1$ and $\|\Box v_j\| < 1/j$. There exists a subsequence $\{v_{j_k}\}$ such that $v_{j_k} \to v$ in $L^2(\Omega; E_2)$ and clearly $v \in \ker \Box \cap (\ker \Box)^\perp$. Notice that $$ 1 \leq C (\|\Box v_{j_k}\|_0 + \|v_{j_k}\|_0)$$ and so
$$ 1 \leq C \|v\|_0$$
which is a contradiction. The proof is completed.
\end{proof}

A standard argument, making use of 2.10 (see, for instance, Lemma 4.1.1 in \cite{hormander1973introduction}), shows the existence of a bounded linear operator $G : (\ker \Box)^\perp \to (\ker \Box)^\perp$ with $\|G\| \leq C$ such that $\Box G f = f$ for $f \in \ker \Box.$ Now, given $g \in L^2(\Omega; E_2)$, we can write $g = f + h$ with $h \in \ker \Box$ and $f \in (\ker \Box)^\perp$. Hence $g = \Box G f + h$.
Moreover, if $g \in \mc C^\infty(\Omega; E_2)$, both $f$ and $h$ also are in $\mc C^\infty(\Omega; E_2)$, and also $Gf \in \mc C^\infty(\Omega; E_2)$ because $\Box$ is hypoelliptic.

We have proved:

\begin{prp} If $\Box$ is $\epsilon$-subelliptic, then
we have the Hodge decomposition
\begin{equation}
	\label{eq:hodge_decomposition}
	\mathcal C^\infty(\Omega;E_2) = \ker \Box \oplus \Box(\mathcal  C^\infty(\Omega;E_2) ) = \ker \Box \oplus A_1 ( \mathcal C^\infty(\Omega;E_1) ) \oplus A_2^* (\mathcal C^\infty(\Omega;E_3))
\end{equation}
and
\begin{equation}
	\label{eq:hodge_decomposition_isomorphism}
	\frac{ \ker A_2 }{ A_1(\mathcal C^\infty(\Omega; E_2))} \cong \ker \Box.
\end{equation}

\end{prp}

\chapter{Involutive structures on compact Lie groups}
\label{chp:involutive_structures_on_compact_lie_groups}

In this chapter, we are going to restrict our attention to involutive structures on compact Lie groups. We introduce the concept of left-invariant involutive structure, which is an involutive structure that encodes some algebraic properties of the Lie group.

Next, we prove that there is an one-to-one correspondence between left-invariant involutive structures on Lie groups and subalgebras of the complexification of the Lie algebra. Also, we show that there are many left-invariant structures with interesting analytical properties.

We want to show that the study of the cohomology of some of these left-invariant involutive structures can be done using only the Lie algebras. In order to explain exactly how to do that, we need to introduce the concepts of left-invariant cohomology of Lie algebras and left-invaraint cohomology of Lie algebras induced by subalgebras.

Then we explain how to use some of the Chevalley and Eilenberg's techniques to study the de Rham cohomology. These techniques can be partially adapted to other left-invariant involutive structures and with them we show that it is possible to include the left-invariant cohomology relative to a subalgebra into the usual cohomology. This inclusion shows that there are algebraic obstructions to solvability. Then we discuss some necessary conditions so that these techniques can be used in general.

Finally, we start using Hodge decomposition and its relation to Lie derivatives to obtain some new results regarding left-invariant cohomologies, namely Theorem \ref{thm:comp_lie:hypo_subelliptic} and Theorem \ref{thm:comp_lie:hypo_subelliptic_torus_2}. Then, we show one application of Serre duality, namely Theorem \ref{thm:left_invariance_on_degree}, which is the main theorem in this chapter.

\section{Left-invariant involutive structures}
\label{sec:basicdef}

Let $G$ be a Lie group with a Lie algebra $\frk g_\R$, that is, the set of all left-invariant vector fields of $G$. We denote by $L_x : G \to G$ the left-multiplication, that is, $L_x(g) = xg$ for $g \in G$.

\begin{dfn}
A vector bundle $\mathcal V \subset \C T G$ is called \emph{left-invariant} if $$(L_x)_*X_g \in \mathcal V_{xg}$$ for all $X_g \in \mathcal V_g$.
\end{dfn}

The next lemma shows how to construct an \emph{involutive} left-invariant subbundle of $\C T G$ using a Lie subalgebra of $\frk g \doteq \frk g_\R \otimes \C.$

\begin{lem}
Let $G$ be a Lie group and let $\frk h \subset \frk g$ be any complex Lie subalgebra. Then $\frk h$ defines an involutive vector subbundle of $\C TG$, which is invariant by the left action of $G$ on itself.
\end{lem}

\begin{proof}
Let $\frk h_x = \{X(x) \in \C T_x G: X \in \frk h \}$ and define $\mathcal V_{\frk h} \doteq \bigcup_{x \in G} \frk h_x$. Notice that $\frk h_x = (L_x)_* \frk h_e$ and, since for each $x$ the map $(L_x)_*$ is an isomorphism of vector spaces, we have that $\mathcal V_{\frk h}$ is a complex vector bundle and $\mathcal V_{\frk h} \cong G \times \frk h_e$.
Now we are going to verify that $\mathcal V_{\frk h}$ is involutive. Notice that any left-invariant vector field of $\frk h$ is a smooth section of $\mathcal V_{\frk h}$. Let $X_1, \ldots, X_k$ be a vector basis for $\frk h$. Now $X_1(g), \ldots, X_k(g)$ is a vector basis for $\frk h_g$ for all $g$. Let $X, Y$ be smooth sections of $\mathcal V_{\frk h}$. We write $X = \sum_i x_i X_i$ and $Y = \sum_j y_j X_j$, in which $x_i, y_j \in \mathcal C^\infty(G).$
Notice that
$$[X,Y] = \sum_j \left( \sum_i x_i X_i(y_j) \right) X_j - \sum_i \left( \sum_j y_j X_j(x_i) \right) X_i + \sum_i \sum_j x_i y_j [X_i,X_j].$$

The first two summations clearly are smooth sections of $\mathcal V_{\frk h}$ and, by using that $[X_i,X_j]$ is a linear combination of $X_1, \ldots, X_k$, it follows that the third summation also is.
\end{proof}

As a corollary, we have that each left-invariant involutive structure $\mathcal V \subset \C T G$ defines a complex Lie subalgebra $\frk h \subset \frk g$. We prove this fact as follows:

\begin{cor} Let $G$ be a Lie group and let $\mathcal V$ be an involutive vector subbundle of $\C TG$, which is left-invariant. Then $\mathcal V$ defines a complex subalgebra of $\frk g.$
\end{cor}

\begin{proof}
Let $X_1(e), \ldots, X_k(e)$ be a vector basis for $\mathcal V_e$ and define $X_j(g) = (L_g)_*X_j(e)$. The vector bundle $\mathcal V$ is left-invariant, so we have $X_j(g) \in \mathcal V_g$ and each $X_j$ is a left-invariant section of $\mathcal V$. We then define $\frk h = \operatorname{span}_\C \{X_1, \ldots, X_k\}$.
\end{proof}

Therefore, on a Lie group $G$, we have that any complex Lie subalgebra $\frk h \subset \frk g$ is in an one-to-one correspondence with left-invariant involutive structures on $G$. We usually denote the vector bundle $\mathcal V_{\frk h}$ by its corresponding Lie algebra $\frk h$.

Notice that, since every Lie group is an analytic manifold, every involutive structure over $G$ is a \emph{locally integrable structure} \citep{berhanu2008introduction, treves1992hypo}.

We are interested in studying properties of the space $H^{p,q}_{\mathcal C^\infty}(G; \frk h)$. Inspired by the local theory, we believe that the properties of this space depend on how the Lie algebra $\frk h$ is included on the algebra $\frk g$. However, we also expect to take into account some intrinsic algebraic properties of $\frk h$ and topological aspects of $G$.

By using the language of Lie algebras, some special locally integrable structures over $G$ can be highlighted.

\begin{dfn}
\label{dfn:basic_types_of_lie_algebras}
We shall say that a Lie subalgebra $\frk h \subset \frk g$ defines:

\begin{itemize}
\item an \emph{elliptic structure} if $\frk h + \bar {\frk h} = \frk g$. The subalgebra $\frk h$ is called an \emph{elliptic subalgebra};
\item a \emph{complex structure} if $\frk h \oplus \bar {\frk h} = \frk g$. The subalgebra $\frk h$ is called a \emph{complex subalgebra};
\item a \emph{Cauchy-Riemann (CR) structure} if $\frk h \cap \bar {\frk h} = 0$. The subalgebra $\frk h$ is called a \emph{CR subalgebra};
\item an \emph{essentially real structure} if $\frk h = \bar {\frk h}$. The subalgebra $\frk h$ is called a \emph{essentially real subalgebra}.

\end{itemize}
\end{dfn}

For example, when $\frk h = \frk g$, clearly the complex \eqref{h_complex} is the de Rham complex. In this case, notice that we have an elliptic structure that is also an essentially real structure. When $\frk h \oplus \overline{\frk h} = \frk g$, we have a complex structure that is also an elliptic structure over $G$. In this case, the operator $\dext'$ is the $\bar \del$ operator and the associated complex \eqref{h_complex} is the Dolbeault complex. Another example showing an interesting behavior is the following:


\begin{exa}
\label{exa:torus_liouvile} Consider the two torus $\mathbb{T}^2$ with coordinates $(x,y)$. Its Lie algebra is given by $\frk t = \operatorname{span}_\C \{ \del/\del x,  \del/\del y\} $. Let $\mu$ be a real number, take $L = \del/\del x - \mu \del/\del y$, and define $\frk h = \operatorname{span}_\C \{ L \}$.




It is easy to see that $$ \mathcal C^\infty(\mathbb{T}^2, \Lambda^{0,1}) \cong \{ f \dext x : f \in \mathcal C^\infty(\mathbb{T}^2)  \} $$
and
$$ \mathcal C^\infty(\mathbb{T}^2, \Lambda^{0,2}) \cong \{0\}.$$

By writing $M = \del / \del y$, we have that $L$ and $M$ form a basis for $\frk t$ and the respective dual basis consists of $\zeta = \dext x$ and $\tau = \mu \dext x + \dext y$. Therefore, the operator $\dext$ can be written as
$\dext u = M(u) \zeta + L(u) \dext x$ and so $\dext'u = (Lu)\dext x$. Thus, we have the following complex
$$\mathcal  C^\infty(\mathbb{T}^2) \xrightarrow[]{ \dext'} \{ f \dext x : f \in \mathcal C^\infty(\mathbb{T}^2)  \} \xrightarrow[]{ \dext'} \{0\}.$$

Let $v \in \mathcal C^\infty(\mathbb{T}^2, \Lambda^{0,1})$ and write $v = f \dext x$, in which $f \in \mathcal C^\infty(\mathbb{T}^2)$. Notice that, if we want to show that $v$ is $\dext'$-exact, we have to find a function $u$ satisfying $\dext' u = v$. This is equivalent to solve the partial differential equation $Lu = f$. It is well known that this equation is globally solvable in $\mathcal C^\infty(\mathbb{T}^2)$ if and only if $\mu$
is not a Liouville number \citep{greenfield1972global} and $\int_{\mathbb{T}^2} f(x,y) \dint x \dint y = 0$.

Thus, if $\mu$ is a non-Liouville number, the dimension of $H^{0,1}(\mathbb{T}^2; \frk h)$ is one. It is also possible to show that, when $\mu$ is rational, the dimension of $H^{0,1}(\mathbb{T}^2; \frk h)$ is infinity and $\dext'$ has closed range. We prove this claim on Lemma \ref{lem:torus_liouville}. When $\mu$ is a Liouville number, we have that $\dext' (C^\infty(\mathbb{T}^2))$ is not closed and, hence, by Proposition \ref{prp:finite_dim_and_closedness}, the dimension
$H^{0,1}(\mathbb{T}^2; \frk h)$ also is infinity. We prove this fact on Lemma \ref{lem:torus_liouville_2}.
\end{exa}

\begin{lem}
\label{lem:torus_liouville} If $\frk h$ is defined as in Example \ref{exa:torus_liouvile} and $\mu$ is a rational number, then the dimension of $H^{0,1}(\mathbb{T}^2; \frk h)$ is infinity and the range of the operator $$\mathcal  C^\infty(\mathbb{T}^2) \xrightarrow[]{ \dext'} \{ f \dext x : f \in \mathcal C^\infty(\mathbb{T}^2) \}$$
is closed.
\end{lem}

\begin{proof}
To show that the dimension of $H^{0,1}(\mathbb{T}^2; \frk h)$ is infinity, we only have to show that there exists an arbitrarily large linearly independent set contained in $H^{0,1}(\mathbb{T}^2; \frk h)$. Let $k > 0$ be an integer. Since $\mu$ is rational, the equation $\xi - \mu \eta = 0$ has infinite zeros.

Thus, let $(\xi_j, \eta_j) \in \Z \times \Z$, $j=0, \ldots, k$ be a sequence of distinct zeros of the equation $\xi - \mu \eta = 0$. We are going to show that the set $\{e^{i (\xi_j x + \eta_j y)} \dext x \}_{j=0}^k$ is linearly independent. Let $c_j \in \C$ be constants. To show that this
set is linearly independent, it is enough to show that $$\sum_{j=0}^k c_j e^{i (\xi_j x + \eta_j y)}\dext x$$ is $\dext'$-exact only when every constant $c_j$ is zero. By using Fourier series, we can easily see that a solution to the equation $$Lf = \sum_{j=0}^k c_j e^{i (\xi_j x + \eta y_j)}$$ exists if, and only if, every constants $c_j$ is zero.

To prove that the range of the operator $\dext'$ is closed, it is enough to see that, if $f_j$ is in the range of $\dext'$ and is a sequence converging to $f$ in $\mathcal C^\infty(\mathbb{T}^2)$,
it holds that $\hat f_j(\xi,\eta) = 0$ whenever $\xi - \mu \eta = 0$. Also, since $f_j \to f$ in $\mathcal C^\infty(\mathbb{T}^2)$, we have that
$\hat f_j(\xi, \eta) \to \hat f(\xi, \eta) $ for all $(\xi, \eta)$ and, thus, $\hat f(\xi, \eta) = 0 $ whenever $\xi - \mu \eta = 0$. This implies that $f$ is in the range of $\dext'$.
In fact, if we take
$$u = \sum_{(\xi,\eta) \in \Z^2,~ \xi - \mu \eta \neq 0 } \frac{ \hat f(\xi, \eta)}{ (\xi - \mu \eta) } e^{i (\xi_j x + \eta_j y)},$$ we have that $u$ is smooth, since $\mu$ is rational, and it holds that $Lu = f$.
\end{proof}

\begin{lem}
  \label{lem:torus_liouville_2}
  If $\frk h$ is defined as in Example \ref{exa:torus_liouvile} and $\mu$ is a Liouville number, then the range of the operator $$\dext' : \mathcal  C^\infty(\mathbb{T}^2) \xrightarrow[]{ } \{ f \dext x : f \in \mathcal C^\infty(\mathbb{T}^2) \}$$ is not closed.
\end{lem}

\begin{proof}
  We again adapt some ideas from \cite{greenfield1972global}. If $\mu$ is a Liouville number, then there exists a sequence $\{(\xi_j, \eta_j)\}_{j \in \N}$ such that $\{\xi_j\}$ and $\{\eta_j\}$ are increasing and
  $|\xi_j - \mu \eta_j| \leq 1/(\xi_j^2 + \eta_j^2)^j$. Let $L$ be as in Example \ref{exa:torus_liouvile}. By defining $f_k = \sum_{j \leq k} \exp(i\xi_jx + i\eta_jy)$, we have that $f_k, Lf_k \in \mathcal C^\infty(\mathbb{T}^2)$ and
  $f = \lim_k f_k \in \mathcal D'(\mathbb{T}^2) \bs \mathcal C^\infty(\mathbb{T}^2)$, but with $Lf \in \mathcal C^\infty(\mathbb{T}^2).$
  Notice that any other solution of $u$ of $Lu = Lf$ is of the form $u = f + c$, with $c \in \C$. Therefore, $Lf$ is not in the range of the operator $\dext' : \mathcal  C^\infty(\mathbb{T}^2) \xrightarrow[]{ } \{ f \dext x : f \in \mathcal C^\infty(\mathbb{T}^2) \}$.
\end{proof}


\subsection{Construction of left-invariant involutive structures}
\label{subsec:existence}

In Definition \ref{dfn:basic_types_of_lie_algebras}, we introduced some types of Lie algebras by defining certain special involutive structures. Let $\Omega$ and $\VV$ be the quadric introduced in Example \ref{exa:3}.
Notice that $\Omega$ has a natural structure of Lie group given by
 $$ (z,w) \circ (z_0, w_0) \doteq (z + z_0 + 2i\mc Q(w, \overline{w}_0), w + w_0), \qquad (z,w), (z_0, w_0) \in \Omega.$$ The involutive structure $\VV$
 is left-invariant (we refer to \cite{boggess1991cr} to a detailed description of this class of CR manifolds and the computation of its Levi form). In this case, the Levi form is left-invariant.

Now we are going to see how to construct some more examples of such Lie algebras.

Let $G$ be a compact connected Lie group and consider any maximal torus $\mathbb{T} \subset G$. Let $\frk t$ be the complexification of the Lie algebra of $\mathbb{T}$. The Lie algebra $\frk t$ is abelian and self-normalizing, that is, if $[X,Y] \in \frk t$ for all $X \in \frk t$, then $Y \in  \frk t$. We endow $\frk g_\R$ with an $ad$-invariant inner product, that is, an inner product satisfying
$$ \ip{ [X, Y], Z} = -\ip{ Y, [X, Z]},$$
for all $X, Y, Z \in \frk g_\R$.
Also, we extend it to an Hermitian inner product on $\frk g$, which is going to satisfy
$$ \ip{ [X, Y], Z} = -\ip{ Y, [\overline X, Z]}$$
for all $X, Y, Z \in \frk g.$

Let $T_1, \ldots, T_r$ be a basis for $\frk t_\R$ and consider the operators $\operatorname{ad}_{T_j} : \frk g \to \frk g$, given by $\operatorname{ad}_{T_j}(X) = [T_j, X]$. These operators have the following properties:

\begin{enumerate}
	\item $\operatorname{ad}_{T_j}$ and $\operatorname{ad}_{T_k}$ commute for every $j$ and $k$ (Jacobi's identity);
    \item $\operatorname{ad}_{T_j}$ is skew-Hermitian (and diagonalizable) for all $j$;
    \item all $\operatorname{ad}_{T_j}$ share the same eigenspaces.
\end{enumerate}

If $\frk g'$ is one of the eigenspaces, we denote by $\alpha = (\alpha_1, \ldots, \alpha_r)$ the ordered set of eigenvalues associated to $\operatorname{ad}_{T_1}, \ldots, \operatorname{ad}_{T_r}$, specifically, $\operatorname{ad}_{T_j}(X) = \alpha_j X$ for all $X \in \frk g'$. We denote this eigenspace by $\frk g_\alpha$. Each $\alpha$ is called a root of the Lie algebra. We denote by $\Delta$ the set of all roots. Notice that every eigenvalue is purely imaginary.

By the spectral theorem, we have a decomposition
$$ \frk g = \frk t \oplus \bigoplus_{\alpha \in \Delta} \frk g_\alpha.$$

It is sometimes convenient to deal with $\alpha$ as an element of $\frk t^*$,
in which $\alpha(T_j)$ is defined by $\operatorname{ad}_{T_j}(X) = \alpha(T_j)X$ for $X \in \frk g_\alpha.$

Suppose that $X \in \frk g_\alpha$ and $Y \in \frk g_\beta$. By Jacobi's identity we have

$$
\begin{aligned}
	\operatorname{ad}_{T_i}([X,Y]) & = -[Y,[T_i,X]] - [X,[Y,T_i]] \\
    	& = -[Y, \alpha(T_i)X] -[X,-\beta(T_i)Y] \\
        & = \left\{\alpha(T_i) + \beta(T_i) \right\}[X,Y].
\end{aligned}
$$

Thus, we conclude that $[X,Y] \in \frk g_{\alpha + \beta}$ if $\alpha + \beta$ is a root and zero otherwise.

Let $\Delta_+ \subset \Delta$ be a maximal set of roots satisfying the following properties:
\begin{enumerate}
	\item for all $\alpha \in \Delta$, exactly one of $\alpha$ or $-\alpha$ is in $\Delta_+$;
	\item if $\alpha, \beta \in \Delta_+$ and $\alpha + \beta$ is a root, then $\alpha + \beta \in \Delta_+$.
\end{enumerate}

\begin{rmk} The dimension $\frk t$, usually called the \emph{rank} of $G$, is independent of the choice of the maximal torus and $\dim G - \dim_\C \frk t$ always is an even number.
\end{rmk}

Now we can define some left-invariant involutive structures.

Since $\frk t$ is abelian, it is trivial to construct essentially real, elliptic, complex, or CR subalgebras of $\frk t$.
In fact, take $\frk t_\R$ the real Lie algebra of $\mathbb{T}$ and let $X_1, \ldots, X_r$ be a basis for $\frk t_\R$. We have that $\frk t = \frk t_\R \otimes \C$. We identify $X_j$ with $X_j \otimes 1$. Let $s,t \geq 0$ and $s + t \leq r$. Consider the subalgebra
$$ \frk u = \operatorname{span}_\C \{ X_s \ldots, X_{s+t}, X_{s+t+1} + i X_{s+it+2}, \ldots X_{r-1} + iX_{r} \}.$$
Notice that these are involutive structures on the maximal torus.

If $s > 0$ and $t = 0$, we have an essentially real subalgebra; if $s, t > 0$ and $s + t = r$, we have an elliptic subalgebra of $\frk u$; if $r$ is even and $s = 0$, we have a complex subalgebra of $\frk u$; and, finally, if $s = 0$ and $t > 0$, we have a CR structure of $\frk u$.

Now consider $\frk h = \frk u \oplus \bigoplus_{\alpha \in \Delta_+} \frk g_\alpha $. By the preceding discussion,
$\frk h$ is a Lie subalgebra and is elliptic if $\frk u$ is elliptic, is complex if $\frk u$ is complex, and is CR if $\frk u$ is CR.

\begin{rmk} We can remove the hypothesis of $G$ being compact by taking a Cartan subalgebra $\frk t$ and by proving the statement without using the $\operatorname{ad}$-invariant metric.
\end{rmk}

Let $\frk h \subset \bigoplus_{\alpha \in \Delta} \frk g_\alpha$ be any elliptic subalgebra of $\bigoplus_{\alpha \in \Delta} \frk g_\alpha$. If $\frk e \subset \frk t$ is an elliptic subalgebra of $\frk t$, then $\frk e \oplus \frk h$ is an elliptic subalgebra of $\frk g$. The most obvious example of such structure is when $\frk h = \bigoplus_{\alpha \in \Delta_+} \frk g_\alpha$. 

Now we are going to construct some concrete examples.

\subsection{Involutive structures on \texorpdfstring{$\operatorname{SU}(2)$}{\operatorname{SU}(2)} and \texorpdfstring{$\operatorname{SU}(3)$}{\operatorname{SU}(3)}}

The group $\operatorname{SU}(2)$ is defined as
$$ \operatorname{SU}(2) \doteq \left\{ \left(  \begin{matrix} z_1 & -\overline{z_2} \\ z_2 & \overline{z_1} \end{matrix} \right) : z_1, z_2 \in \C, ~|z_1|^2 + |z_2|^2 = 1 \right\}.$$ The rank of $\operatorname{SU}(2)$ is 1.

The Lie algebra of $\operatorname{SU}(2)$, denoted by $\frk{su}(2)$ is generated by
$$ X = \left(  \begin{matrix} 0 & i \\ i & 0 \end{matrix} \right), \quad Y = \left(  \begin{matrix} 0 & -1 \\ 1 & 0 \end{matrix} \right), \quad  T = \left(  \begin{matrix} i & 0 \\ 0 & -i \end{matrix} \right).$$

We have the following relation between $X, Y$ and $T$:

\begin{equation}
\label{sec:su2:basic_commutator_relations}
[T, X] = 2Y, \quad [T, Y] = -2X, \quad [X, Y] = 2T.
\end{equation}

\begin{exa}
 Let $L = X - iY$ and consider the structure $\frk h = \operatorname{span}_\C \{L\}$. This is clearly a CR structure. We have that $\{L, \overline{L},T\}$ is a basis for $\frk{su}(2) \otimes \C$ and we denote its dual basis by $\{\zeta, \overline{\zeta}, \tau\}$. We note that  $\frk h^\perp_0$ is a vector space of dimension 1 and that $\tau \in \frk h^\perp_0$. Thus, if
 $\xi \in \frk h^\perp_0$, then we can write $\xi = \lambda \tau$, in which $\lambda \in \R \bs \{0\}$. We note that $\frk h$ also has dimension 1 and so $Z = \alpha L$ for all $Z \in \frk h$. Then, the Levi form for the structure $\frk h$ is given by
$$ \mathscr L_\xi (Z, Z) = \frac{1}{2i} \lambda \alpha \overline{\alpha} \tau([L, \overline{L}]) = \frac{1}{2i} \lambda \alpha \overline{\alpha} 4i \tau(T) = 2\lambda|\alpha|^2. $$
In conclusion, the Levi form is non-degenerate and does not change sign for all $Z \in \frk h$. 
\end{exa}

It is not possible to have a hypocomplex structure of rank 1 on $\operatorname{SU}(2)$. Therefore, $\operatorname{SU}(2)$ has no hypocomplex CR structure. We are going to show that $\operatorname{SU}(2)$ does not admit a Levi degenerate left-invariant structure.

Suppose that $\operatorname{SU}(2)$ admits a Levi degenerate left-invariant CR structure $\frk h$. Since it is left-invariant, we have that the essentially real structure $\frk h + \overline{\frk h}$ is a Lie algebra. Take $\frk k = \frk{su}(2) \cap (\frk h + \overline{\frk h})$ and notice that $\frk k \otimes \C = \frk h + \overline{\frk h}$.

We are going to see that $\frk k$ is abelian. In fact, let $X, Y \in \frk k$ be linearly independent. We have that $\ip{[X,Y],X} = \ip{Y,[X,X]} = 0$ and that $\ip{[X,Y],Y} = -\ip{X,[Y,Y]} = 0$. Therefore, $[X,Y]$ is orthogonal to $\frk k$. Since $\frk k$ is a Lie algebra, we also have that $[X,Y] \in \frk k$ and so the only possibility is that $[X,Y] = 0$. Therefore, there is an abelian subalgebra of rank 2. This is a contradiction because the rank of
$\operatorname{SU}(2)$ is 1.

This CR structure is related to the natural complex structure of $\C^2$. There is a natural diffeomorphism between the group $\operatorname{SU}(2)$ and the 3-sphere $\mathbb{S}^3$, which has a natural CR structure inherited from $\C^2$. This structure, when pulled back to $\operatorname{SU}(2)$, is exactly the left-invariant structure we just described. Now we are going to prove this claim:


The complex structure of $\C^2$ is given by the involutive bundle $\mathcal V = \bigcup_{p\in\C^2} \mathcal V_p$, in which $$\mathcal V_p = \operatorname{span} \left\{ \left. \frac{\del}{\del \overline{z}_j} \right|_p : j = 1, 2 \right\}$$ and
$$ \frac{ \del }{ \del \overline{z}_j } = \frac{1}{2} \left( \frac{ \del }{ \del x_j } + i \frac{ \del }{ \del y_j }\right), \qquad j = 1,2.$$

By writing $\rho(z_1, z_2) = |z_1|^2 + |z_2|^2  - 1$, we have that
$ \mathbb{S}^3 = \rho^{-1}(0)$. Clearly, $\dext \rho_p \neq 0$ if $p \in \bb S^3 $, so we have
$$ \C T_p \bb S^3 = \{ Z \in \C T_p \C^2 : Z_p (\rho) = 0 \}.$$

Notice that $\mathcal W_p = \mathcal V_p \cap \C T_p \bb S^3$ defines a CR structure on $\mathbb S^3$. We want a description of $\mathcal W_p$ to simplify some computations. We note that
$$ \mathcal W_p = \{ Z \in \mathcal V_p : Z_p (\rho) = 0 \}.$$

Notice that if $Z \in \mathcal V_p$, by taking $p_1, p_2 \in \C$ with $p = (p_1, p_2)$, then $$Z = \alpha_1 \left. \frac{\del}{\del \overline{z}_1} \right|_p + \alpha_2 \left. \frac{\del}{\del \overline{z}_2} \right|_p$$ and so $(Z \rho)(p) = 2(a_1p_1 + a_2p_2) = 0.$
Now we have
$$ \mathcal W_p = \left\{ \alpha_1 \left. \frac{\del}{\del \overline{z}_1} \right|_p + \alpha_2 \left. \frac{\del}{\del \overline{z}_2} \right|_p
 : \alpha_1, \alpha_2 \in \C;~ a_1p_1 + a_2p_2 = 0 \right\}.$$

The identification of $\mathbb S ^3$ with $\operatorname{SU}(2)$ is given by the diffeomorphism
$$  (z_1, z_2) \in \mathbb S ^3 \xmapsto{\theta} \left(  \begin{matrix} z_1 & -\overline{z_2} \\ z_2 & \overline{z_1} \end{matrix} \right) \in \operatorname{SU}(2).$$

This map endows $\mathbb S ^3$ with a structure of Lie group with multiplication given by
$$ (\alpha_1,\alpha_2) \cdot (z_1,z_2) \doteq \theta^{-1}(\theta(\alpha_1,\alpha_2) \cdot \theta(z_1,z_2)) = (\alpha_1 z_1 - \overline{\alpha_2} z_2, \alpha_2 z_1 + \overline{\alpha_1} z_2).$$

Now, for $\alpha \in \mathbb S^3$, we have the diffeomorphism $ L_\alpha : \mathbb S^3 \to \mathbb S ^3 $ defined by $L_\alpha(z_1, z_2) = (\alpha_1, \alpha_2) \cdot (z_1, z_2)$. Notice that this map has a natural extension to $\C^2$ and each component is a holomorphic function. More than that, for $\alpha \neq 0$, this map is actually a biholomorphism from $\C^2$ to $\C^2$ and so it preserves the complex structure of $\C^2$. Therefore, it clearly restricts to a CR diffeomorphism from
$\mathbb S^3$ to $\mathbb S ^3$. What we proved is that the involutive structure $\mathcal W$ is invariant by the action of $L_\alpha$, or, in other words, $\mathcal W$ is left-invariant.

Now we have almost everything we need to relate the CR structure we introduced on $\mathbb{S}^3$ with the abstract CR structure we constructed on $\operatorname{SU}(2)$. To simplify even more the exposition, let $p \in \C^2$ and identify $\C T_p \C^2$ with $\C^2$ by
\begin{equation}
\label{eq:adele}
\left. \frac{ \del }{ \del x_1 } \right|_p \mapsto (1,0), \left. \frac{ \del }{ \del y_1 } \right|_p \mapsto (i,0), \left. \frac{ \del }{ \del x_2 } \right|_p \mapsto (0,1) \text{ and } \left. \frac{ \del }{ \del y_2 } \right|_p \mapsto (0,i).
\end{equation}

With this identification, we have that $\mathcal W_p$ is identified with $$\{(z_1, z_2) \in \C^2 : z_1p_1 + z_2p_2 = 0 \}.$$ Notice that $ \theta $ is linear and so its differential can be identified with $\theta$ itself. Also, notice that, on $p = (1,0)$, we have that $Z \in \mathcal W_{(1,0)}$ is of the form $Z = \Lambda \del / \del \overline{z_2}$, in which $\Lambda \in C$. Therefore, all we need to compute is $\theta_*(\del / \del \overline{z_2})$.
By using the identification \eqref{eq:adele},
we have
$$
\begin{aligned}
\theta_*(2\del / \del \overline{z_2}) & = \theta_*((0,1) + i(0,i)) \\
		& = \left(  \begin{matrix} 0 & -1 \\ 1 & 0 \end{matrix} \right) + i\left(  \begin{matrix} 0 & i \\ i & 0 \end{matrix} \right) \\
        & = i(X - iY).
\end{aligned}$$

\begin{exa}
\label{exa:there_are_examples}
 Let $L = X - iY$ and consider the structure $\frk h = \operatorname{span}_\C \{L, T\}$. This is clearly a left-invariant elliptic structure.
\end{exa}


We recall that $\operatorname{SU}(3)$ is the group of all unitary $3 \times 3$ matrices with complex coefficients having determinant 1. It is a real Lie group, compact, and has dimension 8 and rank 2. The Lie algebra of $\operatorname{SU}(3)$ is denoted by $\frk{su}(3)$ and is the set of all traceless skew-Hermitian $3 \times 3$ matrices with complex coefficients. The following matrices form a basis for $\frk{su}(3)$:
$$
\begin{aligned}
T_1 &= \left( \begin{matrix} i & 0 & 0 \\ 0 & -i & 0 \\ 0 & 0 & 0 \end{matrix} \right), & \quad T_2 &= \left( \begin{matrix} i & 0 & 0 \\ 0 & i & 0 \\ 0 & 0 & -2i\end{matrix} \right), \\
X_1 &= \left( \begin{matrix} 0 & i & 0 \\ i & 0 & 0 \\ 0 & 0 & 0 \end{matrix} \right), & \quad Y_1 &= \left( \begin{matrix} 0 & -1 & 0 \\ 1 & 0 & 0 \\ 0 & 0 & 0\end{matrix} \right),  \\
X_2 &= \left( \begin{matrix} 0 & 0 & i \\ 0 & 0 & 0 \\ i & 0 & 0 \end{matrix} \right), & \quad Y_2 &= \left( \begin{matrix} 0 & 0 & -1 \\ 0 & 0 & 0 \\ 1 & 0 & 0\end{matrix} \right),  \\
X_3 &= \left( \begin{matrix} 0 & 0 & 0 \\ 0 & 0 & i \\ 0 & i & 0 \end{matrix} \right), & \quad Y_3 &= \left( \begin{matrix} 0 & 0 & 0 \\ 0 & 0 & -1 \\ 0 & 1 & 0\end{matrix} \right).  \\
\end{aligned}
$$

We have the following relations between the coefficients. Each cell corresponds to the commutator between the first element of the line and the first element of the column. Since the commutator is anti-symmetric, we omitted half of the commutators.


\bgroup
\def\arraystretch{1.5}%
\begin{center}\begin{tabular}{c || c | c | c | c | c | c | c | c}
         & $T_1$ & $T_2$ & $X_1$  & $Y_1$   & $X_2$       &  $Y_2$       & $X_3$     & $Y_3$   \\
\hline
\hline
$T_1$    & $0$   & $0$   & $2Y_1$ & $-2X_1$ & $Y_2$       &  $-X_2$      & $-Y_3$   & $X_3$       \\
\hline
$T_2$    & -     & $0$   & $0$    & $0$     & $3Y_2$      &  $-3X_2$     & $3Y_3$   & $-3X_3$     \\
\hline
$X_1$    & -     &  -    & 0      & $2T_1$  & $Y_3$       &  $-X_3$      & $Y_2$    & $-X_2$      \\
\hline
$Y_1$    & -     & -     & -      & 0       & $X_3$       &  $Y_3$       & $-X_2$   & $-Y_2$      \\
\hline
$X_2$    & -     & -     & -      & -       & 0           &  $T_2 + T_1$ & $Y_1$    & $X_1$       \\
\hline
$Y_2$    & -     & -     & -      & -       & -           & 0            & $-X_1$   & $Y_1$       \\
\hline
$X_3$    & -     & -     & -      & -       & -           & -            & 0        & $T_2 - T_1 $\\
\hline
$Y_3$    & -     & -     & -      & -       & -           & -            & -        & 0           \\

\end{tabular}\end{center}
\egroup

\begin{exa}
Considering the commutators we just computed, we have an obvious CR structure. We define the following vector fields on $\frk{su}(3) \otimes \C$:
$$L_1 = X_1 - iY_1, \quad L_2 = X_2 - iY_2, \quad L_3 = X_3 - iY_3 .$$

Also, from the commutators, it follows that, for each $j,k$, there exists a real number $\Lambda_{jk}$ such that $[T_j,X_k] = \Lambda_{jk} Y_k$ and $[T_j,Y_k] = -\Lambda_{jk} X_k$. By combining these last two equations,
we have $ [T_j, L_k] = -i\Lambda_{jk} L_k $, which means that each $L_k$ is an eigenvector for the map $\operatorname{ad}_{T_j}$.

Notice that $[L_1,L_2] = 0$, $[L_1, L_3] = 2iL_2$, and $[L_2,L_3] = 0$. Thus, $$\frk h = \operatorname{span}_\C \{L_1, L_2, L_3\}$$ is a Lie subalgebra of $\frk{su}(3) \otimes \C$, which is clearly a CR subalgebra.
\end{exa}

We want to compute the Levi form associated to this Lie algebra. Each cell on the following table corresponds to the commutator between the first element of the line and the first element of the column.

\bgroup
\def\arraystretch{1.5}%
\begin{center}\begin{tabular}{c||c | c | c }
                 & $\overline{L_1}$                 & $\overline{L_2}$           & $\overline{L_3} $      \\
\hline
\hline
$L_1$ & $ 4iT_1 $            & $-2iL_3$         &$0$            \\
\hline
$L_2$ & $-2i \overline{L_3}$ & $2i(T_2 + T_1)$   & $-2i \overline{L_1}$ \\
\hline
$L_3$ & 0           &  $-2i L_1$       & $2i(T_2 - T_1)$              \\
\end{tabular}\end{center}
\egroup


A generic element of $\frk h$ is of the form $ Z = \alpha_1 L_1 + \alpha_2 L_2 + \alpha_3 L_3 $, in which $\alpha_1, \alpha_2, \alpha_3 \in \C$. Therefore, we have that
$$[Z, \overline{Z}] = |\alpha_1|^2 4i T_1 + |\alpha_2|^2 2i(T_2 + T_1) + |\alpha_3|^2 2i(T_2 - T_1) + W,$$ in which $W$ is a linear combination of $X_1, Y_1, X_2, Y_2, X_3, Y_3$.

Also, we have that $$\{ L_1, L_2, L_3, \overline{L_1}, \overline{L_2}, \overline{L_3}, T_1, T_2 \}$$ is a basis for $\frk{su}(3) \otimes \C$ and we denote by $$\{\zeta_1,\zeta_2,\zeta_3,\overline{\zeta_1},\overline{\zeta_2}, \overline{\zeta_3}, \tau_1, \tau_2 \}$$
its dual basis. For any $\xi \in \frk h^\perp_0$, we can write $ \xi = \xi_1 \tau_1 + \xi_2 \tau_2 $, in which $\xi_1, \xi_2 \in \R$ and $(\xi_1, \xi_2) \neq 0$.

Notice that the Levi form of $\frk h$ is then
$$
\begin{aligned}
\mathscr L_{\xi} (Z,Z) & = \frac{1}{2i}\xi([Z, \overline{Z}]) \\
	& = 2|\alpha_1|^2 \xi_1 + |\alpha_2|^2(\xi_2 + \xi_1) + |\alpha_3|^2(\xi_2 - \xi_1) \\
    & = \xi_1(2|\alpha_2|^2 + |\alpha_2|^2 - |\alpha_3|^2) + \xi_2(|\alpha_2|^2 + |\alpha_3|^2).
\end{aligned}
$$

If $\xi_1 = 0$ and $\xi_2>0$, then clearly there is no $Z$ such that $\mathscr L_{\xi} (Z,Z) = \xi_2(|\alpha_2|^2 + |\alpha_3|^2) < 0$. Therefore, by \cite[Theorem 6.2]{baouendi1983microlocal}, we conclude that this structure is not hypocomplex.

Now we are going to use this structure to construct a hypocomplex one.







\begin{exa}
 Let $\frk h' = \operatorname{span}_\C\{L_1, L_2, L_2, aT_1 + bT_2\}$. Notice that $L_1, L_2, L_3$, and their conjugates are eigenvectors for $\operatorname{ad}_{T_j}$ and that $[aT_1 + bT_2, aT_1 + bT_2] = 0$. Therefore, this is a Lie subalgebra of $\frk{su}(3) \otimes \C$.
\end{exa}

We claim that this structure is hypocomplex. A generic element of $\frk h'$ is of the form $ Z = \alpha L_1 + \beta L_2 + \gamma L_3 + \delta U$, in which $ U = aT_1 + bT_2$ and $\alpha, \beta, \gamma, \delta$ are complex numbers. Therefore, we have that $$[Z, \overline{Z}] = |\alpha|^2 4i T_1 + |\beta|^2 2i(T_2 + T_1) + |\gamma|^2 2i(T_2 - T_1) + W,$$ in which $W$ is a linear combination of $X_j, Y_j$, for $j = 1, 2, 3$.

For any $\theta \in \frk h^\perp_0$, we can write $ \theta = t(-b \tau_1 + a \tau_2)$, in which $t \in \R$ and $t \neq 0$. To simplify the proof, we are going to take $t=1$.

Hence, the Levi form of $\frk h'$ at $(e,\theta)$ can be computed as:
$$
\begin{aligned}
\mathscr L_{\theta} (Z,Z) & = \frac{1}{2i}\theta([Z, \overline{Z}]) \\
	& =  -2|\alpha|^2 b + |\beta|^2(a - b) + |\gamma|^2(a+b)  \\
    & = b(|\gamma|^2 - 2|\alpha|^2 - |\beta|^2) + a(|\beta|^2 + |\gamma|^2).
\end{aligned}
$$

From this expression, if $b \neq 0$, we can easily see that, for any $\theta \in \frk h^\perp_0$, it is possible to find a $Z$ such that $\mathscr L_{\theta} (Z,Z) < 0$.
Therefore, by \cite[Corollary 6.1]{baouendi1983microlocal}, we conclude that this structure is hypocomplex.

Now we are going to use the structure just defined to construct an elliptic one.

\begin{exa}
 Consider the structure $\frk h'' = \operatorname{span}_\C\{L_1, L_2, L_3, T_1, T_2\}$. Notice that $L_1, L_2, L_3$, and their conjugates are eigenvectors for $\operatorname{ad}_{T_j}$ for $j = 1, 2$ and that $[T_1, T_2] = 0$. This structure is clearly elliptic.
\end{exa}



\section{Left-invariant cohomologies}
\label{sec:ce_tech}

One of the objectives of this work
is to show that, under some reasonable conditions, we can compute the cohomology of the complex \eqref{h_complex} by restricting our attention only to left-invariant forms.

To accomplish this, we recall some basic definitions and expose some techniques that Claude Chevalley and Samuel Eilenberg introduced in \cite{chevalley1948cohomology} to study the cohomology of the de Rham complex.

Then we extend some of these techniques to construct left-invariant cohomologies for the operator $\dext'$. Finally, we discuss conditions for usual cohomology spaces to be isomorphic to the left-invariant cohomology spaces.

\subsection{Left-invariant de Rham cohomology}

Let $G$ be a compact Lie group, $\frk X(G)$ be the set of all smooth vector fields on $G$, and $\SF{G}{p}{}$ be the set of all smooth $p$-forms on $G$. We recall that a vector field $X \in \frk X(G)$ is called left-invariant if, for every $x \in G$, we have
$ (L_x)_* X = X$.
Also, a differential form $u \in \SF{G}{p}{})$ is called left-invariant if, for every $x \in G$, we have $ (L_x)^* u = u $. We denote by $\frk X_L(G)$ the set of all left-invariant vector fields and by $\SF{G}{p}{L}$ the set of all left-invariant $p$-forms.

Since $G$ is a compact Lie group, it can be endowed with a Haar measure, which we denote by $\mu$.

\begin{lem}
 If $u \in \SF{G}{p}{L}$, then $\dext u \in \SF{G}{p+1}{L}$.
\end{lem}

\begin{proof} It follows directly from the fact that the exterior derivative commutes with pullbacks.
\end{proof}

Considering the notation from Chapter \ref{chp:involutive_structures},
denote the $p$-cocycles by $Z^p(G)$, the $p$-coboundaries by $ B^p(G)$, and the $p$-cohomology classes by $H^p(G) = Z^p(G) / B^p(G).$

Also, define in the obvious way the left-invariant $p$-cocycles, left-invariant $p$-coboundary, and left-invariant $p$-cohomology classes, denoted, respectively, by $Z^p_L(G)$, $B^p_L(G)$, and $H^p_L(G)$.

Notice that, for each integer $p$, the inclusion $\SF{G}{p}{L} \hookrightarrow \SF{G}{p}{}$ induces a linear mapping $$i_p : H^p_L(G) \to H^p(G).$$

To show that this linear mapping is an isomorphism, we need some basic tools.

We start by defining an operator which allows us to use a useful averaging trick:


\begin{lem}
\label{averaging_operator_lemma}
 Let $u \in \SF{G}{p}{}$. The differential form
\begin{equation}
\label{averaging_operator}
  Au = \int_G L^*_x u \dint \mu(x)
\end{equation}
has the following properties:
\begin{enumerate}
\item \label{deRhamLeftInvariantCommutative} $\dext (Au) = A(\dext u)$;
\item $Au$ is left-invariant;
\item if $u$ is left-invariant, then $u = Au$.
\end{enumerate}
\end{lem}

This lemma is just a particular case of a more general result stated and proved at \cite[page 90]{chevalley1948cohomology}, so the proof is omitted.

Notice that Lemma \ref{averaging_operator_lemma} shows that the following diagram is commutative:
$$
\xymatrix
{
  \SF{G}{p}{} \ar[r]^{\dext~} \ar[d]_{A} & \SF{G}{p+1}{} \ar[d]^{A} \\
  \SF{G}{p}{L} \ar[r]_{\dext~} & \SF{G}{p+1}{L}
}
$$

The following lemma shows that $i_p$ is injective.

\begin{lem}
\label{averaging_trick_left_invariant_exactness}
Let $u \in \SF{G}{p}{L}$ be such that $u = \dext v$ for some $v \in \SF{G}{p-1}{}$. Then, there exists $w \in \SF{G}{p}{L}$ such that $u = \dext w$.
\end{lem}

\begin{proof}
We have $u = Au = A\dext v = \dext A v \in \dext \SF{G}{p}{L}$, so we can take $w = Av$.
\end{proof}

To prove that $i_p$ is surjective, we need the following theorem, whose proof can be founded in \cite{chevalley1948cohomology}:

\begin{thm}
\label{thm:exact}
 If $u$ is a closed form such that $\int_c u = 0$ for every homology class, then it is exact.
\end{thm}

Since Theorem \ref{thm:exact} cannot be applied to the general case, we omit its proof. 


\begin{thm}
\label{coho_left_equiv} If $u \in \SF{G}{p}{}$ is closed, then there exists $v \in \SF{G}{p-1}{}$ such that $u + \dext v = Au$.
\end{thm}

\begin{proof}
We want to show that the form $u - Au$ is exact. Since it is already closed, we only need to show that $\int_c (u - Au) = 0$ for every $p$-dimensional
homology class $c$.
$$ \int_c Au = \int_c \int_G L_x^* u \dext
\mu(x) = \int_G \int_c L_x^* u \dext \mu(x) = \int_G \int_{L_xc}  u \dext \mu(x) =  \int_G \int_c  u \dext \mu(x) = \int_c  u.$$
Hence, $\int_c u - Au = 0$ for every $c$.
\end{proof}

Therefore, for each $p$, we have an isomorphism of vector spaces: $$ H^p(G) \cong H^p_L(G).$$

These isomorphisms simplify a lot of computations because they show that we can restrict our attention to left-invariant forms. Also, left-invariant forms are
completely determined by their values at the tangent space at the identity. In other words, we reduced the problem to a problem in linear algebra.


\subsection{Left-invariant cohomology relative to a subalgebra}
\label{sec:left_invariant_bidegree_cohomology}

We define $N^{0,q}(G; \frk h) = \mathcal C^\infty(G; \Lambda^q)$ and, for $p > 0$, we define
$$N^{p,q}(G; \frk h) = \{ u \in \mathcal C^\infty(G; \Lambda^{p+q}): u(X_1, \ldots, X_{p+q}) = 0 \text{ when } q+1 \text{ arguments are in } \frk h \}.$$

Also, we define $N^{p,q}_{L}(G; \frk h) \doteq N^{p,q}(G; \frk h) \cap \SF{G}{p+q}{L}$ and, since $$\dext \SF{G}{p}{L} \subset \SF{G}{p+q+1}{L},$$ we can restrict the exterior derivative and obtain
$$ \dext_L : N^{p,q}_{L}(G; \frk h) \to N^{p,q+1}_{L}(G; \frk h).$$
Therefore, $\dext_L$ induces a coboundary operator on the quotient
$$ \SF{G}{p,q}{L} = N^{p,q}_{L}(G; \frk h) / N^{p+1,q-1}_{L}(G; \frk h),$$
which is denoted by $\dext^{p,q}_{\frk h, L}.$ When there is no risk of confusion, we omit the bidegree writing only $\dext_{\frk h, L}$. When the structure $\frk h$ is obvious by the context, we write $\dext'_{L}$.

For $p \geq 0$, we denote the set of the left-invariant $(p,q)$-cocycles elements by $$ Z_{L}^{p,q} (G; \frk h) = \ker \left( \dext'_L : \SF{G}{p,q}{L}) \to \SF{G}{p,q+1}{ L}) \right),$$
the set of left-invariant $(p,q)$-coboundaries by
 $$ B^{p,q}_{L}(G; \frk h) = \img \left( \dext'_L : \SF{G}{p,q-1}{L}) \to \SF{G}{p,q}{L} \right),$$
and the left-invariant $(p,q)$-cohomology classes by
 $$ H^{p,q}_{L}(G; \frk h) = \frac{Z^{p,q}_{L
 }(G; \frk h)}{B^{p,q}_{L}(G; \frk h)}.$$

Our next objective is to construct an averaging operator $$A^{p,q}_{\frk h} : \SF{G}{p,q}{}) \to \SF{G}{p,q}{L})$$ by making commutative the diagram
$$
\xymatrix
{
  \SF{G}{p,q}{} \ar[r]^{ \dext^{p,q}_{\frk h} } \ar[d]_{ A^{p,q}_{\frk h}} & \SF{G}{p,q+1}{} \ar[d]^{ A^{p,q+1}_{\frk h} } \\
  \SF{G}{p,q}{L} \ar[r]_{\dext^{p,q}_{\frk h,L} } & \SF{G}{p,q+1}{ L} .
}
$$
Then, we are going to find the necessary and sufficient conditions so that we have an isomorphism of vector spaces
$$ H^{p,q}(G;\frk h) \cong H^{p,q}_{L}(G; \frk h).$$



Denote by $ \pi^{p,q}_{\frk h} : N^{p,q}(G; \frk h) \to \SF{G}{p,q}{}$ and $ \pi^{p,q}_{\frk h, L} : N^{p,q}_{L}(G; \frk h) \to \SF{G}{p,q}{L}$ the quotient maps.

\begin{lem}
\label{lem:op_A_and_pi_commutes}
The operator $A$ defined on Lemma \ref{averaging_operator_lemma} satisfies:
$$ A(N^{p,q}(G; \frk h)) \subset N^{p,q}(G; \frk h). $$
Thus, there is an unique operator $A^{p,q}_{\frk h}$ making the following diagram commutative:
$$
\xymatrix
{
  N^{p,q}(G; \frk h) \ar[r]^{A} \ar[d]_{\pi^{p,q}} & N^{p,q}_{L}(G; \frk h) \ar[d]^{\pi^{p,q}_L} \\
  \SF{G}{p,q}{} \ar[r]_{A^{p,q}_{\frk h} } & \SF{G}{p,q}{L}
}
$$
\end{lem}

\begin{proof}
In fact, let $X_1, \ldots, X_{p+q} \in \frk X(G)$ and suppose that $q+1$ of these vector fields are in $\frk h$. We can assume that $X_1, \ldots, X_{q+1} \in \frk h$. It follows that
$$
\begin{aligned}
 Au (X_1, \ldots, X_{p+q}) & = \int_G (L^*_x u) (X_1, \ldots, X_{p+q}) \dint \mu(x) \\
	& = \int_G u ((L_x)_*X_1, \ldots, (L_x)_*X_{p+q}) \mu(x) \\
	& = \int_G u (X_1, \ldots, X_{q+1}, (L_x)_*X_{q+2},\ldots, (L_x)_*X_{p+q}) \mu(x)
	&= 0.
\end{aligned}
$$
The operator $$A^{p,q}_{\frk h} : \SF{G}{p,q}{} \to \SF{G}{p,q}{\frk h,L}$$ is defined by $$A^{p,q}_{\frk h} (\pi^{p,q}(u)) =  \pi_{p,q+1,L}(Au).$$ The uniqueness is obvious.
\end{proof}

Notice that the operator $ A^{p,q}_{\frk h} $ satisfies
$$ \dext^{p,q}_{\frk h,L} \circ A^{p,q}_{\frk h} = A^{p,q}_{\frk h} \circ \dext^{p,q}_{\frk h}.$$
This follows directly from the definition of $A^{p,q}_{\frk h}$ and from the fact that $A$ and $\dext$ commute.


\subsection{Necessary and sufficient conditions for the isomorphism}
\label{subsection_necessary_sufficient}

Denote by $\phi^{p,q}_{\frk h}$ the quotient map $$\phi^{p,q}: Z^{p,q}(G; \frk h) \to H^{p,q}(G; \frk h)$$ and by $\phi^{p,q}_{\frk h, L}$ the quotient map $$\phi^{p,q}_{\frk h, L}: Z^{p,q}_{L}(G; \frk h) \to H^{p,q}_{L}(G; \frk h).$$

Notice that, if $u \in \SF{G}{p,q}{L}$, we can write $ u = \phi^{p,q}_{\frk h, L}(u') $, in which $ u' \in N^{p,q}_{L}(G; \frk h) \subset N^{p,q}(G; \frk h) $. Also, we can define $i(u) = \phi^{p,q}_{\frk h}(u') \in \SF{G}{p,q}{}$ so we have a natural map
$$ i : \SF{G}{p,q}{L} \to \SF{G}{p,q}{}.$$
Since this map is injective, we can identify $ \SF{G}{p,q}{L} $ with its image
under $i$ and assume that $\SF{G}{p,q}{L} \subset \SF{G}{p,q}{} $. The map $ i : \SF{G}{p,q}{L} \to \SF{G}{p,q}{}$ satisfies with $\dext_{\frk h} \circ i = i \circ \dext_{\frk h, L}$, so it induces a map
$  i : H^{p,q}_{L}(G; \frk h) \to H^{p,q}(G; \frk h).$

\begin{lem}
\label{inclusion_left_invatiant_cohomology_injective} The map $i : H^{p,q}_{L}(G; \frk h) \to H^{p,q}(G; \frk h)$ is injective.
\end{lem}

\begin{proof}
 Let $u \in H^{p,q}_{L}(G; \frk h)$, with $i(u) = 0$. We write $u = \pi^{p,q}_{\frk h, L}(u')$, in which $u' \in \SF{G}{p,q}{L}$ satisfies $\dext'_L u' = 0$. Since $i(u) = 0$, there exists $v \in \SF{G}{p,q}{}$ with $\dext' v = u$. Thus, we have that $\dext'(Av) = A \dext' v = A u = u$
 and so $u = 0$.
\end{proof}

Since the map $i : H^{p,q}_{L}(G; \frk h) \to H^{p,q}(G; \frk h)$ is always injective, we can already show some algebraic obstructions that can prevent the cohomology space $H^{p,q}(G; \frk h)$ from being zero.
That is, we always have $ \dim H^{p,q}_{L}(G; \frk h) \leq \dim H^{p,q}(G; \frk h)$.

To show that the map $$i : H^{p,q}_{L}(G; \frk h) \to H^{p,q}(G; \frk h)$$ is surjective, given $v \in \SF{G}{p,q}{}$, we have to find $u \in \SF{G}{p,q}{L}$ such that $v - u = \dext' w$ for some $w \in \SF{G}{p,q-1}{}$.
The natural candidate for $u$ is $Av$. Thus, we want to find a $w \in \SF{G}{p,q-1}{}$ satisfying $Av - v = \dext' w$.
In other words, we want to solve the following problem: given $v \in \ker A^{p,q}_{\frk h}$ with $\dext' v = 0$, can we find $u \in \SF{G}{p,q}{}$ satisfying $\dext' u = v$? If $\frk h$ is a involutive structure such that this problem has a solution in degree $(p,q)$, then we say that it has \textit{property (K) in degree $(p,q)$}.

\begin{exa}
 In the case in which $\frk h = \frk g$, we have that $\dext'$ is the usual exterior derivative. Also, Theorem \ref{coho_left_equiv} shows that we always have property (K).
\end{exa}

\begin{lem}
 Suppose that the map $i : H^{p,q}_{\frk h, L}(G) \to H^{p,q}_{\frk h}(G)$ is surjective. Then the involutive structure $\frk h$ has property (K) in degree $(p,q)$.
\end{lem}

\begin{proof}
In fact, let $v \in \SF{G}{p,q}{}$ such that $\dext' v = 0$ and satisfying $Av = 0$. By hypothesis, there exists $u \in \SF{G}{p,q}{\frk h,L}$ satisfying $\dext' u = 0$ and $v - u = \dext' \alpha$. By applying $A$ on both sides, we obtain $Av - Au = \dext'(A\alpha)$. Since $Av = 0$ and $Au = u$, we have $u = -\dext'(A\alpha)$ and so
$v = \dext'(\alpha - A\alpha)$.
\end{proof}

\begin{exa}
 In the complex case, when the group $G$ is semisimple, it is possible to use a result by Bott \cite{bott1957homogeneous} to see that we always have property (K). The application of the Bott's result was made explicit by Pittie on \cite{pittie1988dolbeault}.
\end{exa}

\begin{exa}
 Consider the Example \ref{exa:torus_liouvile} and notice that, if $(\xi_0, \eta_0) \neq (0,0)$, then $u = e^{i (\xi x + \tau t)  } \dext x \in \SF{\T^2}{0,1}{L}$ satisfies $Au = 0$. If  $(\xi_0, \eta_0)$ is a zero of $\xi - \mu\tau = 0$,
 we \emph{cannot} solve the problem $\dext' f = u$ for $f \in \SF{\T^2}{0,1}{L}$. This can be easily seem by using Fourier series.
\end{exa}

\begin{lem}
 If the involutive structure $\frk h$ satisfies property (K) in degree $(p,q)$, then the mapping  $i : H^{p,q}_{L}(G; \frk h) \to H^{p,q}_{}(G; \frk h)$ is surjective.
\end{lem}

Therefore, we have the following theorem:

\begin{thm}
 If the involutive structure $\frk h$ satisfies property (K) in degree $(p,q)$, the cohomologies $H^{p,q}(G; \frk h)$ and $H^{p,q}_{ L}(G; \frk h)$ are isomorphic.
\end{thm}


\section{The Hodge decomposition for the operator \texorpdfstring{$\dext'$}{d'}}

Let $W_1, \ldots, W_N$ be a basis for $\frk g$ with dual basis $w_1, \ldots, w_N$. We endow $\frk g$ with the Hermitian inner product making the basis $W_1, \ldots, W_N$ an orthonormal basis, that is
$$W_j \cdot W_k \doteq \delta_{jk}, \qquad j,k=1,2, \ldots N.$$

Notice that $W_j \cdot W_k = w_j(W_k)$ for $j,k = 1,2, \ldots N.$

We denote by $\Lambda^p \frk g$ the $p$th exterior product of $\frk g$ and we recall that, by using the universal property, we can naturally identify $\Lambda^p \frk g^*$ with $(\Lambda^p \frk g)^*$, the dual of $\Lambda^p \frk g$.

We write $W_I = W_{i_1} \wedge \ldots \wedge W_{i_p} \in \Lambda^p \frk g$ and denote its dual by $ w_I =  w_{i_1} \wedge \ldots \wedge w_{i_p}$. Now we have that $\{w_I: |I| = p\}$ is a basis for $\Lambda^p \frk g^*$ and we can define an Hermitian inner product on $\Lambda^p \frk g^*$ by $$w_I \cdot w_J = w_I(W_J), \quad \text{ for all } I,J \text{ with } |I|=|J| = p.$$

If $u \in \SF{G}{p}{}$, we can write $u = \sum_{|I| = p} u_I w_I$ and, for each $x \in G$, we have that $$u_x = \sum_{|I| = p} u_I(x) w_I \in \Lambda^p \frk g^*.$$
Therefore, we can endow $\SF{G}{p}{}$ with the inner product given by
$$(u,v) \doteq \int_G u_x \cdot v_x \dint \mu(x),$$ 
for $u, v \in \SF{G}{p}{}$.

Without loss of generality, we assume that the first $n$ elements of $W_1, \ldots, W_N$ is a basis for $\frk h$. Also, we denote these elements by $L_1, \ldots, L_n$ and the other elements by $M_1, \ldots, M_m$. That is, $L_j = W_j$, for $j=1,\dots,n$, and $M_j = W_{n+j}$, for $j=1, \ldots, m$.

We denote the dual basis for $\{L_1, \ldots, L_n,M_1, \ldots, M_m\}$ by $\{\tau_1, \ldots, \tau_n,\zeta_1, \ldots, \zeta_m\}$.

With this basis, each element of $\SF{G}{p,q}{}$ can be identified with an element in $\SF{G}{p+q}{}$. Thus, each $\SF{G}{p,q}{}$ is endowed with an Hermitian inner product induced by the inner product in $\SF{G}{p,q}{}$.
%
%
%
%
We finally endow $\SF{G}{p,q}{}$ with the inner product
$$(u,v) \doteq \int_G u(x) \cdot v(x) \dint \mu(x),$$ 
for $u, v \in \SF{G}{p,q}{}$.

With this inner product, for each degree $(p,q)$, we define the formal adjoint of the operator $$\dext^{p,q}_{\frk h} : \SF{G}{p,q}{} \to \SF{G}{p,q+1}{}$$ as the unique operator $$\delta^{p,q}_{\frk h} : \SF{G}{p,q+1}{} \to \SF{G}{p,q}{}$$ such that
$$ (\dext^{p,q}_{\frk h}u, v) = (u, \delta^{p,q}_{\frk h}v),$$
for all $u \in \SF{G}{p,q}{}$ and $v \in \SF{G}{p,q+1}{}$.



For each pair $p,q$ with $p,q \geq 0$, we define the Hodge Laplacian, also known as Box operator, $$\Box^{p,q} : \SF{G}{p,q}{} \to \SF{G}{p,q}{}$$
by $$\Box^{p,q} = \delta^{p,q}_{\frk h} \circ \dext^{p,q}_{\frk h} + \dext^{p,q-1}_{\frk h} \circ \delta^{p,q-1}_{\frk h}.$$


With conditions that are going to be introduced later, we have the decomposition:
$$ \SF{G}{p,q}{} = \ker \Box^{p,q} \oplus \dext^{p,q-1}(\SF{G}{p,q-1}{}) \oplus \delta^{p,q}(\SF{G}{p,q+1}{}).$$

With this decomposition, we have: suppose that $u \in \SF{G}{p,q}{}$ is such that $\dext' u = 0$. Clearly, we have that $(u,v) = 0$ for all $v \in \delta^{p,q}(\SF{G}{p,q+1}{})$. Thus, if we show that $(u,v) = 0$ for all $v \in \ker \Box^{p,q}$, we have that $u$ is $\dext'$-exact.


\subsection{Applications of the Hodge decomposition}

Suppose that, in degree $(p,q)$, the kernel of the Box operator has only the constant sections. In other words, we are supposing that $\ker \Box^{p,q} \subset \ker (A - I)$.

We are going to show that, with this assumption, the involutive structure has property (K) in degree $(p,q)$. In other words, we are going to show that, given $v \in \ker A^{p,q}$ with $\dext' v = 0$, we can find $u \in \SF{G}{p,q}{}$ satisfying $\dext' u = v$.

We already have $\iip{v,u} = 0$ for all $v \in \delta^{p,q}_{\frk h}(\SF{G}{p,q+1}{})$. Thus, we just have to show that $\iip{v,u} = 0$ for all $u \in \ker \Box^{p,q}$. This is easy. We have that $u = Au$ for all $u \in \ker \Box^{p,q}$ and thus
$\iip{v,u} = \iip{v,Au} = \iip{Av,u} = 0.$

Here we used that the operator $A$ is formally self-adjoint. We prove this fact in the following lemma:

\begin{lem}
\label{lem:formally_self-adjoint}
 The operator $A : \SF{G}{p}{} \to \SF{G}{p}{}$ defined in Lemma \ref{averaging_operator_lemma} is formally self-adjoint.
\end{lem}

\begin{proof}
Let $u, v \in \SF{G}{p}{}$ and write
$$ u = \sum_{|I|=p} u_I w_I, \quad v = \sum_{|J|=p} v_J w_J.$$
Notice that, for $X_1, \ldots, X_p \in \frk g$, we have
$$
\begin{aligned}
	(Av)_g(X_1, \ldots, X_p)
		& = \int_G ((L_x)^* v)_g(X_1, \ldots, X_p) \dint \mu(x) \\
		& = \int_G (L_x)^* \left(\sum_{|J|=p} v_J w_J\right)(X_1, \ldots, X_p) \dint \mu(x) \\
		& = \sum_{|J|=p} \int_G ((L_x)^* v_J)(g) \wedge (L_x)^*(w_J) \dint \mu(x) \\
		& = \sum_{|J|=p} \left(\int_G v_J \circ L_x(g) \dint \mu(x) \right) w_J  \\
		& = \sum_{|J|=p} \left(\int_G v_J \dint \mu(x) \right) w_J.  \\
\end{aligned}
$$

Thus, it holds that
$$\begin{aligned}
	(u,Av)	& = \sum_{|I|=p}\sum_{|J|=p} \int_G u_I(g) \left(\int_G \overline{ v_J \circ L_x(g)} \dint \mu(x) \right) w_I \cdot w_J \dint \mu(g) \\
		& = \sum_{|I|=p}\sum_{|J|=p} \left( \int_G \int_G u_I(g) \overline{ v_J \circ L_x(g)} \dint \mu(x)  \dint \mu(g) \right) w_I \cdot w_J \\
		& = \sum_{|I|=p}\sum_{|J|=p} \left( \int_G \int_G u_I \circ L_x(g) \overline{ v_J (x)} \dint \mu(x)  \dint \mu(g) \right) w_I \cdot w_J \\
		& = (Au,v).
\end{aligned}$$
\end{proof}


\subsection{Lie derivatives}

Let $X$ be a real smooth vector field on $G$. We denote by $\Phi^X$ the flow of $X$ in $G$, that is, $\Phi^X(x, t)$ is a smooth curve on $G$ satisfying $\Phi^X(x, 0) = x$ and $ \dext / \dext t |_{t=0} \Phi^X(x, t) = X(x)$. We also use the notation $\Phi^X_t(x) = \Phi^X(x, t)$. If $f \in \mathcal C^\infty(G)$ is a smooth function, we define the Lie derivative of $f$ in $x \in G$ by
$$ (\mathscr{L}_X f)(x) = \left. \frac{ \dext }{ \dext t} \right|_{t=0} f(\Phi^X(t,x)).$$
Also, if $u \in \SF{G}{q}{}$ is a smooth $p$-form, we define its Lie derivative in $x$ by
$$ (\mathscr{L}_X u)_x = \left. \frac{ \dext }{ \dext t} \right|_{t=0} ((\Phi^X_t)^* u)_x.$$

\begin{prp}
The Lie derivative has the following properties: \begin{enumerate}
	\item If $u$ and $v$ are two smooth forms, it holds that
$$ \mathscr{L}_X (u \wedge v) = (\mathscr{L}_X u) \wedge v + u \wedge (\mathscr{L}_X v).$$
	\item The Lie derivative commutes with the exterior derivative, that is, $$ \mathscr{L}_X \circ \dext = \dext \circ \mathscr{L}_X.$$
	\item Cartan's Magic Formula:
$$ \mathscr{L}_X u = \imath_X (\dext u) + \dext (\imath_X u).$$
\end{enumerate}
\end{prp}

Another important property we are going to use is the following:

\begin{lem}
\label{lem:adele}
 Let $u \in \SF{G}{p}{}$ be a real right-invariant form. Then $\mathscr{L}_X u = 0$ for any real left-invariant vector field $X$.
\end{lem}

\begin{proof}
 In fact, the flow of $X$ on $G$ can be written as $\Phi^X_t(g) = g\exp(tX) = R_{\exp(tX)}(g)$. Thus,
$$ (\mathscr{L}_X u)_x = \lim_{t \to 0} \frac{1}{t} \left[ (R_{\exp(tX)})^* u)_x - u_x \right] = 0$$
because $(R_{\exp(tX)})^* u)_x - u_x = 0$ due to the right-invariance of $u$.
\end{proof}

If $Z \in \frk g$, we write $Z = X + i Y$, in which $X, Y \in \frk g_\R$, and we define
$$ \mathscr{L}_Z u = \mathscr{L}_X u + i \mathscr{L}_Y u.$$

If $Z \in \frk h$, we can define the Lie derivative on $\SF{G}{p,q}{}$. In fact, if $Z \in \frk h$ and $u \in N^{p,q}(G; \frk h)$, then notice that $\imath_Z u \in N^{p,q-1}(G; \frk h)$. Thus, we have an operator $$ \imath_X : N^{p,q}(G; \frk h) \to N^{p,q-1}(G; \frk h)$$ satisfying $\imath_Z (N^{p+1,q-1}) \subset N^{p+1,q-2}.$
Also, it defines a map
$$ \imath_X : \SF{G}{p,q}{} \to \SF{G}{p,q-1}{}$$
making the following diagram commutative:
$$
\xymatrix{
N^{p,q}(G; \frk h) \ar[d]_{\pi^{p,q}} \ar[r]^{\imath_Z} & N^{p,q-1}(G; \frk h) \ar[d]^{\pi^{p,q-1}} \\
\SF{G}{p,q}{} \ar[r]^{\imath_Z} & \SF{G}{p,q-1}{}.
}
$$
Therefore, by using Cartan's Magic Formula
$$ \mathscr{L}_Z u = \imath_Z (\dext u) + \dext (\imath_Z u),$$
we obtain
$$
\begin{aligned}
  \pi^{p,q} ( \mathscr{L}_Z u ) = & \pi^{p,q}(\imath_Z (\dext u)) + \pi^{p,q}(\dext (\imath_Z u)) \\
	= & \imath_Z (\pi^{p,q-1}(\dext u)) + \dext' ( \pi^{p,q-1}(\imath_Z u)) \\
	= & \imath_Z (\dext' \pi^{p,q-1}(u)) + \dext' ( \imath_Z \pi^{p,q-1}(u)).
\end{aligned}
$$

The Cartan's Magic Formula preserves elements in $\SF{G}{p,q}{}$ if $Z$ is an element of the involutive structure. Thus, we can define an operator $$ \mathscr{L}'_Z : \SF{G}{p,q}{} \to \SF{G}{p,q}{}$$
by
$$ \mathscr{L}'_Z u = \imath_Z (\dext' u) + \dext' (\imath_Z u). $$

\begin{lem}
 Let $u \in \SF{G}{p,q}{}$ be a right-invariant form. Then $\mathscr{L}_Z u = 0$ for any left-invariant vector field.
\end{lem}

\begin{proof}
 The proof is just a straightforward application of Lemma \ref{lem:adele}.
\end{proof}


\subsection{Relation between \texorpdfstring{$\Box$}{Box} and \texorpdfstring{$\Delta$}{Laplacian}}

With the inner product $( \cdot , \cdot )$, each $\SF{G}{p}{}$ is a pre-Hilbert space and we can define $$\delta : \SF{G}{p+1}{} \to \SF{G}{p}{},$$ the formal adjoint of the operator $$\dext : \SF{G}{p}{} \to \SF{G}{p+1}{}.$$

The Laplace-Beltrami operator is defined by
$$ \Delta \doteq \dext \circ \delta + \delta \circ \dext.$$

\begin{lem}
\label{lem:lbo}
 It holds that $\delta ( N^{p,q+1}(G; \frk h) ) \subset N^{p,q}(G; \frk h).$
\end{lem}

\begin{proof}
 Notice that the restriction of $\dext$ give us an operator $$\tilde {\dext} : N^{p,q}(G; \frk h) \to N^{p,q+1}(G; \frk h).$$ Thus, we can construct the formal adjoint of this operator, denoted by $\tilde{\delta}$. By construction, we have that $\tilde{\delta} (N^{p,q+1}(G; \frk h)) \subset N^{p,q}(G; \frk h)$. So all we have to do now is show that $\tilde{\delta}$ and $\delta$ agree on $ N^{p,q+1}(G; \frk h)$, that is, $\tilde{\delta} = \delta|_{ N^{p,q+1}(G; \frk h) } $.
 Let $u \in N^{p,q+1}(G; \frk h)$. We have that, for all $v \in N^{p,q}(G; \frk h)$,
$$ (\tilde{\delta}u - \delta u, v ) = (\tilde{\delta}u,v)  - (\delta u, v ) = (u, \tilde{\dext} v ) - ( u, \dext v) = 0.$$
Thus, $\tilde{\delta}u = \delta u.$
\end{proof}

By using Lemma \ref{lem:lbo}, we have $\dext \circ \delta + \delta \circ \dext : N^{p,q}(G; \frk h) \to N^{p,q}(G; \frk h),$ that is, the Laplace-Beltrami operator satisfies $$\Delta (N^{p,q}(G; \frk h)) \subset N^{p,q}(G; \frk h).$$ Thus, it can be defined on the quotient $\SF{G}{p,q}{}$. That is, there exists an unique operator $\Box^{p,q}_{\frk h}$
making the following diagram commutative:
$$
\xymatrix{
N^{p,q}(G; \frk h) \ar[r]^{ \Delta } \ar[d]_{ \pi^{p,q} } & N^{p,q}(G; \frk h) \ar[d]^{ \pi^{p,q} } \\
\SF{G}{p,q}{} \ar[r]_{ \Box^{p,q}_{\frk h} } & \SF{G}{p,q}{}.
}
$$

This operator on the quotient is just the operator $\Box$ we introduced earlier.

Notice that, if $[u] \in \SF{G}{p,q}{}$, then $\Box [u] = [\Delta u]$ and, if
$$ \mathscr{L}'_Z \Box [u] = \mathscr{L}'_Z [\Delta u] = [\mathscr{L}_Z \Delta u] = [\Delta \mathscr{L}_Z u],$$ we have the following result:

\begin{prp}
 The operator $\mathscr{L}'_Z$ commutes with $\Box^{p,q}$.
\end{prp}

\begin{cor}
If $\Box u = 0$, then $ \Box \mathscr{L}'_Z u = 0$, $\dext' \mathscr{L}'_Z u = 0$, and $\delta' \mathscr{L}'_Z u = 0.$
\end{cor}

\begin{prp}
 If $u \in \ker \Box$, then $\mathscr{L}'_Z u = 0$ for all $Z \in \frk h$.
\end{prp}
\begin{proof}
 Since $u \in \ker \Box$, we have that $\dext' u = 0$. Thus, $\mathscr{L}'_Z u = \dext' (i_Z u) + i_Z (\dext' u) = \dext' (i_Z u)$ and
$$(\mathscr{L}'_Z u, \mathscr{L}'_Z u) = (\mathscr{L}'_Z u, \dext' (i_Z u)) = (\delta' \mathscr{L}'_Z u, i_Z u) = 0.$$
\end{proof}



\begin{prp}
\label{chp:comp_lie:hypo_p_zero} Let $\frk h$ be a hypocomplex structure. If $u \in \SF{G}{p,0}{}$ is such that  $\mathscr{L}'_Z u = 0$ for all $Z \in \frk h$, then $u$ is left-invariant.
\end{prp}

\begin{proof}
Let $u \in \SF{G}{p,0}{}$ and assume that $\mathscr{L}'_Z u = 0$ for all $Z \in \frk h$. We write $ u = \sum_{|I|=p} u_I \zeta_I $, in which $u_I \in \mathcal C^\infty(G)$. Notice that
 $$\mathscr{L}'_Z \left( \sum_{|I|=p} u_I \zeta_I \right) =  \sum_{|I|=p} \left\{ \mathscr{L}'_Z (u_I) \zeta_I + u_I (\mathscr{L}'_Z \zeta_I) \right\}.$$

By definition, we have $\mathscr{L}'_Z \zeta_I = \dext' (\imath_Z \zeta_I) + \imath_Z(\dext' \zeta_I)$. Since $\imath_Z \zeta_I$ is zero, we have $\dext' (\imath_Z \zeta_I) = 0$. Also, there exist constants $a_{Ij} \in \C$ such that $$\dext' \zeta_I = \sum_{j=1}^n c_{Ij} \zeta_I \wedge \tau_j.$$ Thus, by applying $\imath_Z$ on both sides, we obtain
$$\imath_Z(\dext' \zeta_I) = \imath_Z(\sum_{j=1}^n c_{Ij} \zeta_I \wedge \tau_j) = \sum_{j=1}^n c_{Ij} \imath_Z(\zeta_I \wedge \tau_j) = 0$$
because $\imath_{Z}(\zeta_I \wedge \tau_j) = 0$ when computed in any element of $\frk h$.

Therefore, $\mathscr{L}'_Z (\sum_{|I|=p} u_I \zeta_I) =  \sum_{|I|=p} \mathscr{L}'_Z (u_I) \zeta_I = 0$ and so $\mathscr{L}'_Z u_I = 0$ for all $Z$. Since $\frk h$ is hypocomplex, we conclude that $u_I$ is a constant. Therefore, the form $u$ is left-invariant.
\end{proof}

Notice that the Lie derivative can be extended to work with currents and the same proof yields the following result:

\begin{cor}
\label{cor:hypo_p_zero_in_Eprime} Let $\frk h$ be a hypocomplex structure. If $u \in \mathcal D'(G; \Lambda^{p,0})$ is such that  $\mathscr{L}'_Z u = 0$ for all $Z \in \frk h$, then $u$ is left-invariant.
\end{cor}

\begin{thm}
\label{thm:comp_lie:hypo_subelliptic}
 Let $\frk h$ be an hypocomplex structure on a compact Lie group $G$. Then every cohomology class $u \in H^{p,0}_{\mathcal C^\infty}(G; \frk h)$ has a representative that is left-invariant. That is,
 $$ H^{p,0}_{\mathcal C^\infty}(G; \frk h) = H^{p,0}_{L}(G; \frk h).$$
\end{thm}

\begin{prp}
\label{prp:comp_lie:hypo_subelliptic_torus}
 Let $\frk h$ be an hypocomplex structure on a torus $\mathbb{T}$. If $u \in \SF{\T}{p,q}{}$ is such that $\mathscr{L}'_Z u = 0$ for all $Z \in \frk h$, then $u$ is left-invariant.
\end{prp}

\begin{proof}
 If $u \in \Lambda^{p,q}(\T)$ is such that $\mathcal L'_Z u = 0$ for all $Z \in \frk h$, by writing $$ u = \sum_{|I|=p} \sum_{|J|=q} u_{IJ} \zeta_I \wedge \tau_J,$$
 we have that $$ \mathscr L'_Z u = \sum_{|I|=p} \sum_{|J|=q} \mathscr L'_Z ( u_{IJ} ) \zeta_I \wedge \tau_J = 0. $$
 Thus, $\mathscr L'_Z (u_{IJ}) = 0$. Since $\frk h$ is hypocomplex, we have that $u_{IJ} \in \C$ and the theorem is proved.
\end{proof}

\begin{thm}
\label{thm:comp_lie:hypo_subelliptic_torus_2}
 Let $\frk h$ be a hypocomplex structure on a torus $\mathbb{T}$ such that the associated box operator is subelliptic. Then, every cohomology class $u \in H^{p,q}(\T; \frk h)$ has a representative that is left-invariant.
\end{thm}

\begin{proof}
 Let $[u] \in \Lambda^{p,q}(\T)$. Since we are assuming $\Box$ to be subelliptic, we can assume $\Box u = 0$ and $\mc L'_Z u = 0$ for all $Z \in \frk h$, which implies that $u$ is left-invariant. The proof is concluded.
\end{proof}


\section{Application of Serre duality}

For this application, we need to extend the averaging operator acting on forms to an operator acting on currents, which can be done by transposition. Since the averaging operator is formally self adjoint, as proved at Lemma \ref{lem:formally_self-adjoint}, we can easily see that this extension is well-defined.

The following lemma is going to be useful:

%

\begin{lem}
 Let $u \in \mathcal D'(G; \Lambda^k)$ and write $\Sigma_{|I| = k} u_I w_I$. Then $Au = \Sigma_{|I| = k} (Au_I) w_I.$
\end{lem}

\begin{proof}
 Let $\phi \in \mathcal C^\infty(G; \Lambda^{N-k})$ and write $\phi = \Sigma_{|I'| = N-k} \phi_{I'} w_{I'}$. By definition, we have
$$
\begin{aligned}
	(Au)(\phi) & = u(A \phi) \\
    		& = u(\Sigma_{|I'| = N-k} (A\phi_{I'}) w_{I'}) \\
            & = (\Sigma_{|I| = k} u_I w_I)(\Sigma_{|I'| = N-k} (A\phi_{I'}) w_{I'}) \\
            & = \Sigma_{|I| = k} \Sigma_{|I'| = N-k} u_I(A\phi_{I'}) w_I \wedge w_{I'} \\
            & = \Sigma_{|I| = k} \Sigma_{|I'| = N-k} (Au_I)(\phi_{I'}) w_I \wedge w_{I'} \\
            & = (\Sigma_{|I| = k} (Au_I) w_I)(\phi).
\end{aligned}
$$
Thus, the proof is completed.
\end{proof}

\begin{lem}
 Let $u \in \mathcal D'(G; \Lambda^k)$ and suppose that $Au - u = 0$. Then $u$ is a smooth left-invariant form.
\end{lem}

\begin{proof}
 Let $u \in \mathcal D'(G; \Lambda^k)$ and write $u = \Sigma_{|I| = k} u_I w_I$. We note that, if $0 = Au - u$, then $Au_I - u_I = 0$, which means that each $u_I$ is a constant.
\end{proof}

\begin{thm}
\label{thm:li:thm:duality_of_left} Let $\frk h$ be any involutive structure such that the operators
  \begin{equation}
    \label{eq:two_operators_closed_range}
    \mathcal C^\infty(G; \Lambda^{p,q-1}) \xrightarrow{\dext'} \mathcal C^\infty(G; \Lambda^{p,q}) \xrightarrow{\dext'} \mathcal C^\infty(G; \Lambda^{p,q+1})
\end{equation}
 have closed range. Then, it holds that every cohomology class in degree $(p,q)$ has a left-invariant representative if, and only if, every cohomology class in degree
$(m-p, n-q)$ also has a left-invariant representative.
\end{thm}

\begin{proof}
Since the operators \eqref{eq:two_operators_closed_range} have closed range, by Theorem \ref{chp:inv:thm:serre_duality}, we have $$ H^{p,q}_{\mathcal C^\infty}(G; \frk h)^* \cong H^{m-p,n-q}_{\mathcal D'}(G; \frk h).$$
Let $[u] \in H^{m-p,n-q}_{\mathcal D'}(G; \frk h)$. For every $[v] \in H^{p,q}_{\mathcal C^\infty}(G; \frk h)$, we have
$$ (Au - u)(v) = Au(v) - u(v) = u(Av) - u(v) = u(Av - v).$$

Suppose that every class in $H^{p,q}_{\mathcal C^\infty}(G; \frk h)$ has a representative that is left-invariant. Thus, we can assume that $v$ is left-invariant, that is, $Av - v = 0$. Then, $[Au - u] = 0$ in $H^{m-p,n-q}_{\mathcal D'}(G; \frk h)$.
The other direction of the equivalence follows by applying the exact same argument.
\end{proof}

\begin{rmk}
\label{chp:li:rmk:const} Notice that, since $[Av] = [v]$ in $H^{m-p,n-q}_{\mathcal D'}(G; \frk h)$ and $Av$ is left-invariant, we actually have $$H^{m-p,n-q}_{\mathcal D'}(G; \frk h) \cong  H^{m-p,n-q}_{\mathcal C^\infty}(G; \frk h) \cong H^{m-p,n-q}_{L}(G; \frk h).$$
\end{rmk}

\begin{thm}
\label{thm:left_invariance_on_degree} Let $\frk h$ be a left-invariant hypocomplex structure on a compact Lie group $G$. Suppose that $$\mathcal C^\infty(G; \Lambda^{p,n-1}) \xrightarrow{\dext'} \mathcal C^\infty(G; \Lambda^{p,n})$$ has closed range. Then we have $$H^{m-p,n}_{\mathcal C^\infty}(G; \frk h) = H^{m-p,n}_{L}(G; \frk h).$$
\end{thm}

\begin{proof}
If $[u] \in H^{p,0}_{\mathcal{E}'}(G; \frk h)$, then clearly $\Box u = 0$, which means that $\mathscr L'_Z u = 0$ for all $Z \in \frk h$. Therefore, $u$ is left-invariant by Corollary \ref{cor:hypo_p_zero_in_Eprime}. By combining this with Theorem \ref{thm:li:thm:duality_of_left} and Remark \ref{chp:li:rmk:const}, we have left-invariance of cohomologies in top degree, that is,
$$H^{p,n}_{\mathcal C^\infty}(G;\frk h) = H^{m-p,n}_{L}(G; \frk h).$$
\end{proof}

\begin{rmk}
\label{rmk:adele_labels_everything}
 The operator $\dext' : \mathcal C^\infty(G; \Lambda^{p,q}) \xrightarrow{} \mathcal C^\infty(G; \Lambda^{p,q+1})$ is known to have closed range in many situations. For example, it always has closed range in the case of elliptic structures. For hypocomplex structures, the situation is more delicate, but we covered some sufficient conditions in Chapter \ref{chp:involutive_structures}.
 For example, see Theorem \ref{thm:closed_range} and Proposition \ref{prp:finite_dim_and_closedness}. Notice that, due to the left-invariance of the structures, we only need to verify the hypothesis in the origin.
\end{rmk}

\chapter{Involutive structures on homogeneous manifolds}
\label{chp:homogeneous}
In this chapter, we introduce the concept of invariant involutive structures on homogeneous manifolds and show how to construct some examples. The main reason for introducing this concept is to simplify the study left-invariant involutive structures on Lie groups.

The main idea we want to explore is the following: let $G$ be a compact Lie group endowed with an elliptic involutive structure $\frk h$. Also, let $\frk k_\R = \frk h \cap \frk g_\R$ and consider the group $K = \exp(\frk k_\R)$. We assume that $K$ is closed, so we can define a homogeneous manifold $\Omega = G / K$. This homogeneous manifolds inherits a complex involutive structure from $G$ via the quotient map.

With all these ingredients in hand, we want to see how much we can infer about the cohomology spaces $H^{p,q}(G;\frk h)$ from the De Rham cohomology of $K$ and from the Dolbeault cohomology of $\Omega$. This is done via a spectral sequence. In order to use this spectral sequence, we need some further assumptions. One of the assumptions is that the real orbits of the elliptic structure $\frk h$ is closed and the other one is that $\frk h$ can be decomposed into a direct sum of two Lie subalgebras: one is the algebra $\frk k = \frk h \cap \frk g$ and the other is an ideal in $\frk h$. We show that this tecnique can always be applied in a situation similar to a famous conjecture by Treves.

To conclude, we introduce the concept of Lie algebra cohomology and use it to state a theorem by Bott, which, under certain assumptions, gives us algebraic information about the Dolbeault cohomology of $\Omega$. Under certain algebraic and topological restrictions on $G$, Bott's theorem, along with the spectral sequences, give us a complete algebraic description of $H^{p,q}(G;\frk h)$.


\section{Homogeneous manifolds and invariant involutive structures}

Let $\Omega$ be a homogeneous manifold for a Lie group $G$. We recall that $\Omega$ is a smooth manifold and that there exists a smooth map $$T : G \times \Omega \to \Omega,$$ called \emph{left-action} of $G$ on $\Omega$.
Other notations for $T$ are $T(g, x) = T_g x = g \cdot x$.
The map $T$ satisfies the following properties:
$$
	g \cdot (g' \cdot x) = (gg') \cdot x \text{ and } e \cdot x = x,
$$ in which $e \in G$ is the identity of $G$ and $T$ is assumed to be \emph{transitive}, that is, for any two points $x,y \in \Omega$, there is an element $g$ such that $g \cdot x = y$.

\begin{exa} Any Lie group acts transitively on itself by left multiplication.
\end{exa}

\begin{exa}
If $G$ is any Lie group and $H$ is a closed Lie subgroup, the space $\Omega = G / H$ is a homogeneous manifold with $G$ acting on $\Omega$ by
$$ g \cdot (g' H) = (g g') \cdot H.$$

By Theorem 7.19 on \cite{lee2003smooth}, we know that all examples of homogeneous manifolds are equivalent to this example and, thus, we have that $\Omega$ is an analytic manifold and the map $T$ is an analytic map.
\end{exa}

From now on, $\Omega$ is always assumed to be a connected $G$-homogeneous manifold.

We want to use involutive structures on homogeneous manifolds. In order to take into consideration the symmetries given by the group action, we need to consider involutive structures that carry information about the action. We then introduce the following concept:

\begin{dfn}
 An involutive structure $\mathcal V \subset \C T\Omega$ is said to be \emph{invariant by the action of $G$} if $(T_g)_*X \in \mathcal V_{g \cdot x}$ for all $X \in  \mathcal V_x.$
\end{dfn}

The most obvious example is the following:

\begin{exa}
 If $\mathcal{V} = \C T \Omega$, then clearly $\mathcal V$ is invariant.
\end{exa}


Next, we construct some examples of invariant involutive structures.


Let $G$ be a Lie group and denote its Lie algebra by $\frk g_{\R}$. We denote by $\frk g$ the complexification of $\frk g_{\R}$.


Let $X$ be a vector field on a smooth manifold $\Omega$.
Suppose that $\mathcal V \subset \C T\Omega$ is a vector bundle. We say that $\mathcal V$ is preserved by the vector field $X$ if, for all smooth sections $Y$ of $\mathcal V$, we have that $[X,Y]$ is a smooth section of $\mathcal V$. 
If $F$ is a set of vector fields, we say that $\mathcal V$ is preserved by $F$ if $\VV$ is preserved by all $X \in F$.

%
%

\begin{prp}
\label{prp:inv_defined_on_quotient}
Let $G$ be a Lie group with $K \subset G$ being a connected closed subgroup. Consider $\pi : G \to G/K = \Omega$ the quotient map. Suppose that $\frk h \subset \frk g$ is a complex Lie subalgebra and that $\frk h$ is preserved by $\ker \pi_*$. Then, $\pi_* \frk h$ defines an invariant involutive structure on $\Omega$.
\end{prp}

\begin{proof}
 It follows from Problem 2.57 of \cite{gadea2012analysis} that $\pi_* \frk h$ is an involutive vector bundle. The invariance follows directly from the definition.
\end{proof}

The following definition is going to be useful.

\begin{dfn}
    Let $G$ be a Lie group and let $\frk k \subset \frk g_\R$ be a subalgebra. We say that $\frk k$ is a \emph{closed subalgebra} of $\frk g_\R$ if the group $K = \exp_G(\frk k)$ is a closed subgroup of $G$. When $\frk h$ is a left-invariant involutive structure, we say that $\frk h$ has \emph{closed real orbits} if $\frk k_\R = \frk h \cap \frk g_\R$ is a closed subalgebra.
\end{dfn}

The following lemma is one of the key technical results in this chapter.

\begin{lem}
    Let $G$ be a Lie group endowed with a left-invariant elliptic structure $\frk h$ having closed real orbits. Then, the homogeneous space $\Omega = G/K$, with $K = \exp(\frk h \cap \frk g_\R)$, has a natural complex structure given by $\pi_*(\frk h)$.
\end{lem}

\begin{proof}
    Since $\pi$ is a real map, we have
    $$\pi_*(\frk h) + \overline{ \pi_*(\frk h) } = \pi_*(\frk h) + \pi_*(\overline{ \frk h } ) =  \pi_*(\frk h + \overline{ \frk h } ) = \pi_*(\frk g) = \C T \Omega.$$
    The last equality holds because $\pi_*$ is surjective, implying that $\pi_*(\frk h)$ is elliptic. Notice that, by dimensional reasons, the involutive structure $\pi_*(\frk h)$ is a complex structure.
\end{proof}

\begin{thm}
\label{thm:homogeneous:converges_in_E_2}
 Let $G$ be a compact connected Lie group and let $\frk h \subset \frk g$ be an elliptic involutive structure having closed real orbits. Let $\frk k = \frk h \cap \overline{\frk h}$ and assume that there exists an ideal
$\frk u \subset \frk h$ such that $\frk k \oplus \frk u = \frk h$. Consider the homogeneous space $\Omega = G / K$, with $K = \exp(\frk k \cap \frk g_\R)$, endowed with the complex structure $\mathcal V = \pi_* \frk h$ induced by the quotient map $ \pi : G \to \Omega$. Then, we have
\begin{equation}
\label{eq:decomposition_in_E_2}
 H^{p,q} (G; \frk h) = \sum_{r + s = q} H^{p,r}(\Omega; \mathcal V) \otimes H^s(K; \C).
\end{equation}
\end{thm}

For the proof, we need the following technical lemma:

\begin{lem}
\label{lem:elliptic:important}
Let $G$, $\frk h$, $\frk k$, $K$, and $\frk u$ be as in Theorem \ref{thm:homogeneous:converges_in_E_2}. Let $u$ be a left-invariant closed $s$-form in $K$. Then $u$ can be extended to a $\dext'$-closed form in $G$. This extension, when restricted to a leaf $gK$, is a left-invariant form on that leaf $gK$.
\end{lem}

\begin{proof}
\label{lem:cohomology_extension}
Since $u$ is left-invariant, we can regard it as an element of the dual of $\Lambda^s\frk k$ and we can extend it to the dual of $\Lambda^s\frk g$ by defining it as zero if any of its arguments it in $\frk u$ or $\overline{\frk u}$. Let $w_g = (L_{g^{-1}})^* u$. Notice that $$(L_g)^* : H^s(gK) \to H^s(K)$$ is an isomorphism. Thus, $w$ is a smooth $s$-form in $G$ and the restriction of $w$ to a leaf $gK$ is representive of a cohomology class in that leaf $gK$.

We are going to see that $w$ is also $\dext'$-closed. We just need to show that $w$ satisfies $\dext w \in N^{1,r}(\frk g; \frk h)$. Let $X_1, \ldots, X_{r+1} \in \frk h$. Since the exterior derivative commutes with the pullback, we have
$$ \dext w_g = \dext (L_{g^{-1}})^* u = (L_{g^{-1}})^* u (\dext u).$$
Now, since every term $X_j u (X_1, \ldots, \hat {X_j}, \ldots, X_{r+1})$ is zero, we only have to do the following computation
$$
\begin{aligned}
\dext u  (X_1,  \ldots, X_{r+1})  & = \sum_{j=1}^{r+1} (-1)^{j+1} X_j u (X_1, \ldots, \hat {X_j}, \ldots, X_{r+1}) \\
    & + \sum_{j < k} (-1)^{j+k+1} u ([X_j, X_k], X_1, \ldots, \hat {X_j}, \ldots, \hat {X_k}, \ldots, X_{r+1}) \\
    & = \sum_{j < k} (-1)^{j+k+1} u([X_j, X_k], X_1, \ldots, \hat {X_j}, \ldots, \hat {X_k}, \ldots, X_{r+1}).
\end{aligned}
$$

If we assume that $X_j \in \frk k$ for all $j$, then $(\dext  w_g)  (X_1,  \ldots, X_{r+1}) = 0$ because $\dext u = 0$.

However, if we assume that $X_1 \in \frk u$ and that $X_j \in \frk k$ for j > 1, then
$$
\begin{aligned}
\sum_{j < k} (-1)^{j+k+1} & u([X_j, X_k], X_1, \ldots, \hat {X_j}, \ldots, \hat {X_k}, \ldots, X_{r+1}) \\ &= \sum_{1 < k} (-1)^{1+k+1} u([X_1, X_k], X_2, \ldots, \hat {X_j}, \ldots, \hat {X_k}, \ldots, X_{r+1}) = 0
\end{aligned}
$$
because $\frk u$ is an ideal and so $[X_1, X_k] \in \frk u$.

Finally, if we assume that two or more elements are $X_1,X_2 \in \frk k^\perp$, we can use the fact that $\frk u$ is an ideal and an argument similar to the one above to prove that we have $$(\dext  w_g)  (X_1(g),  \ldots, X_{r+1}(g)) = 0.$$ Therefore, $\dext  w_g \in N^{1,q}(G; \frk h)$.
\end{proof}

\begin{proof}[Proof of Theorem \ref{thm:homogeneous:converges_in_E_2}]



    Let $[\beta] \in H^s(K; \C)$. We can always choose $\beta$ left-invariant, so that, by the last Lemma \ref{lem:cohomology_extension}, there exists an extension of $\beta$ to $G$ such that $\dext' \beta' = 0$. This is a cohomology extension and, according to Leray-Hirsh Theorem \citep[Theorem 9 of Section 7, Chapter 5]{spanier1994algebraic}, the proof is completed.
\end{proof}

Notice that, in Equation \eqref{eq:decomposition_in_E_2}, the terms $H^s(K; \C)$ can always be computed by using algebraic methods. On the other hand, the terms $H^{p,r}(\Omega; \mathcal V)$ can be more complicated to compute.

Our next objective is to find some conditions for which the terms $H^{p,r}(\Omega; \mathcal V)$ can be computed using only algebraic methods. The first thing we are going to do is find conditions so that $\Omega$ is a compact connected Riemann surface. That is exactly what we do in next theorem:

\begin{thm}
\label{thm:omegaIsRieaman}
    Let $G$ be a compact Lie group endowed with an elliptic subalgebra $\frk h \subset \frk g$ such that $\dim \frk h = \dim \frk g - 1$ with $\frk k_\R = \frk h \cap \frk g_\R$ closed. Let $K$, $\Omega$ and $\mathcal V$ be defined as in Theorem \ref{thm:homogeneous:converges_in_E_2}. Then $$ H^{p,q} (G; \frk h) =  H^q(K; \C) + H^{q-1}(K; \C).$$
\end{thm}

\begin{proof}
    Let $\ip{,}$ be a Hermitian extension of any ad-invariant inner product on $\frk g_\R$.
		We write $\frk k = \frk h \cap \overline{\frk h}$ and we define $$\frk k^\perp = \{Z \in \frk h : \ip{Z,W} = 0 ~\forall W \in \frk h\}$$. Since $\dim \frk h = \dim \frk g - 1$, we have that $\dim \overline{\frk k^\perp} = 1$ and we clearly have $$ \frk g = \frk k \oplus \frk k^\perp \oplus \overline{\frk k^\perp}. $$
Since $\frk k^\perp$ has dimension 1, it is abelian. More than that, it is an ideal. In fact, we have that, if $T \in \frk k_\R$ and $Z \in \frk h^\perp$, then
		$$ \ip{[Z,T],U} = \ip{Z,[T,U]} = 0, \qquad \forall U \in \frk k, $$
		which means that $[Z,T] \in \frk k^\perp$, for all $T \in \frk k_\R$ and, by linearity, this is also true for all $T \in \frk k$. Therefore, if we assume that $K = \exp_G(\frk k_\R)$ is closed, we have all the necessary conditions to apply Theorem \ref{thm:homogeneous:converges_in_E_2}. Notice that, in this case, the homogeneous space $\Omega = G / K$ is a compact Riemann surface, which concludes the proof.
\end{proof}

Notice that, in Example \ref{exa:there_are_examples}, we constructed an example satisfying the hypothesis of the Theorem \ref{thm:omegaIsRieaman}.

Now we are going to discuss some techniques to deal with the case in which the complex dimension of $\Omega$ is bigger than 1. For this, we need to introduce the concepts of Lie algebra cohomology and of Lie algebra cohomology relative to a subalgebra. With this new concepts, we can state Bott's theorem, which gives us an algebraic way to compute the Dolbeault cohomology of certain compact homogeneous complex manifolds.


\section{Cohomology of Lie algebras}
\label{sec:cohomology_of_lie_algebras}

In this section, we define algebraic counterparts of the objects related to bidegree cohomology. These objects are the same introduced by Hochschild and Serre \cite{hochschild1953cohomology}, but here we used a more suitable notation. A good part of this section was inspired by Bott's paper \cite{bott1957homogeneous}, but the notation was changed to better suit the necessities of this work. Most of the results in this section were written for convenience of the reader, but, in order to not make the text too long, the proofs were omitted. We suggest the reader to consult the original papers for a detailed description of this theory and proofs.

The next two definitions are fundamental for this section. First, we define what is a representation of a Lie algebra on a vector space. Then, we use that representation to define what is an action of the Lie algebra's elements on the vector space. This notion is necessary to define what is the cohomology associated to Lie algebras.

\begin{dfn} Let $\frk g$ be a Lie algebra over $\C$ and $M$ a vector space over $\C$. A \emph{representation} of $\frk g$ is a Lie algebra homomorphism of $\frk g$ into the Lie algebra of all linear transformations of $M$ into $M$.
\end{dfn}

\begin{dfn}
Let $\frk g$ be a Lie algebra over $\C$ and $M$ a vector space over $\C$. We say that the pair $(M, \phi)$ is a \emph{$\frk g$-module} if $\phi$ is a representation of the Lie algebra $\frk g$ into $M$.
We denote by $X \cdot x$ the action of $\phi(X)$ on the element $x \in M$.
\end{dfn}

By definition, we have the following properties:

\begin{itemize}

\item If $X, Y \in \frk g$, $\lambda, \mu \in \C$, and $x \in M$, then
$$(\lambda X +  \mu Y) \cdot x = (\lambda X) \cdot x + (\mu Y) \cdot x.$$

\item If $X \in \frk g$, $\lambda, \mu \in \C$, and $x, y \in M$, then
$$X \cdot (\lambda x + \mu y)  = X \cdot (\lambda x) + X \cdot (\mu x).$$

\item If $X, Y \in \frk g$ and $x, y \in M$, then
$$[X,Y] \cdot x  = X \cdot ( Y \cdot x ) - Y \cdot (X \cdot x).$$

\end{itemize}

We denote by $M^{\frk g}$ the subspace of $M$ consisted of all $x \in M$ with $X \cdot x = 0$, for all $X \in \frk g$. We denote by $C^p(\frk g; M)$ the set of all alternating multi-linear $p$-forms on $\frk g$ with values in $M$. Also, we denote by $C(\frk g; M)$ the formal sum $\sum_p C^p(\frk g; M)$. We identify $C^0(\frk g; M)$ with $M$. Notice that each $C^p(\frk g; M)$ is a vector space.

If $\frk h$ is an ideal of $\frk g$ (here we allow $\frk h = \frk g$), we can turn $C^p(\frk h; M)$ into a $\frk g$-module. Indeed, since $C^0(\frk h; M) = M$, it already has a structure of a $\frk g$-module. For $p > 0$, $u \in C^p(\frk h; M)$, $X \in \frk g$, and $Y_1, \ldots, Y_p \in \frk h$, we define
$$(\mathcal L_Y u)(Y_1, \ldots, Y_p) = X \cdot u(Y_1, \ldots, Y_p) - \sum_{i=1}^n u(Y_1, \ldots, Y_{i-1}, [X,Y_i],Y_{i+1}, \ldots, Y_p).$$

\begin{lem} It holds that $\mathcal L_Y u \in C^p(\frk h; M)$ for all $u \in C^p(\frk h; M)$ and $\mathcal L$ depends linearly on $Y$ and on $u$. 
\end{lem}

If $ u \in C^{p+1}(\frk h; M) $ and $Y \in \frk h$, we define $\imath_Y u \in C^p(\frk h; M)$ by
$$\imath_Y u (Y_1, \ldots, Y_p) = u (Y, Y_1, \ldots, Y_p).$$

Notice that this is an algebraic definition of the Lie derivative we know from the basic theory of differential geometry. The Lie derivative has the following property: for $X, Y \in \frk g$ and $u \in C^p(\frk g; M)$, we have
$$\imath_X (\mathcal L_Y u) - \mathcal L_Y ( \imath_X u) = \imath_{[X,Y]} u.$$

The coboundary operator is the unique $\frk g$-homomorphism $$\dext : C^p(\frk g; M) \to C^{p+1}(\frk g; M)$$ such that $$\imath_X \circ \dext + \dext \circ \imath_X = \mathcal L_X,$$
for all $X \in \frk g.$

It is a basic exercise to verify that, for all $X \in \frk g$, it holds that $ \dext (\mathcal L_X u) = \mathcal L_X (\dext u)$. It is also easy to verify that, for $p > 0$, $u \in C^p(\frk g; M)$, and $X_1, \ldots, X_{p+1} \in \frk g$ we have the following formula
\begin{equation}
\label{eq:algebraic_diff_op}
\begin{aligned}
	\dext u(X_1, \ldots, X_{p+1})
		= & \sum_{j=1}^{p+1} (-1)^{j+1} X_j \cdot u(X_1, \ldots, \hat{X}_j, \ldots, X_{p+1})  \\
		 + & \sum_{j < k} (-1)^{j+k+1} u([X_j, X_k], X_1, \ldots, \hat X_j, \ldots, \hat X_k, \ldots, X_p).
\end{aligned}
\end{equation}

From the formula \ref{eq:algebraic_diff_op}, we see that the operator $\dext$ is just the algebraic version of the usual exterior differentiation on smooth manifolds. As expected, we have: if $u \in C^0(\frk g; M) = M$, then $(\dext u)(X) = X \cdot u$ and, obviously, it holds that $\dext^2 = 0$.

Therefore, we have a complex with respect to $\dext$. The cohomology space associated to $\dext$, denoted by $H^*(\frk g; M)$, is called cohomology module of $\frk g$ with coefficients in $M$. It is easy to see that $H^*(\frk g; M)$ is a vector space over $\C$.

Now we want to discuss some conditions on modules and on Lie algebras so that we can construct isomorphisms between cohomology spaces. Let $\lambda : \frk{g}' \to \frk{g}$ be a Lie homomorphism and $\eta : M \to M'$ be a linear map of the $\frk g$-module $M$ into the $\frk g'$-module $M'$. Also, let $\lambda^* \otimes \eta : C(\frk g; M) \to C(\frk g'; M')$ be the following linear map:
$$ [(\lambda^* \otimes \eta) u](X_1, \ldots, X_n) = \eta [u (\lambda X_1, \ldots, \lambda X_n)] \quad X_i \in \frk {g}',~u \in C^n(\frk g; M).$$

\begin{dfn} We say that the maps $\lambda$ and $\eta$ are \emph{compatibles} if
$$ \eta( \lambda X \cdot u) = X \cdot \eta u \qquad X \in \frk g',~u \in M.$$
\end{dfn}

With this notion of compatibility, we can state the following basic but important result:

\begin{lem} If $\lambda$ and $\eta$ are compatibles, the map $\lambda^* \otimes \eta $ commutes with the differential operators of the complexes in question, so it induces a homomorphism between $H^*(\frk g; M)$ and $H^*(\frk g'; M')$.
\end{lem}

Now we want to expand a little bit the notions we just introduced in order to be able to cover the cohomologies related to involutive structures. Let $R$ be a commutative algebra over $\C$. The Lie algebra of all derivations of $R$ is then itself a left $R$-module, which we denote by $\mathcal{D}$.

\begin{dfn}
If $M$ is a $R$-module, then $B_R(\mathcal{D}; M)$ denotes the graded $R$-module of $R$-linear alternating functions from $\mathcal{D}$ to $M$.
\end{dfn}

\begin{dfn}
Let $M$ be a $R$-module also having a $\mathcal{D}$-module structure. We say that the $R$-module structure is compatible with the $\mathcal{D}$-module structure if
$$ X \cdot (\lambda u) = (X \cdot \lambda)u + \lambda (X \cdot u) \qquad X \in \mathcal{D},~\lambda \in R,~u \in M.$$
\end{dfn}

Suppose that $M$ is a $R$-module having a compatible $\mathcal{D}$-module structure. In this situation, $C(\mathcal{D}; M)$ is well defined and $B_R(\mathcal{D}; M)$ can be identified with the subspace of $R$-linear maps in $C(\mathcal D; F).$

\begin{lem}
The subspace $B_R(\mathcal{D}; M)$ is closed under $\mathcal L_X,~ \imath_X$ and $\dext$.
\end{lem}

As a consequence, $B_R(\mathcal{D}; M)$ is a subcomplex of $C(\mathcal{D}; M)$. This subcomplex is denoted by $C_R(\mathcal{D}; M)$ and its cohomology module is denoted by $H^*_R(\mathcal D; M)$. If $\mathcal{D}_0$ is a subalgebra of $\mathcal D$, closed under $R$, then $H^*_R(\mathcal D_0; M)$ and $H^*_R(\mathcal D, \mathcal D_0; M)$ are defined analogously from the complexes $C_R(\mathcal D_0; M)$ and $C_R(\mathcal D, \mathcal D_0; M)$, respectively.

Now we are finally able to show how to use the algebraic language introduced in this section to describe one of the most important examples of involutive structures:

\begin{exa}
Let $\Omega$ be a smooth manifold and let $R = C^\infty(\Omega)$ be the algebra of all smooth complex valued functions on $\Omega$. The set of all derivations of $R$ is the Lie algebra of all smooth vector fields $\frk X (\Omega)$. Take $\mathcal D = \frk X (\Omega)$ and $M = C^\infty(\Omega)$. In this case, $C^q_R(\mathcal D; M)$ becomes the complex $\mathcal C^\infty(\Omega; \Lambda^q)$ of all complex valued differential forms on $\Omega$ and $H^*_R(\mathcal D; M)$ is the cohomology of $\Omega$ in the sense of de Rham. We usually write $H^q(\Omega; \C)$
for $H^q_R(\mathcal D; M)$.
\end{exa}

 If $f : \Omega' \to \Omega$ is a smooth map, then $f$ induces a homomorphism $$f^* : C_R(\mathcal D; M) \to C_{R'}(\mathcal D'; M')$$ which commutes with $\dext$ and so determines a homomorphism $$f^* : H^*(\Omega; \C) \to H^*(\Omega'; \C).$$

\begin{rmk} The map $f$ does not induce a Lie homomorphism from $\mathcal D'$ into $\mathcal D$, and this is the disadvantage of treating differential forms on $\Omega$ from this point of view.
\end{rmk}

The following example and theorem are fundamental for this text. They provided the main motivation for this project.

\begin{exa} Let $G$ be a Lie group. Also, consider $R = C^\infty(G; \C)$ and let $\frk g$ be the complex Lie algebra of $G$. In this case, the injection $$\rho : \frk g \to \mathcal D$$ extends to a $R$-linear Lie homomorphism
$$\rho : \frk g \otimes_\C R \to \mathcal D,$$
which is bijective because $\frk g$ is a $R$-base for $\mathcal D$. Hence, the induced map
$$ \rho^* : C_R(\mathcal D; R) \to C_R(\frk g; R) $$ is bijective and a chain homomorphism. We conclude that
$$ H_R^*(\mathcal D; R) \cong H^*(\frk g; R).$$

In this case, we also have $H^*(G; \C) \cong H^*(\frk g; R)$.
\end{exa}

In \cite[Theorem 2.3]{chevalley1948cohomology}, Chevalley and Eilenberg proved the following result:

\begin{thm} Let $G$ be a connected compact Lie group, then
$$ H^*(G) = H^*(\frk g; \C).$$
\end{thm}


\subsection{The cohomology induced by subalgebras}

Now we are going to introduce the concept of Lie algebra cohomology induced by a subalgebra. We reinforce here that this is the algebraic version of the cohomology induced by involutive subbundles of the tangent bundle. This is going to be really important for our results. Let $\frk h \subset \frk g$ be a subalgebra.

We define $N^{0,q}_{\frk h}(\frk g; M) = C^q(\frk g; M)$ and, for $p > 0$,
$$N^{p,q}_{\frk h}(\frk g; M) = \{ u \in C^{p+q}(\frk g; M): u(X_1, \ldots, X_{p+q}) = 0 \text{, when } q+1 \text{ arguments are in } \frk h \}.$$
If $q < 0$, we define $N^{p,q}_{\frk h}(\frk g; M) = \{0\}$.

Since $N^{p,q}_{\frk h}(\frk g; M) \subset N^{p+1,q-1}_{\frk h}(\frk g; M),$
we can define
$$ C^{p,q}_{\frk h}(\frk g; M) \doteq N^{p,q}_{\frk h}(\frk g; M) / N^{p+1,q-1}_{\frk h}(\frk g; M).$$
Also, since $\dext N^{p,q}_{\frk h}(\frk g; M) \subset N^{p,q+1}_{\frk h}(\frk g; M)$, we can define an operator
$$ \dext_{\frk h} : C^{p,q}_{\frk h}(\frk g; M) \to C^{p,q+1}_{\frk h}(\frk g; M)$$
and then we have a cochain complex.

For $p \geq 0$, we denote the set of the $(p,q)$-cocycles elements by $$ Z^{p,q}_{\frk h} (\frk g; M) = \ker \left( \dext : C^{p,q}_{\frk h}(\frk g; M) \to C^{p,q+1}_{\frk h}(\frk g; M) \right),$$
the set of $(p,q)$-coboundaries by
 $$ B^{p,q}_{\frk h}(\frk g; M) = \img \left( \dext : C^{p,q-1}_{\frk h}(\frk g; M) \to C^{p,q}_{\frk h}(\frk g; M) \right),$$
and the $(p,q)$-cohomology classes by
 $$ H^{p,q}_{\frk h}(\frk g; M) = \frac{Z^{p,q}_{\frk h}(\frk g; M)}{B^{p,q}_{\frk h}(\frk g; M)}.$$

Notice that, if $\frk h$ is a subalgebra of $\frk g$, then the space $C^p(\frk g / \frk h; M)$ has a natural structure of $\frk h$-module given by $$(X \cdot u)([X_1], \ldots, [X_p]) = X \cdot u([X_1], \ldots, [X_p]) - \sum_{j=1}^p  u([X_1], \ldots, [ [X, X_j] ], \ldots, [X_p]).$$

 With this structure, we have the following theorem whose proof can be founded in \cite{hochschild1953cohomology}:

\begin{thm}
The complexes $(C^{p,q}_{\frk h}(\frk g; M), \dext_{\frk h})$ and $(C^q(\frk h; C^p(\frk g / \frk h; M)), \dext)$ are isomorphic and $$ H^{p,q}_{\frk h}(\frk g; M) = H^{q}(\frk h; C^p(\frk g / \frk h ; M)). $$
\end{thm}

Notice that the left-invariant cohomologies we introduced in Chapter  \ref{chp:involutive_structures_on_compact_lie_groups} are equivalent to the Lie algebra cohomologies just defined. In fact, given any left-invariant form $u \in C^\infty_L(G; \Lambda^k )$, by restricting it to left-invariant vector fields, we get an element in $C^k(\frk g; \C)$. On the other hand, given any element $u \in C^k(\frk g; \C) $, by extending it $C^\infty(G)$-linearly, we get an element in $C^\infty_L(G; \Lambda^k)$.
This is easily seem to be an isomorphism between the cochains and therefore we have an isomorphism between the cohomology spaces, that is, $H^k(G) \cong H^k(\frk g)$.

Now consider $\frk h \subset \frk g$ a subalgebra and consider $\Lambda^{p,q}$ the vector bundle associated to $\frk h$ as we defined on Chapter \ref{chp:involutive_structures}.
Since the representatives are always left-invariant, we also can easily see that there is an one-to-one correspondence between elements of $C_L^\infty(G; \Lambda^{p,q})$ and elements of $C_{\frk h}^{p,q}(\frk g; C)$. This correspondence also gives us $H^{p,q}(G; \frk h) \cong H^{p,q}_{\frk h}(\frk g; \C)$.

\section{Bott's Theorem}

   In this section, we state a theorem by Bott that is useful when dealing with elliptic involutive structures on compact Lie groups. The following two lemmas are necessary for the statement of Bott's theorem:

\begin{lem}
    Let $G$ be a complex Lie group. If $f$ is (anti-)holomorphic, then $Xf$ is (anti-)holomorphic for all left-invariant vectors $X$.
\end{lem}

\begin{proof}
    First, we assume that $X$ is a real vector field. By definition, we have $$(Xf)(x) = \left. \frac{\dext}{\dext t} \right|_{t=0} f(x \exp(t X(e))) $$ and we write $f_t(x) = f(x \exp(t X(e)))$. Since our complex structure is left-invariant, we have that left translations are holomorphic and so $f_t$ is holomorphic for every $t$. Now, let $(z_1, \ldots, z_N)$ be local coordinates in $U$ and let $x \in U$. We have
    $$  \frac{\del}{\del \overline{z_j}} \frac{\dext}{\dext t} f_t(x) =  \frac{\dext}{\dext t} \frac{\del}{\del \overline{z_j}} f_t(x) = 0$$ and thus $Xf$ is holomorphic. For a complex left-invariant vector field $X$, we write $X = \Re X + i \Im X$ and thus $Xf$ is just the sum of two holomorphic functions, $\Re X f$ and $i\Im X f$.

    With the obvious adaptation of this argument, we have that, if $f$ is anti-holomorphic, then $Xf$ is anti-holomorphic.
\end{proof}

Let $U$ and $G$ be complex Lie groups with $U$ closed in $G$. Let $\Omega = G/U$ and suppose that $G$ is connected and $\Omega$ is compact and simply-connected.

By a theorem of Montgomery \citep{montgomery1950simply}, if $K$ is a maximal compact subgroup of $G$, under the above conditions, $K$ acts transitively on $\Omega$ and therefore $\Omega$ has another description:
$$ \Omega = K/H$$
with $H = U \cap K$.

We denote, respectively, by $\frk g, \frk u, \frk k, \frk h$ the complexified Lie algebras of $G$, $U$, $K$ and $H$. Since $G$ is a complex Lie group, we can decompose $\frk g$ into two ideals, that is,  $\frk g = \frk g_\alpha \oplus \frk g_\beta$, in which $\frk g_\alpha$ is the set of all left-invariant vector fields annihilated by all anti-holomorphic differential forms on $G$ and $\frk g_\beta$ is the set of all left-invariant vector fields annihilated by all holomorphic differential forms on $G$.

The sets $\frk g_\alpha$ and $\frk g_\beta$ are ideals. They are obviously subalgebras, so we only need to prove the following:

\begin{lem}
    If $X \in \frk g_\alpha$ and $Y \in \frk g_\beta$, then $[X,Y] = 0$.
\end{lem}

\begin{proof}
    Let $f$ be holomorphic. By definition, if $X \in \frk g_\alpha$, then $X$ annihilates every holomorphic function, so $Xf = 0$ and, by the preceding lemma, $Yf$ is also holomorphic. Therefore, $$ [X,Y]f = X(Yf) - Y(Xf) = 0$$ and we conclude that $[X,Y] \in \frk g_\alpha$. By an analogous argument, we have that $[X,Y] \in \frk g_\beta$ and thus $[X,Y] = 0$.
\end{proof}

Let $\alpha : \frk g \to \frk g_\alpha$ be the projection and denote by $\imath : \frk k \to \frk g$ the inclusion of $\frk k$ into $\frk g$. Write $\frk u_* = (\alpha \imath)^{-1}(\frk u)$. Clearly, $\frk u_* \subset \frk k$. In the case where $G$ is the complexification of $K$,
which always can be assumed by Montgomery's theorem and is the case we are going to be working on, $u_*$ can be identified with $\frk u_\R$, the real Lie algebra of $U$.

\begin{thm}[Bott's theorem]
Let $G$, $U$, $K$, $H$, and $\Omega$ with the conditions we just established and let $\mathscr O^p$ be the sheaf of local holomorphic $p$-forms on $\Omega$. Then,
$$ H^q(\Omega; \mathscr O^p) = H^q(\frk u_*, \frk h, \Lambda^p(\frk k / \frk u_*)^*)$$
with $u_*$ acting on $\Lambda^p(\frk k / \frk u_*)^*$ via adjoint action.
\end{thm}

\subsection{Applications of Bott's theorem}

Now we are going to prove some propositions that are going to be useful when applying Bott' theorem.


The next proposition shows that any compact homogeneous space $\Omega$ endowed with an invariant complex structure can be represented as a quotient of two complex Lie groups.

\begin{prp}
\label{prp:ext_act_complexified} Let $G$ be a connected compact Lie group acting transitively on a smooth manifold $\Omega$. If $\Omega$ is endowed with a complex structure invariant by the action of $G$, then this action extends to a transitive holomorphic action of $G_\C.$
\end{prp}

\begin{proof} We are assuming that the action of $G$ preserves the complex structure, so each automorphism $T_g : \Omega \to \Omega$ is holomorphic and we have a group homomorphism
\begin{equation}
\label{eq:cont_homo_of_groups}
	g \in G \mapsto T_g \in \Aut_\Ho(\Omega).
\end{equation}

Since $\Omega$ is compact, by \cite[Theorem 1.1, Chapter III]{kobayashi2012transformation}, the set $\Aut_\Ho(\Omega)$ is a complex Lie group. The topology on $\Aut_\Ho(\Omega)$ is the compact open topology and, in this case, it is the topology of uniform convergence over compact sets. Since $T : G \times \Omega \to \Omega$ is smooth, particularly it is uniformly continuous. Therefore, the map \eqref{eq:cont_homo_of_groups} is continuous and is automatically a Lie group homomorphism.

Let $(G_\C, \eta)$ be the universal complexification of $G$. By the universal property, there exists a complex Lie group homomorphism $g \in G_\C \mapsto T'_g \in \Aut_\Ho(\Omega)$ such that $T'_{\eta(g)} = T_g$, for all $g \in G$. Since $G$ is compact, we have that $\eta$ is injective. Therefore, $G$ can be identified with $\eta(G)$ and the action $T'$ can be considered an extension of $T$.
\end{proof}

The next proposition shows that, in order to study the cohomology of a left-invariant elliptic involutive structure defined over a compact semisimple Lie group, it is enough to study the cohomology of the structure lifted to the universal covering. The advantage of such approach is to be able to assume that the group in question is simply-connected, removing some topological barries.

\begin{prp}
\label{ss:prp:sctoss} Let $G$ be a semisimple compact Lie group and suppose that it is endowed with an elliptic involutive structure $\frk h \subset \frk g$. Then, its universal covering group $G'$ is also compact and admits an elliptic involutive structure $\frk h'$ such that $H^{p,q}(G; \frk h) \cong H^{p,q}(G';\frk h')$.
\end{prp}

\begin{proof}
By Weyl's theorem \citep[Theorem 4.26]{knapp2016representation}), the universal covering $G'$ is compact. We know that the covering map $\pi : G' \to G$ gives an isomorphism between $\frk g$ and $\frk g'$. Via this isomorphism, we define $\frk h' \subset \frk g'$. We also have that $\pi^{-1}(x)$ is finite for every $x \in G$. Therefore, by \cite[Theorem 11.1]{bredon2012sheaf}, we have the required isomorphism and the proof is completed.

In order to give a direct proof, we can construct the required isomorphism as follows. We denote by $\mathcal S_{\frk h}$ the sheaf of local solutions of $\frk h$ and by $\mathcal S_{ \frk h'}$ the sheaf of local solutions of $\frk h'$. Let $\mathcal U'$ be a finite covering of $G'$ consisting of sets satisfying the following conditions: for every $U' \in G'$ and $U = \pi(U')$, it holds that $\pi|_{U'}: U' \to U$ is a diffeomorphism and also, for $q > 0$, it holds that $H^q(U'; \mathcal S_{\frk h'}) = 0$.
The condition about the diffeomorphism is possible because $\pi$ is a covering map and the condition about the sheaf cohomology is possible because $\frk h'$ is elliptic. By taking $\mathcal U = \{ U = \pi(U'): U' \in U\}$, we have a finite open covering for $G$ such that $H^q(U; \mathcal S_{\frk h}) = 0$.

Next we construct a cochain isomorphism between $(C^q(\mathcal U', \mathcal S_{\frk h'}), \delta'^q)$ and $(C^q(\mathcal U, \mathcal S_{\frk h}), \delta^q)$. Since $\mathcal U'$ is finite, we can enumerate its elements $U'_1, U'_2, \ldots$ and we have a one-to-one correspondence with elements of $\mathcal U$, namely $U_1, U_2, \ldots$. Thus, for any simplex of $\sigma \in N(\mathcal U)$, there is an unique associated simplex $\sigma' \in N(\mathcal U')$, that is, if
$\sigma = (i_0, \ldots, i_q)$ and $|\sigma| = U_{i_0} \cap \cdots \cap U_{i_q} $, then $\sigma' = (i_0, \ldots, i_q)$ and $|\sigma'| = U'_{i_0} \cap \cdots \cap U'_{i_q} $.

Let $f' \in (C^q(\mathcal U', \mathcal S_{\frk h'}), \delta'^q)$ and let $\sigma'$ be a $q$-simplex. We define $f_\sigma$ in $|\sigma|$ by $f_\sigma(x) = f'(\pi^{-1}(x))$. The function $f_\sigma$ is well defined because $\pi|_{\sigma'}$ is a diffeomorphism and is in $\mathcal S_{\frk h}$.
We define a map $\phi : C^q(\mathcal U', \mathcal S_{\frk h'}) \to C^q(\mathcal U', \mathcal S_{\frk h'})$ which is obviously bijective and also commutes with the restrictions. Thus, we have that $\delta^q \phi(f') = \phi( \delta'^q f')$ and the map $\phi$ induces an isomorphism between $H^q(\mathcal U',\mathcal S_{\frk h})$ and $H^q(\mathcal U;\mathcal S_{\frk h})$.

Now, by using Leray's Theorem \cite[Section D, Theorem 4]{gunning1965analytic}, for every $q \geq 0$, we have that $H^q(G',\mathcal S_{\frk h'}) \cong H^q(\mathcal U';\mathcal S_{\frk h'})$ and $H^q(G,\mathcal S_{\frk h}) \cong H^q(\mathcal U;\mathcal S_{\frk h})$. Finally, we have $H^q(G',\mathcal S_{\frk h'}) \cong H^q(G,\mathcal S_{\frk h})$. 
\end{proof}

Now, we assume $G$ to be a compact semisimple Lie group. Combining Proposition \ref{ss:prp:sctoss}, Proposition \ref{prp:ext_act_complexified} and the long exact sequence of homotopy groups, we see that we can apply Bott's theorem to $\Omega$, which, in connection to Theorem \ref{thm:homogeneous:converges_in_E_2}, gives us a complete algebraic description of $H^{p,q}(G; \frk h)$.
To be explicit, we have the following theorem:

\begin{thm}
\label{thm:homogeneous:homo_on_simply3} Let $G$ be a connected, semisimple and compact Lie group and suppose that $\frk h \subset \frk g$ is an elliptic involutive structure having closed real orbits. Let $\frk k = \frk h \cap \overline{\frk h}$ and assume that there exists an ideal $\frk z \subset \frk h$ such that $\frk k \oplus \frk z = \frk h$.
Consider the homogeneous space $\Omega = G / K$, with $K = \exp_G(\frk k \cap \frk g_\R)$, endowed with the complex structure $\pi_* \frk h$. With $\frk u_*$ defined as in Bott's theorem, we have
$$ H^{p,q}(G; \frk h) \cong \sum_{r + s = q} H^r(\frk u_*, \frk k, \Lambda^p(\frk g / \frk u_*)) \otimes H^s(\frk k).$$
\end{thm}


\begin{proof}
By Proposition \ref{ss:prp:sctoss}, we can assume $G$ to be simply-connected. Notice that $K$ is connected by construction, therefore, we use the long exact sequence of homotopy groups to conclude that $\Omega$ is simply-connected. Since $G$ is compact, we have that $\Omega$ is compact. By combining this with the result obtained in Proposition \ref{prp:ext_act_complexified}, we can apply Bott's theorem. The proof is completed.
\end{proof}

\begin{cor}
    Let $G$ be a semisimple compact Lie group and let $T \subset G$ be a maximal torus. Consider the following elliptic Lie algebra $ \frk h = \frk t \oplus \bigoplus_{\alpha \in \Delta_+} \frk g_\alpha$ with $\frk t$ being the complexification of the Lie algebra of the maximal torus $T$. Then we have $$ H^{0,q}(G) = H^{0,q}(T).$$
\end{cor}

\begin{proof}

    We define the homogeneous manifold $\Omega \doteq G / T$ and we denote the quotient map by $\pi : G \to \Omega$. Notice that we have precisely the hypothesis of Theorem \ref{thm:homogeneous:homo_on_simply3}. Therefore, we have
    $$ H^{0,q}(G) = \sum_{ r + s = q} H^{0,r}(T; \frk t) \otimes H^{0,s}_{\delbar}(\Omega) = H^{0,q}(T)$$ and the last equality follows from the fact that  $H^{p,q}_{\delbar}(\Omega) = 0$ if $p \neq q$ and $H^{0,0}_{\delbar}(\Omega) = \C$. 
\end{proof}

\subsection{Examples}

    Now we are going to discuss a few examples. They are not a direct application of the theorems of this Chapter. All the following examples, except for the last one, are compact Lie groups endowed with a left-invariant elliptic structure whose cohomology can be computed by adapting the techniques we developed in this chapter.

    Let $G$ be a semisimple compact Lie group and let $T \subset G$ be a maximal torus. Let $\Omega \doteq G / T$ be a homogeneous manifold with quotient map denoted by $\pi : G \to \Omega$.

\begin{exa}
    Consider the following elliptic Lie algebra $$ \frk h = \frk e \oplus \bigoplus_{\alpha \in \Delta_+} \frk g_\alpha$$ with $\frk e$ being an elliptic structure on the maximal torus.

    The quotient map induces a complex structure on $\Omega$, which is known to be Kähler.

    Let $[u] \in H^{0,q}(T; \frk e)$. We know that we can choose $u$ to be a left-invariant form, so we can extend it to a cohomology extension and, by the Leray-Hirsch Theorem, we have
    $$ H^{0,q}(G; \frk h) = \sum_{ r + s = q} H^{0,r}(T; \frk e) \otimes H^{0,s}(\Omega; \pi_*(\frk h)) = H^{0,q}(T; \frk e).$$ The last equality follows from the fact that $H^{p,q}(\Omega;  \pi_*(\frk h)) = 0$ if $p \neq q$ and $H^{0,0}_{\delbar}(\Omega) = \C$. 

    Notice that, if $\frk k = \frk h \cap \overline{\frk h}$ is closed, then there exists an abelian ideal $\frk c$ such that we can write $\frk h = \frk k \oplus \frk c \oplus \bigoplus_{\alpha \in \Delta_+} \frk g_\alpha$ with $\frk u = \frk e \oplus \bigoplus_{\alpha \in \Delta_+} \frk g_\alpha$ being an ideal in $\frk h$.
		In this case, we can also apply Theorem \ref{thm:homogeneous:homo_on_simply3}.
\end{exa}

\begin{exa}
    Now let $$ \frk h' = \frk e \oplus \bigoplus_{\alpha \in \Delta} \frk g_\alpha$$ with $\frk e$ being an elliptic structure on the maximal torus.

    The projection map induces the usual De Rham cohomology on $\Omega$.

    Let $[u] \in H^{0,q}(T; \frk e)$. Since we can choose $u$ to be a left-invariant form, we can extend it to define a cohomology extension and then we have
    $$ H^{0,q}(G; \frk h') = \sum_{ r + s = q} H^{0,r}(T; \frk e) \otimes H^{s}(\Omega; \frk \C T \Omega).$$

    On the torus, we proved that we can compute the cohomology using only algebraic methods. On the homogeneous manifold, we know, by the Chevalley Eilenberg Theorem, that it can also be computed using only algebraic methods. Thus, the cohomology on $G$ can also be computed using only algebraic methods.


\end{exa}

\begin{exa}
    Now let $$ \frk h'' = \frk e \oplus \frk f,$$ in which $\frk f + \overline{\frk f} = \bigoplus_{\alpha \in \Delta} \frk g_\alpha$ and $\frk e$ is an elliptic structure on the maximal torus.

    The quotient map induces an elliptic structure invariant by the action of $G$ on $\Omega$. We are going to denote this structure by $\mathcal V_{\frk f}$.
    Here we also have a cohomology extension and thus $$ H^{0,q}(G; \frk h'') = \sum_{ r + s = q} H^{0,r}(T; \frk e) \otimes H^{0,s}(\Omega; \pi_*(\frk h'')).$$

    On the torus, we proved that we can compute the cohomology using only left-invariant forms. Unfortunately, it is not known if the cohomology on the homogeneous manifold, namely, $$ H^{0,s}(\Omega; \pi_*(\frk h'')),$$
		can be computed using only algebraic methods.
\end{exa}


















\appendix



\backmatter \singlespacing   
\bibliographystyle{alpha}
\bibliography{bibliography}

\printindex

\end{document}